\documentclass[a4paper,12pt]{article}
\usepackage{amsmath,amsthm,amscd,amsfonts,amsxtra,amssymb}
\usepackage{graphicx}
\usepackage{psfrag}
\usepackage{marvosym}
\usepackage{pstricks}
\usepackage{hyperref}

\newcommand{\ZZ}{\mathbb{Z}}

\newcommand{\RR}{\mathbb{R}}

\newcommand{\MA}{\mathcal{A}}

\newcommand{\fs}[6]{\left[F^{#1 #2 #3}_{#4}\right]_{#5 #6}}

\renewcommand{\fs}[4]{%
F^{#1,#2,#3}_{#4}
}

\newcommand{\cG}{{\mathcal G}}

\newcommand{\qnum}[1]{\lfloor #1 \rfloor}
\newcommand{\id}{{\mathbf 1}}

\newcommand{\qnr}[2][q]{\lfloor #2 \rfloor_{#1}}

\newcommand{\qbin}[3][q]{\left\lfloor  \begin{array}{c}  #2\\#3
                      \end{array}\right\rfloor_{#1}}

\newcommand{\cg}[7]{%
\begin{bmatrix}
#1 & #3 & #5 \\
#2 & #4 & #6 \\
\end{bmatrix}_{#7}
}

\newcommand{\qcg}[6]{%
\cg{#1}{#2}{#3}{#4}{#5}{#6}{q}
}

\newcommand{\lpa}{L^{+}_{1}}
\newcommand{\lpb}{L^{+}_{2}}
\newcommand{\lma}{L^{-}_{1}}
\newcommand{\lmb}{L^{-}_{2}}
\newcommand{\lm}{L^{-}}
\newcommand{\lp}{L^{+}}

\begin{document}

\begin{center}
{\Large\bf  Clebsch-Gordan and $6j$-coefficients for rank two quantum groups}\\[1cm]

{\large Eddy Ardonne} \\[1mm]
{\it Nordita, Roslagstullsbacken 23, 106-91 Stockholm, Sweden}\\[4mm]
             
{\large Joost Slingerland} \\[1mm]
{\it Department of Mathematical Physics, National University of Ireland\\
Maynooth, Ireland
\\and\\
School of Theoretical Physics, Dublin Institute for Advanced Studies\\
10 Burlington Rd, Dublin 4, Ireland}\\[2mm]
 
\end{center}

\begin{abstract}{\small\noindent
We calculate ($q$-deformed) Clebsch-Gordan and $6j$-coefficients for rank two
quantum groups. We explain in detail how such
calculations are done, which should allow the reader to perform similar calculations
in other cases. Moreover, we tabulate the $q$-Clebsch-Gordan and $6j$-coefficients
explicitly, as well as some other topological data associated with theories corresponding
to rank-two quantum groups. Finally, we collect some useful
properties of the fusion rules of particular conformal field theories.}
\end{abstract}

\section{Introduction}

Clebsch-Gordan coefficients and $6j$-symbols have long played an important
role in representation theory and also in quantum mechanics, where they traditionally appear in the addition of angular momenta and of internal quantum numbers of particles, such as isospin and color. Recently, models with so called \emph{topological phases} have come under intense investigation in condensed matter physics and these provide another arena where explicit knowledge of generalized Clebsch-Gordan and $6j$-coefficients is of great importance. Part of the current excitement is generated by the prospect that the unusual properties of the excitations of these phases can be harnessed and put to use in fault tolerant `topological quantum computers'~\cite{Kitaev03,Nayak08}. In view of this, we believe it is useful to provide the explicit data needed to compute time evolutions in a number of topological phases, with a description of the methods employed, so as to assist the reader in calculating similar data for other systems.  

In the low energy limit, all physical information in a topological  phase of matter can be expressed in the framework of topological field theory. Practical expressions for any physical amplitude in a topological field theory can be obtained in terms of generalized (often $q$-deformed) $6j$-symbols, usually called $F$-symbols, and a further set of numbers called the $R$-symbols. These data can be obtained from the representation theory of quantum groups. For an important class of topological models, namely those associated with $su(2)_k$
Chern-Simons theory, there exists an explicit expression for the $F$-symbols, which
was derived from the associated quantum group $U_{q}(sl(2))$. However, for most other topological models, no general formulas are available. Nevertheless, whenever the quantum group that describes a topological model is known, the $F$-symbols can be calculated explicitly. We will explain in detail how such calculations can be done, focusing mostly on the quantum groups arising from $q$-deformations of simple Lie algebras, which describe the Chern-Simons gauge theories based on the corresponding simple Lie groups. We provide explicit tables of $q$-Clebsch-Gordan coefficients and $F$-symbols, as well as $R$-symbols for rank two quantum groups. 

Although we will mainly focus on multiplicity free cases (that is, cases in which there are no fusion multiplicities greater than $1$), we will include one example of a theory which does contain a nontrivial fusion multiplicity, namely the $su(3)_3/Z_3$ orbifold theory. We are aware of only one other instance in which the $F$-symbols were calculated for a theory with a fusion multiplicity, namely in \cite{hh07}, for a theory with three types of particles. In contrast to the example treated here, that theory does not allow for a consistent braiding, and hence does not have an associated $R$-matrix.

This paper is organized as follows. In section \ref{topmodels}, we will give an overview of
the structure of topological models in general. The 'topological data' specifying topological
models will be introduced, and a quick introduction on how these can be used to calculate
physical observables will be given. Section \ref{quantumgroups} introduces basic representation theory of quantum groups which will be used throughout the paper to calculate topological data.
In particular, the definition of the $q$-Clebsch-Gordan coefficients and $F$-symbols will be given.
In section \ref{appfusgen}, we collect useful data on anyon theories. In particular, we list those theories (based on affine Lie algebras), which do not have fusion multiplicities. It turns out that many of the multiplicity free theories at high rank have the same fusion rules as other, `simpler', or better known theories, such as theories at lower rank or orbifolds of the chiral boson. We give such identifications with simpler theories for all affine Lie algebras of types $C$, $D$, $E$, $F$ and $G$.

A rather detailed description of the calculation of $q$-Clebsch-Gordan coefficients will
be given in section \ref{qclebschgordan}, which should enable the readers to carry out
such calculations themselves. Using these methods, we will give a simple expression
for the $R$-symbols in section \ref{rsymbols}.
A summary of the topological data considered in this paper is provided in
section \ref{tqftdata}. Finally, in the last part of the main text (section~\ref{outlook}), we give some applications of these topological data.

About half of this paper consists of appendices, in which we tabulate, amongst other
things, the $q$-Clebsch-Gordan coefficients and the $F$ and $R$ symbols for various
theories. These include $su(3)_2$ (appendix~\ref{appsu3}), $su(3)_3/Z_3$ (appendix~\ref{appsu33z3}) and the non-simply laced cases $so(5)_1$ (appendix \ref{appso51}) and
$g_{2,1}$ (appendix~\ref{appg21}). In particular, this covers all the rank two algebras at
the lowest non-trivial (i.e.~non-abelian) level. We also cover one of the simplest theories
containing a fusion multiplicity, namely $su(3)_3/Z_3$. For completeness, we give explicit
expressions for the case $su(2)_k$ as well, in appendix~\ref{appsu2k}. The $q$-Clebsch-Gordan
coefficients and $F$-symbols have been collected in four mathematica notebooks, which are
available via
\href{http://arxiv.org/src/1004.5456/anc}{this link to the mathematica notebooks}.

\section{Topological models}
\label{topmodels}

Topological phases are phases of matter characterized by the property that their low energy sectors may be described in the language of topological field theory \cite{Turaev,Kitaev06,Preskill-lectures,mslectures}. 
In particular, these phases have a finite number of types of gapped low energy excitations which may exhibit nontrivial behavior under fusion and exchanges. In planar systems, the exchange behavior is governed by a nontrivial representation of the braid group and the excitations are called \emph{anyons}. The different types of anyons are said to have different \emph{topological charges}.  
The amplitude for low-energy processes involving anyons are calculated using diagrams that may be interpreted as spacetime diagrams for the processes. The type, or topological charge, of each anyonic particle is indicated as a label on its worldline in such a diagram. As particles fuse or split, this gives rise to vertices where three worldlines meet. As the particles move around and exchange positions, this induces braiding of their world lines. The amplitude associated with a diagram is invariant under continuous deformations and may be calculated by the application of certain moves, most notably the so called $F$-moves. Using the $F$-moves one may express a diagram as a linear combination of diagrams in which one four particle process in the history depicted by the diagram has been recoupled. The coefficients that appear in these linear combinations are called the $F$-symbols and they essentially determine the value of any observable quantity in the theory. Most of this paper is devoted to an explicit calculation and tabulation of the $F$-symbols for a number of theories.

Before we get started with the definition of the $F$-symbols, we first quickly introduce the concept of fusion rules. The fusion rules of a topological theory state how many times the particle type $c$ appears in the decomposition of the (fusion) product of two particles of type $a$ and $b$,
\begin{equation}
a \times b = \sum_c {n_{a,b}}^c c.
\end{equation}
Here, the fusion coefficients ${n_{a,b}}^c$ are non-negative integers. When the TQFT is described using a quantum group (see also section~\ref{quantumgroups}), the fusion coefficient ${n_{a,b}}^c$ is just the multiplicity of the quantum group representation labeled $c$ in the decomposition of a suitably defined  tensor product of the representations labeled $a$ and $b$. 

The fusion rules are required to be associative. 
One also requires that 
${n_{a,b}}^{c} = {n_{b,a}}^{c}$, and that ${n_{1,a}}^{a'} = \delta_{a,a'}$, where $1$ denotes the trivial
particle. Each particle $a$ has a unique anti-particle $\bar{a}$, for which ${n_{a,\bar{a}}}^{1}=1$
(note that $a$ can be its own anti-particle, as is the case for all representations of $su(2)_k$).

We now define our standard set of $F$-symbols as the coefficients appearing in the diagrammatic equation
%~(\ref{Fdef}).
\begin{equation}
\label{Fdef}
\psset{unit=.5mm,linewidth=.2mm,dimen=middle}
\begin{pspicture}(0,0)(135,35)
\psset{arrowsize=2.5pt 2}
\multips(0,30)(10,-10){3}{\psline{->}(0,0)(6,-6)\psline(5,-5)(10,-10)}
\psline{->}(10,30)(10,23)
\psline(10,25)(10,20)
\psline{->}(20,30)(20,17)
\psline(20,20)(20,10)
%\rput(4,22){$a$}
%\rput(13,27){$b$}
%\rput(23,27){$c$}
%\rput(24,2){$d$}
%\rput(14,12){$e$}
\rput[B](0,31){$a$}
\rput[B](10,31){$b$}
\rput[B](20,31){$c$}
\rput(32,0){$d$}
\rput(14,12){$e$}
\rput(5,0){
\rput(70,15){$=\sum_{f} (\fs{a}{b}{c}{d})_{e,f}$}
}
\rput(20,0){
\multips(95,0)(10,10){3}{\psline{-<}(0,0)(7,7)\psline(5,5)(10,10)}
\psline{->}(105,30)(105,17)
\psline(105,20)(105,10)
\psline{->}(115,30)(115,23)
\psline(115,25)(115,20)
%\rput(101,27){$a$}
%\rput(111,27){$b$}
%\rput(121,22){$c$}
%\rput(101,2){$d$}
%\rput(111,12){$f$}
\rput[B](105,31){$a$}
\rput[B](115,31){$b$}
\rput[B](125,31){$c$}
\rput(93,0){$d$}
\rput(113,12){$f$}
}
\end{pspicture}
\end{equation}
There is a number of remarks to be made before we continue.
First of all, when the TQFT is described using a quantum group, the three lines labeled $a,b$ and $c$ represent three irreducible quantum group representations that are tensored together in two different orders, leading to different intermediate representations with irreducible components $e$ and $f$. The line labeled $d$ indicates that an irreducible component of type $d$ is selected from the three-fold tensor product as the \lq\lq overall topological charge''. The $F$-symbols then give a map between these components of the three-fold tensor product (again see section~\ref{quantumgroups} for more detail).
   
Secondly, the $F$-symbols as defined here are strictly speaking not numbers, but matrices which map between two topological vector spaces characterized by the fusion trees in the diagrams. We will mainly restrict our attention to \emph{multiplicity free} TQFTs, for which all such fusion spaces are either one-dimensional or zero dimensional. In the second case, the diagram has amplitude zero (we say that it is not allowed by the fusion rules). In the first case, all the vertices correspond to
one-dimensional `fusion spaces', and the $F$-symbols can be viewed as numbers (they are $1\times 1$-matrices). We can then view $\left[F^{abc}_{d}\right]$ as a matrix, whose rank is given by the number of consistent choices of $e$ (or $f$, which is equivalent), such that the fusion rules are satisfied. 

Since the $F$-symbols are the coefficients of a transformation between two topological vector spaces determined by fusion trees, they are only fully determined after the bases for these spaces are chosen. A change of basis will transform the $F$-symbols and since every choice of basis is equally good this gives us a gauge freedom in the $F$-symbols. Since we will (mostly) deal with theories without fusion multiplicities and we will use only orthonormal bases, this gauge freedom reduces to a choice of a phase for every vertex in the two diagrams that appear in the definition of the $F$-symbols.

It is not too difficult to see that repeated application of recoupling moves such as that shown in eq.~(\ref{Fdef}) allows us to reduce any diagram that involves only fusions and splittings to a standard form, in effect fixing the time evolution of the state of the system (see for instance~\cite{Kitaev06} for more details). More generally, diagrams may have crossings of the particle lines, corresponding to exchanges of the particles. These may by removed by $R$-moves, such as the one shown below. 
\begin{equation}
\pspicture[shift=-0.65](-0.1,-0.2)(1.5,1.2)
  \small
  \psset{linewidth=0.9pt,linecolor=black,arrowscale=1.5,arrowinset=0.15}
  \psline{->}(0.7,0)(0.7,0.43)
%  \psset{linestyle=dashed}
  \psline(0.7,0)(0.7,0.5)
% \psarc(0.9,0.846410){0.4}{120}{240}
% \psarc(0.5,0.846410){0.4}{-60}{40}
 \psarc(0.8,0.6732051){0.2}{120}{240}
 \psarc(0.6,0.6732051){0.2}{-60}{35}
  \psline (0.6134,0.896410)(0.267,1.09641)
  \psline{->}(0.6134,0.896410)(0.35359,1.04641)
  \psline(0.7,0.846410) (1.1330,1.096410)	
  \psline{->}(0.7,0.846410)(1.04641,1.04641)
  \rput[bl]{0}(0.4,0){$c$}
  \rput[br]{0}(1.35,0.85){$a$}
  \rput[bl]{0}(0.05,0.85){$b$}
% \scriptsize
%  \rput[bl]{0}(0.82,0.35){$\mu$}
  \endpspicture
=
%\sum\limits_{\nu }\left[ 
R_{c}^{ab}
%\right] _{\mu \nu}
\pspicture[shift=-0.65](-0.1,-0.2)(1.5,1.2)
  \small
  \psset{linewidth=0.9pt,linecolor=black,arrowscale=1.5,arrowinset=0.15}
  \psline{->}(0.7,0)(0.7,0.45)
  \psline(0.7,0)(0.7,0.55)
  \psline(0.7,0.55) (0.25,1)
  \psline{->}(0.7,0.55)(0.3,0.95)
  \psline(0.7,0.55) (1.15,1)	
  \psline{->}(0.7,0.55)(1.1,0.95)
  \rput[bl]{0}(0.4,0){$c$}
  \rput[br]{0}(1.4,0.8){$a$}
  \rput[bl]{0}(0,0.8){$b$}
%\scriptsize
%  \rput[bl]{0}(0.82,0.37){$\nu$}
  \endpspicture
\end{equation}
The $R$-symbol $R_{c}^{ab}$ that appears in this equation is in principle a unitary matrix, but in a fusion multiplicity free theory, this reduces to a number (in fact, a phase, since $R_{c}^{ab}$ will be a unitary $1\times 1$ matrix in this case). The $R$-symbols, like the $F$-symbols, depend on the choice of bases in the topological vector spaces (an exception are the $R$-symbols $R_{b}^{aa}$ which are gauge invariant), but often this gauge freedom is exhausted once the $F$-symbols are fixed -- this will be the case for the theories examined in this paper.

Using both $R$-symbols and $F$-symbols, all diagrams that can occur may be reduced to standard form and hence a knowledge of these symbols completely fixes the gauge invariant physical observables of the theory. Some important gauge invariant quantities characterizing topological models are quantum dimensions, twist factors, Frobenius-Schur indicators and the central charge. Formulas for these quantities in terms of the $F$ and $R$ symbols are given in section~\ref{tqftdata}.

\section{Quantum groups, CG-coefficients and $6j$-symbols}
\label{quantumgroups}
Quantum groups are algebras which have a number of structures that make sure that it is possible to define a tensor product on their representations and associated Clebsch-Gordan coefficients and $6j$-symbols. This makes them useful in physics especially in the theory of integrable models and models of anyonic systems. In the quantum group description of anyonic systems the irreducible representations of the quantum group correspond to the topological superselection sectors of the anyon model. Decomposition of tensor products of representations corresponds to fusion and $6j$-symbols correspond to the $F$-symbols. There is also a structure called the universal $R$-matrix which provides for braiding, and in particular gives the values of the $R$-symbols. The structures of the quantum group are defined in such a way that the $F$-symbols and $R$-symbols obtained from the quantum group's representation theory have all the required properties to fulfill their role in the corresponding anyon model. In particular, fusion is associative and the $F$-symbols and $R$-symbols satisfy the pentagon and hexagon equations (see Appendix \ref{apppenthex} for the most general form).
For a detailed introduction to quantum groups the reader may for instance consult \cite{Kassel95}.
Here will not even review all of the structure of a quantum group. Instead, we will discuss just enough of the structure to be able to define Clebsch-Gordan coefficients, $6j$-symbols and $R$-symbols. We will explicitly give the relevant structures for the quantum groups we are interested in, the $q$-deformed universal enveloping algebras $U_{q}(\mathfrak{g})$ based on the semisimple Lie algebras.

Let $\mathfrak{g}$ be a semisimple Lie algebra and denote its simple roots by $\alpha_{i}$. To each simple root, we associate three generators $H_{i},\lp_{i}$ and $\lm_{i}$. These generate $U_{q}(\mathfrak{g})$ as an algebra, subject to the relations
\begin{align}
[H_{i},H_{j}] &= 0 & [H_{i},L^{\pm}_{j}] &= \pm A_{ij}L^{\pm}_{j} &
[\lp_{i},\lm_{j}]=\delta_{ij}\qnr[q_i]{H_{i}} 
\end{align}
%\begin{equation}
%[H_{i},H_{j}]=0 ~~~~ [H_{i},L^{\pm}_{j}]= \pm A_{ij}L^{\pm}_{j} ~~~~ 
%[\lp_{i},\lm_{j}]=\delta_{ij}\qnr[q_i]{H_{i}} 
%\end{equation}
and
\begin{equation}
\sum_{s=0}^{1-A_{ij}} (-1)^s
\mbox{\small $\qbin[q_i]{\textstyle 1-A_{ij}}{\textstyle s}$}
(L^{\pm}_{i})^{1-A_{ij}-s}L^{\pm}_{j}(L^{\pm}_{i})^{s}=0 ~~~~(\mathrm{for~}i\neq j)
\end{equation}
Here, $A_{ij}= \frac{2(\alpha_{i},\alpha_{j})}{(\alpha_{j},\alpha_{j})}$
are the elements of the Cartan matrix of $g$ and we have defined $q_i=q^{1/t_i}$, where $t_{i}=\frac{2}{(\alpha_{i},\alpha_{i})}$ are integers in the set $\{1,2,3\}$ such that the matrix with
$(i,j)$-element $t_i A_{ij}$ is symmetric. In addition, these are the elements of the inverse of the quadratic form matrix.  When $q=1$, these relations reduce to the relations for the Chevalley-Serre basis of the universal enveloping algebra $U(\mathfrak{g})$.
The $q_i$-number $\qnr[q_i]{n}$ is given by 
\[
\qnr[q_i]{n} = \frac{q_i^{n/2}-q_i^{-n/2}}{q_i^{1/2}-q_i^{-1/2}} = \sum_{m=1}^{n} q_i^{\frac{n+1}{2}-m}\]
and the $q_i$-binomials that appear are defined by 
\[ 
\qbin[q_i]{n}{m}= \frac{\qnr[q_i]{n}!}{\qnr[q_i]{m}!\qnr[q_i]{n-m}!},
\]
where, for $n\ge 1$ we introduced 
\[
\qnr[q_i]{n}! = \prod_{m=1}^{n} \qnr[q_i]{m}
\] 
and for $n=0$ we take $\qnr[q_i]{0}!=1$. When $t_i=1$, we will often drop the subscript $i$ from $q_i$, and the subscript $q_i$ from the $q$-numbers and factorials altogether, i.e. $\qnr[q_i]{n} = \qnr[q]{n} = \qnum{n} $.

Any quantum group $\MA$ has a \emph{coproduct}, usually denoted $\Delta$, which is a homomorphism of algebras from $\MA$ into $\MA\otimes \MA$, which is coassociative, that is
\begin{equation}
\label{coasseq}
(\Delta\otimes {\rm id})\Delta=({\rm id}\otimes\Delta)\Delta.
\end{equation} 
Given two representations $\pi^{1},\pi^{2}$ of $\MA$, one can define a
tensor product representation $\pi^1\otimes\pi^2$ by the formula
\begin{equation}
\label{hopftp}
\pi^{1}\otimes\pi^{2}:~ x \rightarrow
(\pi^{1}\otimes\pi^{2})(\Delta(x)). 
\end{equation} 
Since $\Delta$ is a homomorphism, this will indeed be a representation and since $\Delta$ is coassociative, this tensor product will be associative, that is, different orders of tensoring multiple representations lead to the same overall representation of $\MA$.

For $U_{q}(\mathfrak{g})$, the coproduct is given on the generators by the formulae
\begin{eqnarray}
\Delta(H_{i})&=& 1\otimes H_{i} + H_{i}\otimes 1 \nonumber\\
\Delta(L^{\pm}_{i})&=&L^{\pm}_{i}\otimes
q_{i}^{H_{i}/4}+q_{i}^{-H_{i}/4}\otimes L^{\pm}_{i}.
\end{eqnarray}
If $q$ is not a root of unity, the representation theory of $U_{q}(\mathfrak{g})$ is very similar to the representation theory of $\mathfrak{g}$. The irreducible representations of $U_{q}(\mathfrak{g})$ are labeled by dominant integral weights of $\mathfrak{g}$. The module of the representation labeled by the weight $\lambda$ has a basis of eigenstates of the $H_{i}$, so that each such state is itself labeled by a weight of $\mathfrak{g}$. The action of the generators $L^{\pm}_{i}$ on a vector of weight $\mu$ sends this state either to zero or to a state with weight $\mu\pm\alpha_{i}$. The weight of a state is in general not enough to fix the state up to a constant, as there will be cases where the common eigenspaces of the $H_i$ are higher dimensional. When we describe representations in detail in subsequent sections, we will have to make some arbitrary choices to fix a basis for these higher dimensional eigenspaces (see section~\ref{sec:multiplicities} for this). However when no such choices are required, the matrix elements of the generators $L^{\pm}_{i}$ can be written in a form which directly mirrors the undeformed case. For example, the action of $L^{\pm}$ for the states in the irreducible representations of $U_{q}(sl(2))$ is given by  
\begin{equation}
L^{\pm} | j, m \rangle = \sqrt{\qnr[q]{j\mp m} \qnr[q]{j\pm m+1}} | j , m \pm 1 \rangle \ .
\end{equation}
Here $j\in \frac{1}{2} \ZZ$, and $m\in\{-j,-j+1,\ldots j-1,j\}$ label the allowed $z$-components of the `q-spin', or equivalently, the eigenvalues of $H/2$ in the representation. This formula easily translates to the general multiplicity free case.

The difference between the quantum group and the classical algebra it is based on is
best visible for $q$ a root of unity.  
Namely, when $q$ is a root of unity, one finds that all these representations
remain well-defined, but many are no longer irreducible and in particular there are indecomposable representations. In the following, we will focus on the irreducible representations. When taking
a tensor product of two irreducible representations, it often happens that indecomposable
representations occur. To deal with this problem, one has to introduce a `truncated' tensor
product, in which the unwanted, indecomposable representations no longer appear. This
`truncated' tensor product precisely corresponds to the fusion product (as introduced in section
\ref{topmodels}), for physical systems described by the Chern-Simons theory based on $\mathfrak{g}$ at level $k$, where $k$ is related to $q$ through
\[
q=e^{\frac{2\pi i}{k+g}}.
\]
Here $g$ is the dual Coxeter number of $\mathfrak{g}$.

We can now introduce the $q$-Clebsch-Gordan coefficients for this (truncated)
tensor product. 
%
%The quantum group $U_{q}(\mathfrak{g})$
%Quantum groups as a way to construct topological models (F-symbols)
%Quickly review $U_q(su(N))$ (prob. not $U_q(G)$). 
%Def. 6j-symbols, CG-coefficients, say something about the amount of freedom in their definition (gauge). Note that everything is simplified for cases without fusion multiplicities and that we will be concentrating on this case. 
The $q$-Clebsch-Gordan coefficients express the tensor product representations in
terms of the direct sum of irreducible representations. To write explicit CG-coefficients, we need to introduce notation for the states in an irreducible $U_{q}(\mathfrak{g})$-module.
As noted before, each irreducible module is labeled by a weight of $\mathfrak{g}$ and has a basis of eigenstates of the $H_{i}$, which are themselves each labeled by a weight of $\mathfrak{g}$. If all weight spaces which occur in the module have dimension equal to $1$, we can denote the basis states for the module simply as $\lvert j , m\rangle$, where $j$ is the weight labeling the module and $m$ is the weight labeling the weight space in the module. One may now define the Clebsch-Gordan coefficients by
\begin{equation}
\label{qcg_def}
\lvert j , m\rangle = \sum_{m_1,m_2}
\qcg{j_1}{m_1}{j_2}{m_2}{j}{m} \lvert j_1,m_1\rangle \lvert j_2,m_2 \rangle
\end{equation}
In other words the $CG$-coefficients $\qcg{j_1}{m_1}{j_2}{m_2}{j}{m}$ are the coefficients of the decomposition of the state $\lvert j , m\rangle$ in the irreducible $j$-component in the tensor product representation $j_1\otimes j_2$ into the standard product basis of this tensor product representation.

Of course, the basis states that appear above are only fixed up to a phase by the labels $j$ and $m$, so their phases may still be chosen to give `nicer' CG-coefficients. For tensor products that are not multiplicity free, more labels will be needed to distinguish the various irreducible representations of the same weight that may occur in the tensor product. This is dealt with in some detail in section~\ref{sec:multiplicities}, but it does not change the structure of what follows here (we can imagine the extra labels to be implicit in the $m$-labels).

The $q$-Clebsch-Gordan coefficients are orthogonal for $q \in R$ and we can use analytic continuation to obtain an orthogonality relation for arbitrary $q$, which we checked to hold for the $q$-CG coefficients we calculated,
\begin{equation}
\label{qcg_ortho}
\sum_{m_1,m_2}
\qcg{j_1}{m_1}{j_2}{m_2}{j}{m} 
\qcg{j_1}{m_1}{j_2}{m_2}{j'}{m'} = \delta_{j,j'}\delta_{m,m'} \ .
\end{equation}
We will use the `inner product' which is defined, like the usual inner product on $\RR^n$, as the sum of the products of the coefficients of the two vectors, \emph{without complex conjugation}. This inner product is positive definite for all real $q\ge0$, but for other values of $q$ this is not necessarily true. With this inner product the formula above guarantees that states in different irreducible subrepresentations of a tensor product are orthogonal \emph{for all} $q$. 

We calculate the $F$-symbols, or 
$q-6j$ coefficients, by making use of the $q$-Clebsch-Gordan coefficients.
The graphical representation of these symbols in our current notation is
\begin{center}
\psset{unit=.5mm,linewidth=.2mm,dimen=middle}
\begin{pspicture}(0,0)(135,30)
\psset{arrowsize=2.5pt 2}
\multips(0,30)(10,-10){3}{\psline{->}(0,0)(6,-6)\psline(5,-5)(10,-10)}
\psline{->}(10,30)(10,23)
\psline(10,25)(10,20)
\psline{->}(20,30)(20,17)
\psline(20,20)(20,10)
%\rput(4,22){$a$}
%\rput(13,27){$b$}
%\rput(23,27){$c$}
%\rput(24,2){$d$}
%\rput(14,12){$e$}
\rput(-2,30){$j_1$}
\rput(14,30){$j_2$}
\rput(24,30){$j_3$}
\rput(32,0){$j$}
\rput(14,12){$j_{12}$}
\rput(5,0){
\rput[lB](35,15){$=\sum_{j_{23}} (\fs{j_1}{j_2}{j_3}{j})_{j_{12},j_{23}}$}
}
\rput(20,0){
\multips(95,0)(10,10){3}{\psline{-<}(0,0)(7,7)\psline(5,5)(10,10)}
\psline{->}(105,30)(105,17)
\psline(105,20)(105,10)
\psline{->}(115,30)(115,23)
\psline(115,25)(115,20)
%\rput(101,27){$a$}
%\rput(111,27){$b$}
%\rput(121,22){$c$}
%\rput(101,2){$d$}
%\rput(111,12){$f$}
\rput(101,30){$j_1$}
\rput(111,30){$j_2$}
\rput(128,30){$j_3$}
\rput(93,0){$j$}
\rput(113,12){$j_{23}$}
}
\end{pspicture}
\end{center}
Thus, at real positive $q$, each F-symbol can be obtained as the inner product between two states obtained from the two different ways of fusing the three particles (or representations) $j_1$, $j_2$ and $j_3$.
These states themselves can be written in terms of the $q$-Clebsch-Gordan coefficients. Hence we obtain the formula 
\begin{equation}
\begin{split}
&(\fs{j_1}{j_2}{j_3}{j})_{j_{12},j_{23}} = \\
&\sum_{m_1,m_2,m_3,m_{12},m_{23}}
\qcg{j_1}{m_1}{j_2}{m_2}{j_{12}}{m_{12}}
\qcg{j_{12}}{m_{12}}{j_3}{m_3}{j}{j}
\qcg{j_2}{m_2}{j_3}{m_3}{j_{23}}{m_{23}}
\qcg{j_1}{m_1}{j_{23}}{m_{23}}{j}{j}
\end{split}
\end{equation}
The inner product (without complex conjugation!) used to derive this formula is really only well defined (positive definite) for $q>0$ real, but by analytic continuation, the formula nevertheless remains true for arbitrary $q$.

It would be natural to continue with a description of the $R$-symbols in terms of quantum group data. However, these are most easily expressed in terms of particular $q$-Clebsch-Gordan
coefficients. Therefore, we will wait with the $R$-symbols until section \ref{rsymbols},
which immediately follows the section in which we give a detailed explanation of the
calculation of the $q$-Clebsch-Gordan coefficients.

%As an example, we find that the number of $6j$-coefficients for $su(3)_2$
%is 405. Luckily, these are not all independent. It turns out that, up to signs,
%there are only 21 independent coefficients for arbitrary $q$, while for
%$q=e^{\frac{2 \pi i}{5}}$, there are only three possible values, again
%up to signs. Moreover, in the gauge we are using,
%an F-symbol with an incoming vacuum representation is always equal
%to one, and this is also true when the outgoing particle is the identity.

\section{On fusion multiplicities}
\label{appfusgen}

While working on this manuscript, we often wondered about properties of
the fusion rules of particular CFTs. Which theories have fusion multiplicities,
and what is the structure of those that don't? Although this is an issue which is
slightly off topic, we do take this paper as an opportunity to gather this
undoubtedly known information here. We should note that in gathering
this information, we benefited much from the program {\em Kac}, by
N. Schellekens \cite{kac}.

The table \ref{fusmult}, we list the theories, based on affine Lie algebra's,
or WZW CFT's, which do have fusion multiplicities.

\begin{table}[ht]
\begin{center}
\begin{tabular}{|rcl|c||c|c|}
\hline
$A_{r,k}$ & $\equiv$ & $su(r+1)_k$ & $k\geq 3\wedge r\geq 2$ & $E_{6,k}$ & $k\geq 3$ \\
$B_{r\geq 2,k}$ & $\equiv$ & $so(2r+1)_k$ & $k\geq 3$ & $E_{7,k}$ & $k\geq 3$ \\
$C_{r\geq 3,k}$ & $\equiv$ & $sp(2r)_k$ & $k\geq 2$ & $E_{8,k}$ & $k\geq 4$ \\
$D_{r\geq 4,k}$ & $\equiv$ & $so(2r)_k$ & $k\geq 3$ & $F_{4,k}$ & $k\geq 3$ \\
&& &  & $G_{2,k}$ & $k\geq 3$ \\
\hline
\end{tabular}
\end{center}
\caption{Theories with fusion multiplicities}
\label{fusmult}
\end{table}

\begin{table}[ht]
\begin{center}
\begin{tabular}{| c | c | c | c |}
\hline
& $k=1$ & $k=2$ & $k=3$ \\
\hline
$A_{r\geq 2}$ & $\mathbf{Z}_{r+1}$ & $su(r+1)_2$ & \\ 
$B_{r}$ & $su(2)_2$ & $so(2r+1)_2$ & \\
$C_{r}$ & $su(2)_r$ & & \\
$D_{r}$ &
$
\left\{
\begin{matrix}
\mathbf{Z}_2 \times \mathbf{Z}_2 & \mbox{($r$ even)}\\
\mathbf{Z}_4 \qquad & \mbox{($r$ odd)}\\
\end{matrix}
\right.
$ & $Z_2$ orbifold $R=2r$ & \\
$E_6$ & $\mathbf{Z}_3$ & $so(3)_5 \times \mathbf{Z}_3$ &\\
$E_7$ & $\mathbf{Z}_2$ & $so(3)_3 \times su(2)_2$  &\\
$E_8$ & $\mathbf{Z}_1$ & $su(2)_2$ & $so(3)_9$ \\
$F_4$ & $so(3)_3$ & $so(3)_9$ &\\
$G_2$ & $so(3)_3$ & $so(3)_7$ &\\
\hline
\end{tabular}
\end{center}
\caption{Identification of the fusion rules without fusion multiplicities, in terms of perhaps `better' known
fusion rings. The $k=2$ entries for $A_{r\geq 2}$ and $B_r$ are a tautology.}
\label{fusident}
\end{table}
In table \ref{fusident}, we give the list of WZW theories, without multiplicities, for which we
can identify the fusion rules as a `known' fusion ring. Two of the entries are a
tautology, namely $A_{r\geq 2,2}$ and $B_{r,2}$. For completeness, we note that of course,
the theories $A_{1,k}\equiv su(2)_k$, do not have fusion multiplicities for arbitrary $k$.
Again, we used {\em Kac} to gather
this information. For example, the fusion rules of $D_{r,k=2}$ are
the same as the fusion rules of the $Z_2$ orbifold of the chiral boson at radius $R=2r$,
which we checked explicitly up to $r=20$.
We note that with $so(3)_k$, we denote the integer spin sector of the
$su(2)_k$ theory. In particular, $so(3)_3$ corresponds to the Fibonacci theory (sometimes
denoted as ${\rm Fib}$), which consists of two anyon types ${\bf 1}$ and $\tau$,
with the non-trivial fusion rule $\tau \times \tau = {\bf 1} + \tau$.

\section{Calculating the $q$-Clebsch-Gordan coefficients}
\label{qclebschgordan}

In this section, we will show how to calculate the $q$-Clebsch-Gordan coefficients
in full detail. This quantum group calculation closely follows the calculation of
Clebsch-Gordan coefficients in the `classical' case, see for instance \cite{bl81}. 
We will use the case of $su(3)$ as an example, and hence, we will
be dealing with the quantum group $U_q (su(3))$. For details on the representation
theory of the finite dimensional Lie algebra $su(3)$, we refer to appendix \ref{appsu3},
and for instance the books \cite{fms,fs,cftqgconnection}. 

We will start with some straightforward examples, which will
explain the structure of the calculation, without having to worry about additional
complications. The complications consist of two type of multiplicities, namely
weight-space multiplicities and fusion-multiplicities, which will be dealt with after
we have completed explaining the structure of the calculations. 

We first introduce some notation needed in this discussion.
The representations of the rank two algebra $su(3)$ are labeled by
the Dynkin labels of the highest weight, $(\Lambda_1,\Lambda_2)$. In addition, we will
also denote the representations by their dimensions in boldface (and an over-line, to denote
the conjugate representations). The weights of the
representations are labeled by $(\lambda_1,\lambda_2)$. For example, the eight-dimensional
adjoint representation ${\bf 8}$ has highest weight $(\Lambda_1,\Lambda_2)=(1,1)$.
The other weights in this representation are
$(-1,2)$, $(2,-1)$, $(0,0)_+$, $(0,0)_-$, $(-2,1)$, $(1,-2)$ and $(-1,-1)$, which we include
here, in order to show the notation of the two-dimensional weight space $(0,0)$. We note that for
$q=e^{2\pi i/5}$ (or $k=2$), the ${\bf 8}$ is the only representation with a weight-space multiplicity.
In addition, the following representations ${\bf 1} = (0,0)$, ${\bf 3}=(1,0)$, $\bar{{\bf 3}}=(0,1)$,
${\bf 6}=(2,0)$, $\bar{{\bf 6}}=(0,2)$ and ${\bf 8} = (1,1)$ are present at $k=2$. 

As stated in section \ref{quantumgroups}, the $q$-Clebsch-Gordan coefficients express the tensor product representations in terms of the direct sum of irreducible representations 
\begin{equation}
\lvert j , m\rangle = \sum_{m_1,m_2}
\qcg{j_1}{m_1}{j_2}{m_2}{j}{m} \lvert j_1,m_1\rangle \lvert j_2,m_2 \rangle
\end{equation}

\subsection{The structure of calculating $q$-Clebsch-Gordan coefficients}

We will explain the calculation by taking a simple concrete example, based
on the $su(3)$ tensor product ${\bf 3} \otimes {\bf 3} = \bar{\bf 3} \oplus {\bf 6}$, which for
$k\geq 2$ corresponds to the fusion product as well: ${\bf 3} \times {\bf 3} = \bar{\bf 3} + {\bf 6}$. 
Will will start with the Clebsch-Gordan coefficients for the representation ${\bf 6}$ in the tensor product
decomposition of  ${\bf 3} \otimes {\bf 3}$.

The highest weight of the representation ${\bf 6}$ is $(2,0)$, which is the sum of the highest weights
of the two representations ${\bf 3}$, namely $(1,0)$. Thus, the highest weight of the ${\bf 6}$
uniquely decomposes, and we find the first (trivial) Clebsch-Gordan coefficient:
\begin{equation}
\label{td336}
|(2,0)(2,0)\rangle = |(1,0)(1,0)\rangle \otimes |(1,0)(1,0)\rangle \ ,
\end{equation}
namely
\begin{equation}
\qcg{(1,0)}{(1,0)}{(1,0)}{(1,0)}{(2,0)}{(2,0)} = 1 \ .
\end{equation}
This coefficient has norm one, and we choose to set the phase to zero by convention.
By acting with the lowering operator $\lma$ on both sides of eq.~\eqref{td336}, we find
other Clebsch-Gordan coefficients. On the right hand side, we actually need to use
$\Delta(\lma)$, which is given by $\Delta(\lma)= \lma \otimes q^{H_1/4} + q^{-H_1/4}\otimes \lma$.
This leads to the following
\begin{equation}
\label{sixex}
\begin{split}
|(2,0)(0,1)\rangle = \frac{\lma}{\sqrt{\qnum{2}}} |(2,0)(2,0)\rangle =
\frac{\Delta(\lma)}{\sqrt{\qnum{2}}}  |(1,0)(1,0)\rangle \otimes |(1,0)(1,0)\rangle =  \\
\frac{q^{1/4}}{\sqrt{\qnum{2}}} | (1,0)(-1,1) \rangle \otimes |(1,0)(1,0)\rangle +
\frac{q^{-1/4}}{\sqrt{\qnum{2}}} | (1,0)(1,0) \rangle \otimes |(1,0)(-1,1)\rangle
\end{split}
\end{equation}
We put in the factor $1/\sqrt{\qnum{2}}$ `by hand', to make sure that the states are normalized
throughout the calculation. Of course, one could also ignore these factors at first, and normalize
the states later on. To be able to discuss the normalization, 
we recall that we defined the inner product in the usual way for $q\in \RR$, that is, by
simple multiplication of the coefficients, and for arbitrary $q$ by analytic continuation of this definition. We note that this inner product does not involve complex conjugation of the coefficients. Recalling in addition that we used the following convention for the $q$-numbers
\begin{equation}
\qnum{n} = \frac{q^{\frac{n}{2}}-q^{-\frac{n}{2}}}{q^{\frac{1}{2}}-q^{-\frac{1}{2}}} \ ,
\end{equation}
it easily follows that the right hand side of equation~\eqref{sixex} is indeed normalized.

So, we found the additional $q$-Clebsch-Gordan coefficients
\begin{align}
\qcg{(1,0)}{(1,0)}{(1,0)}{(-1,1)}{(2,0)}{(0,1)} &= \frac{q^{-1/4}}{\sqrt{\qnum{2}}} \\
\qcg{(1,0)}{(-1,1)}{(1,0)}{(1,0)}{(2,0)}{(0,1)} &= \frac{q^{1/4}}{\sqrt{\qnum{2}}} 
\end{align}

We can now continue this procedure by acting with additional lowering operators.
In fact, in this case, we can either act with $\lma$ to obtain the state $|(2,0)(-2,2)\rangle$,
or with $\lmb$ to find the coefficients for $|(2,0)(1,-1)\rangle$. To see with which lowering
operators one can act, we refer to the appendix \ref{appsu3} containing the structure of
the representations we consider in this paper.
All the $q$-Clebsch-Gordan coefficients for the ${\bf 6}$ in ${\bf 3} \times {\bf 3}$ are
collected in appendix \ref{su3example}.

We will now explain how to obtain the $q$-Clebsch-Gordan coefficients for the representation
$\bar{\bf 3}$ in the tensor product decomposition of ${\bf 3}\otimes {\bf 3}$. To start with, we need
to express the (highest weight) state $|(0,1)(0,1)\rangle$ in terms of the states
$|(1,0)(1,0)\rangle \otimes |(1,0)(-1,1)\rangle$ and $|(1,0)(-1,1)\rangle \otimes |(1,0)(1,0)\rangle$.
We will do this by using that $|(0,1)(0,1)\rangle$ is a highest weight state, i.e
$\lpa |(0,1)(0,1)\rangle = \lpb |(0,1)(0,1)\rangle = 0$. In this case, we only need to consider
\begin{equation}
\begin{split}
\lpa |(0,1)(0,1)\rangle =& \Delta(\lpa)
\bigl(
a |(1,0)(1,0)\rangle \otimes |(1,0)(-1,1)\rangle + b |(1,0)(-1,1)\rangle \otimes |(1,0)(1,0)\rangle
\bigr) \\
 = &(a q^{-1/4} + b q^{1/4} ) |(1,0)(1,0)\rangle \otimes |(1,0)(1,0) \rangle = 0\ ,
%0 =\frac{\lpa}{\sqrt{\qnum{2}}} |(0,1)(0,1)\rangle =& \frac{\Delta(\lpa)}{\sqrt{\qnum{2}}}
%\bigl(
%a |(1,0)(1,0)\rangle \otimes |(1,0)(-1,1)\rangle + b |(1,0)(-1,1)\rangle \otimes |(1,0)(1,0)\rangle
%\bigr) \\
%& \frac{1}{\sqrt{\qnum{2}}} (a q^{-1/4} + b q^{1/4} ) |(1,0)(1,0)\rangle \otimes |(1,0)(1,0) \rangle \ ,
\end{split}
\end{equation}
from which it follows we should write
\begin{equation}
|(0,1)(0,1)\rangle =
\frac{q^{1/4}}{\sqrt{\qnum{2}}}  |(1,0)(1,0)\rangle \otimes |(1,0)(-1,1)\rangle 
- \frac{q^{-1/4}}{\sqrt{\qnum{2}}}  |(1,0)(-1,1)\rangle \otimes |(1,0)(1,0)\rangle \ ,
\end{equation}
where the overall sign (or better, phase) is our convention. In general, we choose the
overall phase to be such that for $q=1$, the Clebsch-Gordan coefficient with the highest possible
weight on the left hand side of the tensor product is real and positive
(we will spell out all our conventions in detail below).
We should note that this state $|(0,1)(0,1)\rangle$ is orthogonal to the state
$|(2,0)(0,1)\rangle$ (eq.~\eqref{sixex}), as it should. Thus, we find the following $q$-Clebsch-Gordan
coefficients
\begin{align}
\qcg{(1,0)}{(1,0)}{(1,0)}{(-1,1)}{(0,1)}{(0,1)} &= \frac{q^{1/4}}{\sqrt{\qnum{2}}} \\
\qcg{(1,0)}{(-1,1)}{(1,0)}{(1,0)}{(0,1)}{(0,1)} &= -\frac{q^{-1/4}}{\sqrt{\qnum{2}}} 
\end{align}

We can find the other coefficients related to the states $|(0,1)(1,-1)\rangle$ and $|(0,1)(-1,0)\rangle$
by first acting with the lowering operator $\lmb$ and subsequently with $\lma$ on the state
$|(0,1)(0,1)\rangle$. Without spelling out all the details, this gives the following $q$-Clebsch-Gordan
coefficients
\begin{align}
\qcg{(1,0)}{(1,0)}{(1,0)}{(0,-1)}{(0,1)}{(1,-1)} &= \frac{q^{1/4}}{\sqrt{\qnum{2}}} &
\qcg{(1,0)}{(0,-1)}{(1,0)}{(1,0)}{(0,1)}{(1,-1)} &= -\frac{q^{-1/4}}{\sqrt{\qnum{2}}} \\
\qcg{(1,0)}{(-1,1)}{(1,0)}{(0,-1)}{(0,1)}{(-1,0)} &= \frac{q^{1/4}}{\sqrt{\qnum{2}}} &
\qcg{(1,0)}{(0,-1)}{(1,0)}{(-1,1)}{(0,1)}{(-1,0)} &= -\frac{q^{-1/4}}{\sqrt{\qnum{2}}} \ .
\end{align}

With the information we have given so far, it is possible to obtain all the $q$-Clebsch-Gordan
coefficients which do not involve the representation ${\bf 8}$. In the following subsection, we
will explain the subtleties which arise when dealing with this representation.

\subsection{Dealing with multiplicities}
\label{sec:multiplicities}

In this section, we will explain how to deal with the two kinds of multiplicities, namely weight
space multiplicity and fusion multiplicity.

The eight-dimensional adjoint representation of $su(3)$ has the property that the `weight' $(0,0)$
corresponds to a two-dimensional weight space. In other words, the two states\\
$\lmb\lma |(1,1)(1,1)\rangle$ and $\lma\lmb |(1,1)(1,1)\rangle$ are linearly independent. Thus, we will have
to choose a basis for this weight space. We pick the following basis, which is orthonormal.
\begin{equation}
|(1,1)(0,0)_{\pm} \rangle =
\frac{\bigl(\lma\lmb \pm \lmb\lma\bigr)}{\sqrt{2(\qnum{2}\pm 1)}} |(1,1)(1,1)\rangle \ .
\end{equation}
We now have the following relations for the action of the lowering operators on the states \lq above' the weight $(0,0)$ space
\begin{align}
\lma |(1,1)(2,-1)\rangle &=
\sqrt{\frac{\qnum{2}+1}{2}} |(1,1)(0,0)_{+} \rangle + \sqrt{\frac{\qnum{2}-1}{2}} |(1,1)(0,0)_{-} \rangle \\ 
\lmb |(1,1)(-1,2)\rangle &=
\sqrt{\frac{\qnum{2}+1}{2}} |(1,1)(0,0)_{+} \rangle - \sqrt{\frac{\qnum{2}-1}{2}} |(1,1)(0,0)_{-} \rangle  \ , 
\end{align}
while the action of raising and lowering operators on $|(1,1)(0,0)_{\pm} \rangle $ is given by 
\begin{align}
\lpa |(1,1)(0,0)_{\pm} \rangle &= \sqrt{\frac{\qnum{2} \pm 1}{2}} |(1,1)(2,-1)\rangle &
\lpb |(1,1)(0,0)_{\pm} \rangle &= \pm \sqrt{\frac{\qnum{2} \pm 1}{2}} |(1,1)(-1,2)\rangle \\ 
\lma |(1,1)(0,0)_{\pm} \rangle &= \sqrt{\frac{\qnum{2} \pm 1}{2}} |(1,1)(-2,1)\rangle &
\lmb |(1,1)(0,0)_{\pm} \rangle &= \pm \sqrt{\frac{\qnum{2} \pm 1}{2}} |(1,1)(1,-2)\rangle  
\end{align}

The fact that one has to choose a basis for the two-dimensional weight space $|(1,1)(0,0)_{\pm}\rangle$ is not the only subtlety which arises in conjunction with the representation ${\bf 8}$.
Consider the $su(3)$ tensor product
${\bf 8}\otimes {\bf 8} =
{\bf 1} \oplus {\bf 8} \oplus {\bf 8} \oplus {\bf 10} \oplus \overline{\bf 10} \oplus {\bf 27}$.
One finds that in the decomposition, the eight dimensional representation appears twice. At level $k=2$, or in other words, for $q=e^{2 \pi i/5}$, only one of these two eight-dimensional representations is present in the fusion product , which reads ${\bf 8}\times {\bf 8} = {\bf 1} + {\bf 8}$.
One is thus led to the question, how can one decide which eight-dimensional representation to pick?

In solving the highest weight conditions $\lpa |(1,1)(1,1)\rangle = \lpb |(1,1)(1,1)\rangle = 0$, one finds a
two-dimensional space of solutions, as one should. A convenient way of writing the two
solutions is as follows
{\small
\begin{eqnarray}
|(1,1)(1,1)\rangle_{1} =
\frac{q^{3/4}}{\sqrt{\qnum{4}+1}} |(1,1)(1,1)\rangle \otimes |(1,1)(0,0)_{+}\rangle
+\frac{q^{-3/4}}{\sqrt{\qnum{4}+1}}|(1,1)(0,0)_{+}\rangle \otimes |(1,1)(1,1)\rangle \nonumber \\
- \sqrt{\frac{\qnum{2}+1}{2(\qnum{4}+1)}} \Bigl( |(1,1)(-1,2)\rangle \otimes |(1,1)(2,-1)\rangle +|(1,1)(2,-1)\rangle \otimes |(1,1)(-1,2)\rangle \Bigr) \\
|(1,1)(1,1)\rangle_{2} 
=
\frac{q^{3/4}}{\sqrt{\qnum{4}-1}} |(1,1)(1,1)\rangle \otimes |(1,1)(0,0)_{-}\rangle
-\frac{q^{-3/4}}{\sqrt{\qnum{4}-1}}|(1,1)(0,0)_{-}\rangle \otimes |(1,1)(1,1)\rangle \nonumber \\
+ \sqrt{\frac{\qnum{2}-1}{2(\qnum{4}-1)}} \Bigl( |(1,1)(2,-1)\rangle \otimes |(1,1)(-1,2)\rangle -|(1,1)(-1,2)\rangle \otimes |(1,1)(2,1)\rangle \Bigr) 
\end{eqnarray}}
We observe that for generic $q$, the states are orthonormal. However, for $q=e^{2 \pi i/5}$, we
have $\qnum{4}=1$, which means that the state $|(1,1)(1,1)\rangle_{2}$ is not normalizable\footnote{One might be tempted to multiply the given expression for this state by $\sqrt{\qnum{4}-1}$ to obtain a state with a good limit as $q\rightarrow e^{2 \pi i/5}$. This indeed yields a state with finite coefficients, but it will have norm zero with respect to our complex conjugation free inner product.}. Thus,
we conclude that at level $k=2$, the state $|(1,1)(1,1)\rangle_{1}$ corresponds to the eight-dimensional
representation which is present in fusion product ${\bf 8}\times {\bf 8} = {\bf 1}+{\bf 8}$, while the state
$|(1,1)(1,1)\rangle_{2}$ corresponds to the representation which is present in the tensor product, but not
in the fusion product. 

At level $k=3$, i.e. $q=e^{2 \pi i/6}$, both eight dimensional representations are present in the
fusion product: $su(3)_3$ is one of the simplest theories which exhibits a fusion multiplicity.
Therefore, we will also give the $q$-Clebsch-Gordan coefficients for the $su(3)_3/Z_3$ theory,
which contains four representations, $\id$, ${\bf 8}$, ${\bf 10}$ and $\overline{\bf 10}$,
see appendix \ref{appsu33z3} for more details.
In particular, we have the following fusion product:
\begin{equation}
{\bf 8}\times {\bf 8} = \id + {\bf 8} + {\bf 8'} + {\bf 10} + \overline{\bf 10} \ ,
\end{equation}
where the first ${\bf 8}$ in the fusion product corresponds to the one present at level $k=2$, while
the second ${\bf 8'}$ corresponds to the one appearing for the first time at $k=3$. So, in the case
at hand, there is a natural choice of a basis for the two-dimensional space corresponding to this
fusion multiplicity, namely in `order of appearance' when the level $k$ increases.
In fact, for $su(3)$, it is always possible to make such a choice, because the
so called `threshold levels' for the fusion coefficients of $su(3)$ differ by one. We note, however,
that this property is special for $su(3)$. For more on this, we refer to paragraph 16.4 of \cite{fms}.
%%%%%%
% Example: su(4)_3 has the field (1,1,1), appearing for the first time at level 3, but nevertheless,
% (1,1,1) \times (1,1,1) = (0,0,0) + (0,1,2) + (2,1,0) + (0,2,0) + 2 (1,0,1), with fusion multiplicity!
%%%%%%

With this information, we dealt with all the subtleties arising from the adjoint representation of $su(3)$.

\subsection{Gauge convention for the $q$-Clebsch-Gordan coefficients}
\label{sec:gauge}

In the calculation of the $q$-Clebsch-Gordan coefficients, one must choose conventions
for the overall phases of the coefficients. In this section, we will explicitly describe our choice.
To fix the gauge for the symbols $\qcg{j_1}{m_1}{j_2}{m_2}{j}{m}$, for a particular choice of
$j_1$, $j_2$ and $j$, it suffices to fix the phase of one of the symbols. It is convenient to fix the
phase of the symbol with $m=j$, and the highest possible weight $m_1=m_{\rm max}$.
We have chosen this phase in such a way that in the limit $q\rightarrow1$, this symbol
$\qcg{j_1}{m_{\rm max}}{j_2}{j-m_{\rm max}}{j}{j}$ is real and positive.

This gauge choice for the $q$-Clebsch-Gordan coefficients also completely specifies the
gauge degrees of freedom for the $q-6j$ symbols, because they are fixed by the $q$-Clebsch-Gordan
coefficients. Our gauge choice is such that the $q-6j$ symbols with a trivial particle on any of the incoming lines is equal to $1$, when no vertices with a
fusion multiplicity are present. When $q$ is the primitive root of unity corresponding to the associated Chern-Simons theory, the $q-6j$ symbols with a trivial particle on the outgoing line is also equal to $1$ (again assuming there are no vertex multiplicities) .

\subsection{Symmetries of the $q$-Clebsch-Gordan coefficients}
\label{symandconv}

In this section, we give the symmetries of the $q$-Clebsch-Gordan coefficients. Due to these
symmetries, we do not have to calculate all the $q$-Clebsch-Gordan coefficients.
We will first state these symmetry relations, and explain them afterwards.
\begin{align}
\qcg{j_2}{m_2}{j_1}{m_1}{j}{m} &= (-1)^{s_1} \cg{j_1}{m_1}{j_2}{m_2}{j}{m}{\frac{1}{q}}
\label{cswap}\\
\qcg{j_1}{-\overline{m_1}}{j_2}{-\overline{m_2}}{j}{-\overline{m}} &=
(-1)^{s_2+\#(0,0)_-}\cg{j_1}{m_1}{j_2}{m_2}{j}{m}{\frac{1}{q}}
\label{lowbar} \\
\qcg{\overline{j_1}}{\overline{m_1}}{\overline{j_2}}{\overline{m_2}}{\overline{j}}{\overline{m}} &=
(-1)^{s_3+\#(0,0)_-} \qcg{j_1}{m_1}{j_2}{m_2}{j}{m} \ ,
\label{allbar} 
\end{align}
where, as before, we used the $j$'s and $m$'s to denote the highest weights and the weights of the states in the representations respectively. By the conjugate weights (denoted by the overline), we mean the following. In the case
of $su(3)$ the weight conjugate to $m=(\lambda_1,\lambda_2)$ is
$\overline{m} = \overline{(\lambda_1,\lambda_2)} = (\lambda_2,\lambda_1)$. However, for
$so(5)$ and $G_2$, weights are self conjugate
$\overline{m} = \overline{(\lambda_1,\lambda_2)} = (\lambda_1,\lambda_2)$.
In the equations above $s_1$, $s_2$ and $s_3$ are either $0$ or $1$, depending on the case at hand
(see the tables in appendices \ref{su32abc} and \ref{su33abc}). 
In addition, with $\#(0,0)_-$, we mean the number of $m_1$, $m_2$, $m$ which are equal to
$(0,0)_-$ in the representation $\bf{8}$ of $su(3)$, see below for an explanation.

The first relation eq.~\eqref{cswap} is related to the symmetry of the tensor products. When $j_1=j_2$, the value of $s_1$ is gauge invariant and can not be chosen. However, for $j_1\neq j_2$, the value of $s_1$ is a gauge choice. As noted in section~\ref{sec:gauge}, we always choose the overal phase of the states in such a way that the $q$-Clebsch-Gordan coefficient with $m=j$ and with the highest possible value of $m_1$ is positive in the limit $q\rightarrow 1$. Note that the highest possible value of $m_1$ is not always equal to $j_1$. This choice completely fixes the gauge and hence also $s_1$. In the tables in appendices \ref{su32abc} and \ref{su33abc}, we give the values of $s_1$ in our gauge.

The relation eq.~\eqref{lowbar} relates the $q$-Clebsch-Gordan coefficients `within' a certain representation. This relation stems from the fact that one can equally well obtain the whole representation by acting with raising operators on the lowest weights, instead of acting with lowering operators on the highest weights, as we do here. The value of $s_2$ depends on the representation at hand; this value will be specified in the tables in appendices \ref{su32abc} and \ref{su33abc}.
Finally, the last symmetry relation \eqref{allbar} corresponds to the symmetry under complete conjugation.

There is one important issue, which occurs for both eqs \eqref{lowbar} and \eqref{allbar}, which is related to the conjugation of the weight $(0,0)_{-}$ in the eight-dimensional representation of $su(3)$. Loosely speaking, we have the relation $\overline{(0,0)_{-}} = - (0,0)_{-}$. This means that for every time the (conjugated) weight $\overline{(0,0)_{-}}$ appears in eqns. \eqref{lowbar} and \eqref{allbar}, we get an additional sign.

Finally, we also note that the only case in this paper
for which $s_3=1$ is the ${\bf 8'}$ representation which appears in the fusion product 
${\bf 8}\times {\bf 8} = \id + {\bf 8} + {\bf 8'} + {\bf 10} + \overline{\bf 10}$ of $su(3)_3/Z_3$.

\section{Obtaining the $R$-symbols}
\label{rsymbols}

After having obtained the $q$-Clebsch-Gordan coefficients, it is rather straightforward
to obtain the $R$-symbols as well. We will make use of the explicit expression for
the $R$-matrix. This expression is rather cumbersome for arbitrary quantum groups, but
we will only need a very limited amount of information contained in this expression to
extract the $R$-symbols. 
In particular, the expression contains a product over
all positive roots. The ordering of these positive roots is important, but we will not
specify this ordering, because it will be immaterial for our purposes.
In addition to this product, there is a pre-factor, which we will need to obtain the $R$-matrix
elements. 

With these caveats in place, we give the {\em structure} of the $R$ matrix:
\begin{equation}
R = q^{(F_{\rm qf})_{i,j}\frac{h_i \otimes h_j}{2}} \prod_{\alpha > 0} 
E_{q^{-1}}^{(q^{\frac{1}{2}}-q^{-\frac{1}{2}})
(e_\alpha q^{\frac{h_\alpha}{4}} \otimes q^{-\frac{h_\alpha}{4}} f_\alpha)} \ ,
\end{equation}
where $E_q^x = \sum_{n=0}^\infty \frac{q^{\frac{n}{4}(n-1)}}{\qnum{n}!} x^n$.

The overall factor, $q^{(F_{\rm qf})_{i,j}\frac{h_i \otimes h_j}{2}}$, contains the 
quadratic form matrix $F_{\rm qf}$, whose components in terms of the inverse
Cartan matrix elements and the integers $t_j$ read $(F_{\rm qf})_{ij} = (A^{-1})_{ij}/t_j$,
where there is no summation over $j$. Recall that
$t_{i}=\frac{2}{(\alpha_{i},\alpha_{i})}$, with $\alpha_i$ the simple roots of the algebra.
For simply-laced Lie algebras, such as $su(r+1)$, all the $t_i=1$. The short roots of
$B_r$, $C_r$ and $F_4$ have $t=2$, while for the short root of $G_2$, one has $t=3$.

The product is over all positive roots of the algebra and 
the factors consist of exponentials of raising operators $e_\alpha$ acting on the left tensor factor, and
lowering operators $f_\alpha$ acting on the right tensor factors. It is this structure, which will allow us
to obtain the $R$-symbols from the pre-factor alone, as we will explain below, and
give an explicit example for clarification.

In the calculation of the $R$-symbols one first acts with the $R$-matrix
on a particular tensor decomposition, followed by swapping the tensor factors.
This combined operation will be denoted by $\sigma R$. This will give rise to the
tensor decomposition in which the representations are swapped, up to an overall
phase factor. The $R$-symbol is precisely this phase factor. We will now
explain the calculation of the $R$-symbols $R^{j_1,j_2}_{j}$ in detail.

First of all, the $R$-symbols $R^{j_1,j_2}_{j}$ depend only on the
representations $j_1$, $j_2$ and $j$, and not on the particular weights within these
representations, as is the case for the the $F$-symbols. We can thus pick a
particular weight in the representation $j$ to obtain $R^{j_1,j_2}_{j}$, which we will take to
be $m=j$. We will restrict ourselves further by only considering a suitably chosen term in the
decomposition.

Namely, we will pick that component of the tensor decomposition which has a lowest
possible weight in the left factor of the tensor product. By this we mean that one can
not subtract a root from this weight, and obtain another term in the tensor decomposition.

Because all the raising operators in the $R$-matrix act on the left factor of the tensor product,
the only contribution to the component in the tensor product with this lowest possible weight
on the right (after performing the swap $\sigma$!) comes from the component with the
lowest possible weight on the left. In addition, it is only the identity term
in all exponents which contributes, because all the other terms contain raising operators on the
left tensor factor. Thus, the knowledge of the two corresponding
$q$-Clebsch-Gordan coefficients and the factor $q^{(F_{\rm qf})_{i,j}\frac{h_i \otimes h_j}{2}}$
suffices to obtain the $R$-symbols.

To show how this works, we give an explicit example. Let us consider the $\overline{\bf 3}$
in ${\bf 3}\times{\bf 3}$, and take the highest weight of $\overline{\bf 3}$:
\begin{equation}
\label{rcalc}
|(0,1)(0,1)\rangle =
\frac{q^{1/4}}{\sqrt{\qnum{2}}}  |(1,0)(1,0)\rangle \otimes |(1,0)(-1,1)\rangle 
- \frac{q^{-1/4}}{\sqrt{\qnum{2}}}  |(1,0)(-1,1)\rangle \otimes |(1,0)(1,0)\rangle \ .
\end{equation}
We let $\sigma R$ act on the term
$- \frac{q^{-1/4}}{\sqrt{\qnum{2}}}  |(1,0)(-1,1)\rangle \otimes |(1,0)(1,0)\rangle$
and find
\begin{equation}
\begin{split}
\sigma R \bigl(- \frac{q^{-1/4}}{\sqrt{\qnum{2}}}  |(1,0)(-1,1)\rangle \otimes |(1,0)(1,0)\rangle \bigr) &=
 \frac{-q^{-1/6} q^{-1/4}}{\sqrt{\qnum{2}}}  |(1,0)(1,0)\rangle \otimes |(1,0)(-1,1)\rangle \\
&= \bigl(-q^{-2/3}\bigr) \frac{q^{1/4}}{\sqrt{\qnum{2}}}  |(1,0)(1,0)\rangle \otimes |(1,0)(-1,1)\rangle
\end{split}
\end{equation}
Comparing this result to the first term in eq.~\eqref{rcalc}, we find
that $R^{{\bf 3},{\bf 3}}_{\overline{\bf 3}} = - q^{-2/3}$. In general (including the
non-simply laced case), we obtain the following expression for the R-matrix elements
\begin{equation}
R^{j_1,j_2}_{j} = q^{\frac{1}{2}m_{\rm min}\cdot F_{\rm qf} \cdot (j-m_{\rm min})}
\frac{\qcg{j_1}{m_{\rm min}}{j_2}{j-m_{\rm min}}{j}{j}}{\qcg{j_2}{j-m_{\rm min}}{j_1}{m_{\rm min}}{j}{j}} \ ,
\end{equation}
where $m_{\rm min}$ is the lowest possible weight mentioned earlier, i.e.~$m_{\rm min}$ is  chosen so that
$\qcg{j_1}{m_{\rm min}-\alpha}{j_2}{j-m_{\rm min}+\alpha}{j}{j}=0$
for all simple roots $\alpha$, and nevertheless the CG-coefficients in the formula above are nonzero.
Note that we assume that the weights are given in terms of their Dynkin labels (hence the appearance of the quadratic form matrix $F_{\rm qf}$).
%and note that there might be more than one $m_{\rm min}$, which
%would all give the same result.
%In the pre-factor, $F_{\rm qf}$ is the quadratic
%form matrix, whose components in terms of the inverse Cartan matrix elements and the
%integers $d_j$ read $(F_{\rm qf})_{ij} = (a^{-1})_{ij}/d_j$, where there is no summation over $j$.
%Note that for simply-laced Lie algebras, such as $su(n)$, all the $d_j=1$.
Finally, let us note two symmetries of the $R$-symbols which exists in our gauge,
\begin{align}
R^{a,b}_c &= R^{b,a}_c & R^{a,b}_c &= R^{\bar{a},\bar{b}}_{\bar{c}} \ .
\end{align}

\section{TQFT data}
\label{tqftdata}

In this section, we briefly describe some of the topological data associated with
a topological quantum field theory. We will not go into a lot of detail to explain the
theory behind this data. For this, we refer the reader to, for instance, \cite{Kitaev03} or
\cite{bonderson-thesis}.
We will frequently use the data of $su(3)_2$ as an example, and refer to the appendices
for the corresponding data of the other models. We note that the formulas for the
various quantities hold in the case of unitary anyon models, which corresponds here to the cases where $q$ is a primitive root of unity. 

\subsection{Quantum dimensions}

The quantum dimension $d_a$ associated with a particle of type $a$ will
be defined in the following way:
\begin{equation}
d_a = \frac{1}{|\bigl(\fs{a}{\overline{a}}{a}{a}\bigr)_{\id,\id}|}
\end{equation}
Though this is not obvious from this definition, the Hilbert space of $n$ particles of
type $a$ grows as $d_a^n$. 
The quantum dimensions of the six representations of $su(3)_2$ are given by
\begin{align}
d_{\id} &=1 & 
d_{{\bf 3}} &= d_{\bar{\bf 3}} = \qnum{3} \rightarrow \phi&
d_{{\bf 6}} &= d_{\bar{\bf 6}} = 1+ \qnum{5} \rightarrow 1& 
d_{{\bf 8}} &= \qnum{3} + \qnum{5} \rightarrow \phi \ ,
\end{align}
where the numerical values are obtained for $q=e^{\frac{2 \pi i}{5}}$ and $\phi = \frac{1+\sqrt{5}}{2}$.
In particular, the total quantum dimension
$D^2 = \sum_a d^2_a = 3(2+\phi)$.

\subsection{Frobenius-Schur indicator}

The phase, or in our gauge choice, the sign of the F-symbol
$\bigl(\fs{a}{\overline{a}}{a}{a}\bigr)_{\id,\id}$,
is called the Frobenius-Schur indicator
${\rm fb}_a = d_a \bigl(\fs{a}{\overline{a}}{a}{a}\bigr)_{\id,\id}$.
In the case of a self-dual particle,
$\overline{a}=a$, this is a gauge invariant quantity, which takes the values $\pm 1$. 
In the case of non-self dual particles, the Frobenius Schur indicator is often defined to be zero, but we will not do this here and we simply keep the definition above. Concretely, this means we have
${\rm fb}_{\id} = {\rm fb}_{\bf 3} = {\rm fb}_{\bf \overline{3}} =
{\rm fb}_{\bf 6} = {\rm fb}_{\bf \overline{6}} = {\rm fb}_{\bf 8} = 1$.

In section \ref{twistfactors}, we will give an expression for the Frobenius-Schur indicator
of the self-dual particles (of modular theories) in terms of the so-called twist factors.

\subsection{Theta-symbols}
In this section, we will give the value of the `theta' symbols. These depend on our choice of gauge.
By the `theta' symbols, mean the following diagrams.
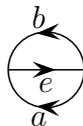
\begin{center}
\psset{unit=.5mm,linewidth=.2mm,dimen=middle}
\begin{pspicture}(0,-5)(20,25)
\psset{arrowsize=4pt 4}
\psarc{->}(10,10){10}{0}{105}
\psarc(10,10){10}{90}{180}
\psline{->}(0,10)(12.5,10)
\psline(10,10)(20,10)
\psarc(10,10){10}{180}{270}
\psarc{<-}(10,10){10}{255}{360}
\rput(8,-3){$a$}
\rput(8,24){$b$}
\rput(10,6){$e$}
\end{pspicture}
\end{center}
We will denote the value of the theta symbols
by $\vartheta (a,b,e)$.
These symbols can easily be evaluated, by applying the appropriate $F$ symbol, and
by noting that tadpoles give zero contributions
\begin{equation}
\vartheta (a,b,e) = (\fs{a}{b}{\bar{b}}{a})_{e,\id} d_a d_b \ .
\end{equation}
Our choice of gauge implies that
$(\fs{a}{b}{\bar{b}}{a})_{e,\id} = \pm\sqrt{\frac{d_e}{d_a d_b}}$, so we find that
$\vartheta (a,b,e) = \pm \sqrt{d_a d_b d_e}$.
The theta symbols have the following symmetries
\begin{align}
\vartheta (a,b,e) &= \vartheta(\bar{a},\bar{b},\bar{e}) &
\vartheta (a,b,e) & = \vartheta (b,a,e) &
\vartheta (a,b,e) & = \vartheta (\bar{e},b,\bar{a})  \ .
\end{align}
Thus, by specifying the following values, all the theta symbols for $su(3)_2$ are determined
\begin{align}
\vartheta ({\bf 3},\bar{\bf 3},{\bf 3}) &= - d_{\bf 3}^{\frac{3}{2}} &
\vartheta (\bar{\bf 3},{\bf 3},{\bf 6}) &= d_{\bf 3} \sqrt{d_{\bf 6}} &
\vartheta ({\bf 3},{\bf 3},{\bf 8}) &= d_{\bf 3} \sqrt{d_{\bf 8}} &
\vartheta ({\bf 3},\bar{\bf 6},{\bf 8}) &= - \sqrt{d_{\bf 3} d_{\bf 6} d_{\bf 8}} \\
\vartheta ({\bf 6},\bar{\bf 6},{\bf 6}) &= d_{\bf 6}^{\frac{3}{2}} &
\vartheta ({\bf 8},{\bf 8},{\bf 8}) &= d_{\bf 8}^{\frac{3}{2}} \ .&
\end{align}
Note that if any of the labels is the vacuum representation, the theta symbol
is just the quantum dimension of the remaining representation.

\subsection{Twist factors}
\label{twistfactors}

The twist factors, or topological spin of a particle $a$ will be denoted by $\theta_a$. In general,
the twist factors are given by
\begin{equation}
\theta_{a} = \theta_{\overline{a}} = {\rm fb}_{a} (R^{\overline{a},a}_{\id})^{*}
\end{equation}

In those cases where the TQFT's are associated with conformal field theories, $\theta_a$ is also
given by $\theta_a = e^{2 \pi i h_a}$, where $h_a$ is the scaling dimension of the corresponding
primary field in the CFT. Again, taking $su(3)_k$ as an example, 
\begin{align}
\theta_{\id} &=1 & 
\theta_{{\bf 3}} &= \theta_{\bar{\bf 3}} = e^{2 \pi i h_3} = e^{\frac{8 \pi i}{3(k+3)}}= q^{4/3} \\
\theta_{{\bf 6}} &= \theta_{\bar{\bf 6}} = e^{2 \pi i h_6} = e^{\frac{20 \pi i}{3(k+3)}} = q^{10/3}& 
\theta_{{\bf 8}} &= e^{2 \pi i h_8} = e^{\frac{6 \pi i}{(k+3)}} = q^{3}
\end{align}
In addition, one can express the Frobenius-Schur indicator for self-dual particles in terms of
the fusion-coefficients and twist factors as follows:

\begin{equation}
{\rm fb}_c = \frac{1}{D^2} \sum_{a,b} {n_{a,b}}^{c} \Bigl(\frac{\theta_a}{\theta_b}\Bigr)
d_a d_b \ ,
\end{equation}
where the summation variables $a$ and $b$ run over all particle types, and
${n_{a,b}}^{c}$ denote the fusion coefficient.
In particular, we find that both ${\rm fb}_\id = {\rm fb}_{\rm 8} = 1$, for $su(3)_2$.

\subsection{Central charge}
The central charge of the theory can be determined (modulo $8$) by using the result
\begin{equation}
e^{\frac{2 \pi i c}{8}} = \frac{1}{D} \sum_a d_a^2 \theta_a \ ,
\end{equation}
which leads to the expected result $c = \frac{16}{5}$ for $su(3)_2$.

\subsection{The tetrahedral symbol}

Finally, we mention a relation between the F-symbols and the so-called
tetrahedral symbols. The tetrahedral symbols are related to the F-symbols
one-to-one, and are proportional the the F-symbols. The advantage of the
tetrahedral symbols is that they satisfy more symmetry properties. In actual
calculations, however, one always needs the F-symbols, so will tabulate
those instead in the appendices.

The tetrahedral symbol is represented as
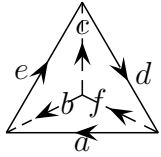
\begin{center}
\psset{unit=.5mm,linewidth=.2mm,dimen=middle,arrowsize=4pt 4}
\begin{pspicture}(0,0)(40,40)
\rput{60}(0,0){\psline{->}(0,0)(22,0)
\psline(20,0)(40,0)}
\rput{180}(40,0){\psline{->}(0,0)(22,0)
\psline(20,0)(40,0)}
\rput{120}(40,0){\psline{-<}(0,0)(22,0)
\psline(20,0)(40,0)}
\psline[linestyle=dashed]{-<}(0,0)(13,6.5)
\psline[linestyle=dashed](10,5)(20,10)
\psline[linestyle=dashed]{->}(40,0)(27,6.5)
\psline[linestyle=dashed](30,5)(20,10)
\psline[linestyle=dashed]{->}(20,10)(20,24)
\psline[linestyle=dashed](20,22)(20,34.641)
\rput(20,-3){$a$}
\rput(16,8){$b$}
\rput(20,28){$c$}
\rput(36,16){$d$}
\rput(4,16){$e$}
\rput(24,8){$f$}
\end{pspicture}
\end{center}
To calculate the value of the tetrahedral symbols, one applies the F-symbol
$(\fs{a}{b}{c}{d})_{e,f'}$. The only non-vanishing term in the sum is the one with
$f'=f$, and the resulting graph is
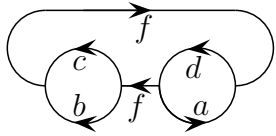
\begin{center}
\psset{unit=.5mm,linewidth=.2mm,dimen=middle,arrowsize=4pt 4}
\begin{pspicture}(0,0)(70,30)
\psarc{->}(20,10){10}{0}{105}
\psarc(20,10){10}{90}{180}
\psarc(20,10){10}{180}{270}
\psarc{<-}(20,10){10}{255}{360}
\psarc{->}(50,10){10}{0}{105}
\psarc(50,10){10}{90}{180}
\psarc(50,10){10}{180}{360}
\psarc{->}(50,10){10}{180}{285}
%\psline(0,10)(10,10)
\psline(30,10)(40,10)
%\psline(60,10)(70,10)
\psarc(10,20){10}{90}{270}
\psarc(60,20){10}{-90}{90}
\psline(10,30)(60,30)
\psline{->}(20,30)(38,30)
\psline{-<}(30,10)(38,10)
\rput(36,25){$f$}
\rput(34,4){$f$}
\rput(19,5){$b$}
\rput(19,16){$c$}
\rput(49,15){$d$}
\rput(51,4){$a$}
\end{pspicture}
\end{center}
The last graph is easily found to give
$(\fs{f}{\bar{f}}{f}{f})_{\id,\id} \vartheta(a,d,f) \vartheta(b,c,f)$.
Putting everything together, and using 
$\left|(\fs{f}{\bar{f}}{f}{f})_{\id,\id}\right| = \frac{1}{d_f}$, we find that
\begin{equation}
\cG(a,b,c,d,e,f) =
{\rm sgn}(\vartheta(a,d,f)) {\rm sgn}(\vartheta(b,c,f)) {\rm fb}_{f}
(\fs{a}{b}{c}{d})_{e,f} \sqrt{d_a d_b d_c d_d} \ ,
\end{equation} 
where ${\rm sgn}(\vartheta(a,b,c)) = \frac{\vartheta(a,b,c)}{|\vartheta(a,b,c)|}$, etc.

As we stated before, the tetrahedral symbols have more symmetry than the
F-symbols. This symmetry arises from the symmetry of the tetrahedron.
Thus, the following relations generate the full symmetry. Note that we have to take
the conjugate representations in those cases where the arrows are reversed before
and after `rotation' of the tetrahedron.
\begin{align}
\label{gsyms}
\cG (a,b,c,d,e,f) &= \cG (d,\bar{f},b,e,a,\bar{c}) = \cG (e,c,\bar{f},a,d,\bar{b}) \\
\cG (a,b,c,d,e,f) &= \cG (c,\bar{d},a,\bar{b},\bar{f},\bar{e}) \nonumber \\
\cG (a,b,c,d,e,f) &= \cG (\bar{c},\bar{b},\bar{a},\bar{d},\bar{f},\bar{e}) \nonumber \ .
\end{align}

\section{Applications, discussion and outlook}
\label{outlook}

The topological data calculated in this paper can be applied to calculate various
properties of physical systems. Particularly interesting examples of such systems occur in the context of the fractional quantum Hall effect. It has been proposed that (some of) the topological properties of fractional quantum Hall states can be measured by means of interference experiments, see \cite{topintmr} for the Moore-Read state \cite{mr91} case, \cite{topintrr} for the Read-Rezayi states \cite{rr99} and \cite{geninterferometry} for an account of the general case. Recent and promising experiments in this area are described in \cite{Willett09,Dolev08,Radu08}.

Of course, apart from what may be experimentally measurable in the immediate future, it is also of interest to determine the full representation of the braid group which governs the exchange properties of the non-abelian anyons in these Hall states, not in the least because a knowledge of these braid group representations or indeed of the full TQFTs describing these states is needed in the design of gates and algorithms for future topological quantum computers based on these states. In \cite{sb01}, the braiding of quasiparticles in the the Read-Rezayi states was investigated using quantum groups and we hope the present paper makes it obvious how similar calculations can be done for more complicated systems. Braiding properties can also be deduced by considering the full CFT correlation functions describing the non-abelian excitations, as was done in \cite{nw96} for the Moore-Read state. The $U_q(su(3))$ coefficients calculated in this paper match up with CFT-correlator calculations in \cite{as07}, where the $su(3)_k$ based spin-singlet state of Ardonne and Schoutens \cite{as99} was considered, as well as with the Read-Rezayi states.

An interesting recent development is the study of anyonic quantum spin chains \cite{anyonicchains,disorderedanyons}. Non-abelian anyons are the basic constituents of these chains, and the $F$-symbols are essential in the construction of the hamiltonian. For more details, we refer to \cite{ttw08}.

On a more formal level, the coefficients calculated can be used to study relations between different anyon models, such as those based on coset-theories. For instance, from the explicit calculation of the $F$-symbols of the $Z_3$-parafermion theory in \cite{bs-pentagon}, where the pentagon equations were solved directly, we know that these symbols are equivalent those calculated here for $su(3)_2$. 

A potentially interesting mathematical problem would be to calculate $q$-deformed Clebsch Gordan and $6j$-coefficients for the full set of irreducible representations of $U_{q}(sl(3))$, at all levels. The multiplicities could be dealt with by the method of threshold levels discussed in section~\ref{sec:multiplicities}. This should also yield a particularly  interesting set of bases for the CG-coefficients and $6j$-symbols of $SU(3)$ itself in the limit where $q$ approaches $1$.  

\vspace*{2mm}\noindent \textbf{\em Acknowledgements} -- We gratefully acknowledge many useful discussions with Tobias Hagge, Zhenghang Wang and,
in particular, Lukasz Fidkowski.\\
The authors would also like to acknowledge the institutions they worked for since the
start of this project, namely, Microsoft Station Q, Caltech, and UC Riverside (JS). We also thank the KITP in Santa Barbara, where this project was initiated, for its hospitality.\\
JS was supported by Science Foundation Ireland through PI Award 08/IN.1/I1961.

\appendix

\section{The case $su(3)$}
\label{appsu3}

In this appendix, we will give the topological data, $q$-Clebsch-Gordan coefficients and
$F$-symbols in the case of $su(3)_2$. However, we will first summarize the data for the
finite-dimensional Lie-algebra $su(3)$. For much more detail on the theory of finite (and
infinite) dimensional Lie algebras, see for instance, \cite{fms,fs,cftqgconnection}.

The $q$-Clebsch-Gordan coefficients, i.e. 
$$
\qcg{j_1}{m_1}{j_2}{m_2}{j_3}{m_3}
$$
will be tabulated in the following way. 
First, the various sections will be labeled by the explicit fusion rule $j_1 \times j_2$. In
this section, one finds the coefficients of all the possible fusion outcomes $j_3$.
The coefficients for a particular $j_3$ are given in a table, whose rows are labeled by
$m_1$, and the columns by $m_2$. In the cases that there are no weight space multiplicities,
this uniquely specifies $m_3 = m_1+m_2$. In the case of the eight dimensional representation
of $su(3)$, ${\bf 8}$, the weight $(0,0)$ corresponds to a two-dimensional space, whose basis
states we denote by $(0,0)_+$ and $(0,0)_-$. Often we can compactly combine the symbols
for $m_3 = (0,0)_\pm$. In cases where this is not easily possible, the 
the $q$-CG coefficients with $m_3=(0,0)_{+}$ and $m_3 = (0,0)_{-}$ are given in separate
tables, while the locations in the original table for which $m_1 + m_2 = (0,0)$ are marked
by $X$.

\subsection{Generalities on $su(3)$}

The Cartan matrix of the simply laces lie algebra $su(3)$, and it's inverse, read
\begin{align}
A &= \begin{pmatrix} 2 & -1 \\ -1 & 2 \end{pmatrix} & 
A^{-1} &=  F_{\rm qf}= \frac{1}{3}\begin{pmatrix} 2 & 1 \\ 1 & 2 \end{pmatrix} & 
\end{align}

The representations relevant for the $su(3)_2$ theory are
$\id = (0,0)$, ${\bf 3} = (1,0)$, ${\bf \overline{3}} = (0,1)$,
${\bf 6} = (2,0)$, ${\bf \overline{6}} = (0,2)$ and ${\bf 8} = (1,1)$.
In figures \ref{fig:33b} and \ref{fig:66b8}, we give the structure of these representations
for completeness.
Here and below, arrows pointing `to the left' correspond to the subtraction of
$\alpha_1$, while arrows pointing `to the right' correspond to the subtraction
of $\alpha_2$.
\begin{figure}[ht]
\begin{center}
\psset{unit=1mm,linewidth=.2mm,dimen=middle,arrowsize=4pt 4}
\begin{pspicture}(-5,-2)(15,22)
\rput(10,20){$(1,0)$}
\rput(0,10){$(-1,1)$}
\rput(10,0){$(0,-1)$}
\psline{->}(10,18)(0,12)
\psline{->}(0,8)(10,2)
\end{pspicture}
\hspace{2 cm}
\begin{pspicture}(-5,-2)(15,22)
\rput(0,20){$(0,1)$}
\rput(10,10){$(1,-1)$}
\rput(0,0){$(-1,0)$}
\psline{->}(0,18)(10,12)
\psline{->}(10,8)(0,2)
\end{pspicture}
\end{center}
\caption{The weights of the $su(3)$ representations ${\bf 3}$ and ${\bf \overline{3}}$.}
\label{fig:33b}
\end{figure}

\begin{figure}[h]
\begin{center}
\psset{unit=1mm,linewidth=.2mm,dimen=middle,arrowsize=4pt 4}
\begin{pspicture}(-5,-2)(25,42)
\rput(20,40){$(2,0)$}
\rput(10,30){$(0,1)$}
\rput(0,20){$(-2,2)$}
\rput(20,20){$(1,-1)$}
\rput(10,10){$(-1,0)$}
\rput(20,0){$(0,-2)$}
\psline{->}(20,38)(10,32)
\psline{->}(10,28)(0,22)
\psline{->}(10,28)(20,22)
\psline{->}(0,18)(10,12)
\psline{->}(20,18)(10,12)
\psline{->}(10,8)(20,2)
\end{pspicture}
\hspace{2 cm}
\begin{pspicture}(-5,-2)(25,42)
\rput(0,40){$(0,2)$}
\rput(10,30){$(1,0)$}
\rput(20,20){$(2,-2)$}
\rput(0,20){$(-1,1)$}
\rput(10,10){$(0,-1)$}
\rput(0,0){$(-2,0)$}
\psline{->}(0,38)(10,32)
\psline{->}(10,28)(0,22)
\psline{->}(10,28)(20,22)
\psline{->}(0,18)(10,12)
\psline{->}(20,18)(10,12)
\psline{->}(10,8)(0,2)
\end{pspicture}
\hspace{2 cm}
\begin{pspicture}(-5,-2)(25,42)
\rput(10,40){$(1,1)$}
\rput(0,30){$(-1,2)$}
\rput(20,30){$(2,-1)$}
\rput(11.5,20){$(0,0)_{\pm}$}
\rput(0,10){$(-2,1)$}
\rput(20,10){$(1,-2)$}
\rput(10,0){$(-1,-1)$}
\psline{->}(10,38)(0,32)
\psline{->}(10,38)(20,32)
\psline{->}(0,28)(10,22)
\psline{->}(20,28)(10,22)
\psline{->}(10,18)(0,12)
\psline{->}(10,18)(20,12)
\psline{->}(0,8)(10,2)
\psline{->}(20,8)(10,2)
\end{pspicture}
\end{center}
\caption{The weights of the $su(3)$ representations ${\bf 6}$, ${\bf \overline{6}}$ and ${\bf 8}$.}
\label{fig:66b8}
\end{figure}

\subsection{The topological data of $su(3)_2$}

In this section, we give the relevant topological data for $su(3)_2$.
The fusion rules of the theory are given in table \ref{su32fr}. The numerical
values of the data are obtained from the general ones in terms of $q$ by
specifying $q=e^{2 \pi i/5}$. 

\begin{table}[h]
\begin{center}
\begin{tabular}{|c|c|c|c|c|c|}
\hline
$\times$ & ${\bf 3}$ & ${\bf \overline{3}}$ & ${\bf 6}$ & ${\bf \overline{6}}$ & ${\bf 8}$ \\
\hline
${\bf 3}$ & ${\bf \overline{3}} + {\bf 6}$ &&&&\\
${\bf \overline{3}}$ & $\id + {\bf 8}$ &  ${\bf 3} + {\bf \overline{6}}$ &&&\\
${\bf 6}$ & ${\bf 8}$ & ${\bf 3}$ & ${\bf \overline{6}}$ &&\\
${\bf \overline{6}}$ & ${\bf \overline{3}}$ & {\bf 8} & $\id$ & ${\bf 6}$ &\\
${\bf 8}$ & ${\bf 3} + {\bf \overline{6}}$ & ${\bf \overline{3}} + {\bf 6}$ & ${\bf \overline{3}}$ &
${\bf 3}$ & $\id + {\bf 8}$ \\
\hline
\end{tabular}
\caption{The fusion rules of $su(3)_2$.}
\label{su32fr}
\end{center}
\end{table}

The quantum dimensions are given by
\begin{align}
d_{\id} &=1 & 
d_{{\bf 3}} &= d_{\bar{\bf 3}} = \qnum{3} \rightarrow \phi&
d_{{\bf 6}} &= d_{\bar{\bf 6}} = 1+ \qnum{5} \rightarrow 1& 
d_{{\bf 8}} &= \qnum{3} + \qnum{5} \rightarrow \phi \ ,
\end{align}
which gives $D^2 = \sum_a d_a^2 = 3\phi + 6$.
The twist factors are given by
\begin{align}
\theta_{\id} &=1 & 
\theta_{{\bf 3}} &= \theta_{\bar{\bf 3}} = e^{2 \pi i h_3} = e^{\frac{8 \pi i}{3(k+3)}}= q^{4/3} \\
\theta_{{\bf 6}} &= \theta_{\bar{\bf 6}} = e^{2 \pi i h_6} = e^{\frac{20 \pi i}{3(k+3)}} = q^{10/3}& 
\theta_{{\bf 8}} &= e^{2 \pi i h_8} = e^{\frac{6 \pi i}{(k+3)}} = q^{3} \ .
\end{align}
The Frobenius-Schur indicators of the self dual particles are ${\rm fb}_\id = {\rm fb}_{\bf 8} = 1$ and
finally, the central charge is $c=\frac{16}{5}$.

\subsection{The $q$-CG coefficients relevant for $su(3)_2$}

In this section, we will give the $q$-CG coefficients in the case of $su(3)_2$.

\subsubsection{${\bf 3}\times{\bf 3}=\bar{\bf 3}+{\bf 6}$}
\label{su3example}

\begin{tabular}{c|ccc}
& $(1,0)$ & $(-1,1)$ & $(0,-1)$\\
\hline
$(1,0)$ & & $\frac{q^{\frac{1}{4}}}{\sqrt{\qnum{2}}}$ & $\frac{q^{\frac{1}{4}}}{\sqrt{\qnum{2}}}$ \\
$(-1,1)$ & $-\frac{q^{-\frac{1}{4}}}{\sqrt{\qnum{2}}}$ & & $\frac{q^{\frac{1}{4}}}{\sqrt{\qnum{2}}}$ \\
$(0,-1)$ & $-\frac{q^{-\frac{1}{4}}}{\sqrt{\qnum{2}}}$ & $-\frac{q^{-\frac{1}{4}}}{\sqrt{\qnum{2}}}$ &  \\
\multicolumn{4}{c}{\phantom{$5$}}
\\
\multicolumn{4}{c}{The $q$-CG coefficients for the $\bar{\bf 3}$ in ${\bf 3}\times{\bf 3}$.}
\end{tabular}
\begin{tabular}{c|ccc}
& $(1,0)$ & $(-1,1)$ & $(0,-1)$\\
\hline
$(1,0)$ & 1 & $\frac{q^{-\frac{1}{4}}}{\sqrt{\qnum{2}}}$ & $\frac{q^{-\frac{1}{4}}}{\sqrt{\qnum{2}}}$ \\
$(-1,1)$ & $\frac{q^{\frac{1}{4}}}{\sqrt{\qnum{2}}}$ & 1 & $\frac{q^{-\frac{1}{4}}}{\sqrt{\qnum{2}}}$ \\
$(0,-1)$ & $\frac{q^{\frac{1}{4}}}{\sqrt{\qnum{2}}}$ & $\frac{q^{\frac{1}{4}}}{\sqrt{\qnum{2}}}$ & 1 \\
\multicolumn{4}{c}{\phantom{$5$}}
\\
\multicolumn{4}{c}{The $q$-CG coefficients for the ${\bf 6}$ in ${\bf 3}\times{\bf 3}$.}
\end{tabular}

\subsubsection{${\bf 3}\times\bar{\bf 3}=\id+{\bf 8}$}

\begin{tabular}{c|ccc}
& $(0,1)$ & $(1,-1)$ & $(-1,0)$\\
\hline
$(1,0)$ & & & $\frac{q^{\frac{1}{2}}}{\sqrt{\qnum{3}}}$ \\
$(-1,1)$ & & $-\frac{1}{\sqrt{\qnum{3}}}$ & \\
$(0,-1)$ & $\frac{q^{-\frac{1}{2}}}{\sqrt{\qnum{3}}}$ & &  \\
\multicolumn{4}{c}{\phantom{$5$}}
\\
\multicolumn{4}{c}{The $q$-CG coefficients for the $\id$ in ${\bf 3}\times\bar{\bf 3}$.}
\end{tabular}
\begin{tabular}{c|ccc}
& $(0,1)$ & $(1,-1)$ & $(-1,0)$\\
\hline
$(1,0)$ & 1 & 1 & $\frac{q^{-\frac{1}{4}}}{\sqrt{2(\qnum{2}\pm1)}}$ \\
$(-1,1)$ & 1 & $\frac{q^{\frac{1}{4}}\pm q^{-\frac{1}{4}}}{\sqrt{2(\qnum{2}\pm1)}}$ & 1 \\
$(0,-1)$ & $\pm\frac{q^{\frac{1}{4}}}{\sqrt{2(\qnum{2}\pm1)}}$ & 1 & 1 \\
\multicolumn{4}{c}{\phantom{$5$}}
\\
\multicolumn{4}{c}{The $q$-CG coefficients for the ${\bf 8}$ in ${\bf 3}\times\bar{\bf 3}$.}
\end{tabular}

\subsubsection{${\bf 3}\times {\bf 6} = {\bf 8}$}
\begin{tabular}{c|cccccc}
& $(2,0)$ & $(0,1)$ & $(-2,2)$ & $(1,-1)$ & $(-1,0)$ & $(0,-2)$ \\
\hline
$(1,0)$ &
&
$\frac{q^{\frac{1}{2}}}{\sqrt{\qnum{3}}}$ &
$q^{\frac{1}{4}}\sqrt{\frac{\qnum{2}}{\qnum{3}}}$ &
$\frac{q^{\frac{1}{2}}}{\sqrt{\qnum{3}}}$ &
$\pm q^{\frac{1}{4}}\sqrt{\frac{\qnum{2}\pm 1}{2\qnum{3}}}$ &
$q^{\frac{1}{4}}\sqrt{\frac{\qnum{2}}{\qnum{3}}}$
\\ 
$(-1,1)$ &
$-q^{-\frac{1}{4}}\sqrt{\frac{\qnum{2}}{\qnum{3}}}$ &
$-\frac{q^{-\frac{1}{2}}}{\sqrt{\qnum{3}}}$ &
&
$\frac{q^{\frac{3}{4}}\mp q^{-\frac{3}{4}}}{\sqrt{2(\qnum{2}\pm 1)\qnum{3}}}$ &
$\frac{q^{\frac{1}{2}}}{\sqrt{\qnum{3}}}$ &
$q^{\frac{1}{4}}\sqrt{\frac{\qnum{2}}{\qnum{3}}}$
\\
$(0,-1)$ &
$-q^{-\frac{1}{4}}\sqrt{\frac{\qnum{2}}{\qnum{3}}}$ &
$-q^{-\frac{1}{4}}\sqrt{\frac{\qnum{2}\pm 1}{2\qnum{3}}}$ &
$-q^{-\frac{1}{4}}\sqrt{\frac{\qnum{2}}{\qnum{3}}}$ &
$-\frac{q^{-\frac{1}{2}}}{\sqrt{\qnum{3}}}$ &
$-\frac{q^{-\frac{1}{2}}}{\sqrt{\qnum{3}}}$ &
\\
\multicolumn{7}{c}{\phantom{$5$}}
\\
\multicolumn{7}{c}{The $q$-CG coefficients for the ${\bf 8}$ in ${\bf 3}\times{\bf 6}$.}
\end{tabular}

\subsubsection{${\bf 3}\times \bar{\bf 6} = \bar{\bf 3}$}

\begin{tabular}{c|cccccc}
& $(0,2)$ & $(1,0)$ & $(-1,1)$ & $(2,-2)$ & $(0,-1)$ & $(-2,0)$ \\
\hline
$(1,0)$ &
&
&
$\frac{q^{\frac{3}{4}}}{\sqrt{\qnum{4}}}$ &
&
$\frac{q^{\frac{3}{4}}}{\sqrt{\qnum{4}}}$ &
$q^{\frac{1}{2}} \sqrt{\frac{\qnum{2}}{\qnum{4}}}$ \\ 
$(-1,1)$ &
&
$-\frac{q^{\frac{1}{4}}}{\sqrt{\qnum{4}}}$&
&
$-\sqrt{\frac{\qnum{2}}{\qnum{4}}}$&
$-\frac{q^{-\frac{1}{4}}}{\sqrt{\qnum{4}}}$&
\\
$(0,-1)$ &
$q^{-\frac{1}{2}} \sqrt{\frac{\qnum{2}}{\qnum{4}}}$ &
$\frac{q^{-\frac{3}{4}}}{\sqrt{\qnum{4}}}$ &
$\frac{q^{-\frac{3}{4}}}{\sqrt{\qnum{4}}}$ &
&
&
\\
\multicolumn{7}{c}{\phantom{$5$}}
\\
\multicolumn{7}{c}{The $q$-CG coefficients for the $\bar{\bf 3}$ in ${\bf 3}\times\bar{\bf 6}$.}
\end{tabular}

\subsubsection{${\bf 3}\times {\bf 8} = {\bf 3}+\bar{\bf 6}$}

\begin{tabular}{c|ccccccc}
& $(1,1)$ & $(-1,2)$ & $(2,-1)$ & $(0,0)_\pm$ & $(-2,1)$ & $(1,-2)$ & $(-1,-1)$ \\
\hline
$(1,0)$ &
&
&
&
$\frac{q^{\frac{3}{4}}\sqrt{\qnum{2} \mp 1}}{\sqrt{2\qnum{2}\qnum{4}}}$&
$\frac{q^{\frac{1}{2}}\sqrt{\qnum{3}}}{\sqrt{\qnum{2}\qnum{4}}}$&
&
$\frac{q^{\frac{1}{2}}\sqrt{\qnum{3}}}{\sqrt{\qnum{2}\qnum{4}}}$\\ 
$(-1,1)$ &
&
&
$-\frac{\sqrt{\qnum{3}}}{\sqrt{\qnum{2}\qnum{4}}}$&
$\frac{(\mp q^{\frac{1}{4}}-q^{-\frac{1}{4}})\sqrt{\qnum{2}\mp 1}}{\sqrt{2\qnum{2}\qnum{4}}}$&
&
$-\frac{\sqrt{\qnum{3}}}{\sqrt{\qnum{2}\qnum{4}}}$&
\\
$(0,-1)$ &
$\frac{q^{-\frac{1}{2}}\sqrt{\qnum{3}}}{\sqrt{\qnum{2}\qnum{4}}}$&
$\frac{q^{-\frac{1}{2}}\sqrt{\qnum{3}}}{\sqrt{\qnum{2}\qnum{4}}}$&
&
$\pm \frac{q^{-\frac{3}{4}}\sqrt{\qnum{2} \mp 1}}{\sqrt{2\qnum{2}\qnum{4}}}$&
&
&
\\
\multicolumn{7}{c}{\phantom{$5$}}
\\
\multicolumn{7}{c}{The $q$-CG coefficients for the ${\bf 3}$ in ${\bf 3}\times{\bf 8}$.}
\end{tabular}\\ ~~\\ ~~\\
\begin{tabular}{c|ccccccc}
& $(1,1)$ & $(-1,2)$ & $(2,-1)$ & $(0,0)_\pm$ & $(-2,1)$ & $(1,-2)$ & $(-1,-1)$ \\
\hline
$(1,0)$ &
&
$\frac{q^{\frac{1}{4}}}{\sqrt{\qnum{2}}}$&
&
$\pm \frac{q^{\frac{1}{4}}\sqrt{\qnum{2}\pm1}}{\sqrt{2}\qnum{2}}$&
$\frac{1}{\qnum{2}}$&
$\frac{q^{\frac{1}{4}}}{\sqrt{\qnum{2}}}$&
$\frac{1}{\qnum{2}}$\\ 
$(-1,1)$ &
$-\frac{q^{-\frac{1}{4}}}{\sqrt{\qnum{2}}}$&
&
$-\frac{q^{-\frac{1}{2}}}{\qnum{2}}$&
$\frac{(\pm q^{\frac{1}{4}}-q^{-\frac{1}{4}})\sqrt{\qnum{2}\pm 1}}{\sqrt{2}\qnum{2}}$&
&
$\frac{q^{\frac{1}{2}}}{\qnum{2}}$&
$\frac{q^{\frac{1}{4}}}{\sqrt{\qnum{2}}}$\\
$(0,-1)$ &
$-\frac{1}{\qnum{2}}$&
$-\frac{1}{\qnum{2}}$&
$-\frac{q^{-\frac{1}{4}}}{\sqrt{\qnum{2}}}$&
$-\frac{q^{-\frac{1}{4}}\sqrt{\qnum{2}\pm1}}{\sqrt{2}\qnum{2}}$&
$-\frac{q^{-\frac{1}{4}}}{\sqrt{\qnum{2}}}$&
&
\\
\multicolumn{7}{c}{\phantom{$5$}}
\\
\multicolumn{7}{c}{The $q$-CG coefficients for the $\bar{\bf 6}$ in ${\bf 3}\times{\bf 8}$.}
\end{tabular}\\

\subsubsection{${\bf 6}\times{\bf 6}=\bar{\bf 6}$}

\begin{tabular}{c|cccccc}
& $(2,0)$ & $(0,1)$ & $(-2,2)$ & $(1,-1)$ & $(-1,0)$ & $(0,-2)$ \\
\hline
$(2,0)$ &
&
&
$\frac{q^{\frac{1}{2}}}{\sqrt{\qnum{3}}}$&
&
$\frac{q^{\frac{1}{2}}}{\sqrt{\qnum{3}}}$&
$\frac{q^{\frac{1}{2}}}{\sqrt{\qnum{3}}}$\\ 
$(0,1)$ &
&
$-\frac{1}{\sqrt{\qnum{3}}}$&
&
$-\frac{q^{-\frac{1}{4}}}{\sqrt{\qnum{2}\qnum{3}}}$&
$\frac{q^{\frac{3}{4}}}{\sqrt{\qnum{2}\qnum{3}}}$&
$\frac{q^{\frac{1}{2}}}{\sqrt{\qnum{3}}}$\\ 
$(-2,2)$ &
$\frac{q^{-\frac{1}{2}}}{\sqrt{\qnum{3}}}$&
&
&
$-\frac{1}{\sqrt{\qnum{3}}}$&
&
$\frac{q^{\frac{1}{2}}}{\sqrt{\qnum{3}}}$\\ 
$(1,-1)$ &
&
$-\frac{q^{\frac{1}{4}}}{\sqrt{\qnum{2}\qnum{3}}}$&
$-\frac{1}{\sqrt{\qnum{3}}}$&
$-\frac{1}{\sqrt{\qnum{3}}}$&
$-\frac{q^{-\frac{1}{4}}}{\sqrt{\qnum{2}\qnum{3}}}$&
\\ 
$(-1,0)$ &
$\frac{q^{-\frac{1}{2}}}{\sqrt{\qnum{3}}}$&
$\frac{q^{-\frac{3}{4}}}{\sqrt{\qnum{2}\qnum{3}}}$&
&
$-\frac{q^{\frac{1}{4}}}{\sqrt{\qnum{2}\qnum{3}}}$&
$-\frac{1}{\sqrt{\qnum{3}}}$&
\\ 
$(0,-2)$ &
$\frac{q^{-\frac{1}{2}}}{\sqrt{\qnum{3}}}$&
$\frac{q^{-\frac{1}{2}}}{\sqrt{\qnum{3}}}$&
$\frac{q^{-\frac{1}{2}}}{\sqrt{\qnum{3}}}$&
&
&
\\ 
\multicolumn{7}{c}{\phantom{$5$}}
\\
\multicolumn{7}{c}{The $q$-CG coefficients for the $\bar{\bf 6}$ in ${\bf 6}\times{\bf 6}$.}
\end{tabular}\\

\subsubsection{${\bf 6}\times\bar{\bf 6}=\id$}

\begin{tabular}{c|cccccc}
& $(0,2)$ & $(1,0)$ & $(-1,1)$ & $(2,-2)$ & $(0,-1)$ & $(-2,0)$ \\
\hline
$(2,0)$ &
&
&
&
&
&
$\frac{q}{\sqrt{\qnum{5}+1}}$\\ 
$(0,1)$ &
&
&
&
&
$-\frac{q^{\frac{1}{2}}}{\sqrt{\qnum{5}+1}}$&
\\ 
$(-2,2)$ &
&
&
&
$\frac{1}{\sqrt{\qnum{5}+1}}$&
&
\\ 
$(1,-1)$ &
&
&
$\frac{1}{\sqrt{\qnum{5}+1}}$&
&
&
\\ 
$(-1,0)$ &
&
$-\frac{q^{-\frac{1}{2}}}{\sqrt{\qnum{5}+1}}$&
&
&
&
\\ 
$(0,-2)$ &
$\frac{q^{-1}}{\sqrt{\qnum{5}+1}}$&
&
&
&
&
\\ 
\multicolumn{7}{c}{\phantom{$5$}}
\\
\multicolumn{7}{c}{The $q$-CG coefficients for the $\id$ in ${\bf 6}\times\bar{\bf 6}$.}
\end{tabular}\\

\subsubsection{${\bf 6}\times {\bf 8} = \bar{\bf 3}$}

\begin{tabular}{c|ccccccc}
& $(1,1)$ & $(-1,2)$ & $(2,-1)$ & $(0,0)_\pm$ & $(-2,1)$ & $(1,-2)$ & $(-1,-1)$ \\
\hline
$(2,0)$ &
&
&
&
&
$\frac{q^{\frac{3}{4}}}{\sqrt{\qnum{4}}}$&
&
$\frac{q^{\frac{3}{4}}}{\sqrt{\qnum{4}}}$\\ 
$(0,1)$ &
&
&
&
$-\frac{q^{\frac{1}{4}}\sqrt{\qnum{2}\pm1}}{\sqrt{2\qnum{2}\qnum{4}}}$&
&
$-\frac{1}{\sqrt{\qnum{2}\qnum{4}}}$&
$\frac{q}{\sqrt{\qnum{2}\qnum{4}}}$\\ 
$(-2,2)$ &
&
&
$\frac{q^{-\frac{1}{4}}}{\sqrt{\qnum{4}}}$&
&
&
$-\frac{q^{\frac{1}{4}}}{\sqrt{\qnum{4}}}$&
\\ 
$(1,-1)$ &
&
$\frac{q^{-\frac{1}{2}}}{\sqrt{\qnum{2}\qnum{4}}}$&
&
$\frac{(-q^{\frac{1}{4}}\pm q^{-\frac{1}{4}})\sqrt{\qnum{2}\pm1}}{\sqrt{2\qnum{2}\qnum{4}}}$&
$-\frac{q^{\frac{1}{2}}}{\sqrt{\qnum{2}\qnum{4}}}$&
&
\\ 
$(-1,0)$ &
$-\frac{q^{-1}}{\sqrt{\qnum{2}\qnum{4}}}$&
&
$\frac{1}{\sqrt{\qnum{2}\qnum{4}}}$&
$\pm\frac{q^{-\frac{1}{4}}\sqrt{\qnum{2}\pm1}}{\sqrt{2\qnum{2}\qnum{4}}}$&
&
&
\\ 
$(0,-2)$ &
$-\frac{q^{-\frac{3}{4}}}{\sqrt{\qnum{4}}}$&
$-\frac{q^{-\frac{3}{4}}}{\sqrt{\qnum{4}}}$&
&
&
&
&
\\ 
\multicolumn{7}{c}{\phantom{$5$}}
\\
\multicolumn{7}{c}{The $q$-CG coefficients for the $\bar{\bf 3}$ in ${\bf 6}\times{\bf 8}$.}
\end{tabular}\\

\subsubsection{${\bf 8}\times {\bf 8} = \id + {\bf 8}$}

\hspace*{-8mm}
\rotatebox{90}{%
\begin{tabular}{c|cccccccc}
& $(1,1)$ & $(-1,2)$ & $(2,-1)$ & $(0,0)_+$ & $(0,0)_-$ & $(-2,1)$ & $(1,-2)$ & $(-1,-1)$ \\
\hline
$(1,1)$ &
&
&
&
$\frac{q^{\frac{3}{4}}}{\sqrt{\qnum{4}+1}}$&
$0$&
$\frac{q^{\frac{1}{2}}\sqrt{\qnum{2}+1}}{\sqrt{2(\qnum{4}+1)}}$&
$\frac{q^{\frac{1}{2}}\sqrt{\qnum{2}+1}}{\sqrt{2(\qnum{4}+1)}}$&
$X$\\
$(-1,2)$ &
&
&
$-\frac{\sqrt{\qnum{2}+1}}{\sqrt{2(\qnum{4}+1)}}$&
$\frac{-q^{-\frac{1}{4}}-q^{\frac{1}{4}}+q^{\frac{3}{4}}}{2\sqrt{\qnum{4}+1}}$&
$-\frac{q^{\frac{1}{4}}\sqrt{\qnum{3}}}{2\sqrt{\qnum{4}+1}}$&
&
$X$&
$\frac{q^{\frac{1}{2}}\sqrt{\qnum{2}+1}}{\sqrt{2(\qnum{4}+1)}}$\\
$(2,-1)$ &
&
$-\frac{\sqrt{\qnum{2}+1}}{\sqrt{2(\qnum{4}+1)}}$&
&
$\frac{-q^{-\frac{1}{4}}-q^{\frac{1}{4}}+q^{\frac{3}{4}}}{2\sqrt{\qnum{4}+1}}$&
$\frac{q^{\frac{1}{4}}\sqrt{\qnum{3}}}{2\sqrt{\qnum{4}+1}}$&
$X$&
&
$\frac{q^{\frac{1}{2}}\sqrt{\qnum{2}+1}}{\sqrt{2(\qnum{4}+1)}}$\\
$(0,0)_+$ &
$\frac{q^{-\frac{3}{4}}}{\sqrt{\qnum{4}+1}}$&
$\frac{-q^{-\frac{1}{4}}-q^{\frac{1}{4}}+q^{-\frac{3}{4}}}{2\sqrt{\qnum{4}+1}}$&
$\frac{-q^{-\frac{1}{4}}-q^{\frac{1}{4}}+q^{-\frac{3}{4}}}{2\sqrt{\qnum{4}+1}}$&
$X$&
$X$&
$\frac{-q^{-\frac{1}{4}}-q^{\frac{1}{4}}+q^{\frac{3}{4}}}{2\sqrt{\qnum{4}+1}}$&
$\frac{-q^{-\frac{1}{4}}-q^{\frac{1}{4}}+q^{\frac{3}{4}}}{2\sqrt{\qnum{4}+1}}$&
$\frac{q^{\frac{3}{4}}}{\sqrt{\qnum{4}+1}}$\\
$(0,0)_-$ &
$0$&
$-\frac{q^{-\frac{1}{4}}\sqrt{\qnum{3}}}{2\sqrt{\qnum{4}+1}}$&
$\frac{q^{-\frac{1}{4}}\sqrt{\qnum{3}}}{2\sqrt{\qnum{4}+1}}$&
$X$&
$X$&
$\frac{q^{\frac{1}{4}}\sqrt{\qnum{3}}}{2\sqrt{\qnum{4}+1}}$&
$-\frac{q^{\frac{1}{4}}\sqrt{\qnum{3}}}{2\sqrt{\qnum{4}+1}}$&
$0$\\
$(-2,1)$ &
$\frac{q^{-\frac{1}{2}}\sqrt{\qnum{2}+1}}{\sqrt{2(\qnum{4}+1)}}$&
&
$X$&
$\frac{-q^{-\frac{1}{4}}-q^{\frac{1}{4}}+q^{-\frac{3}{4}}}{2\sqrt{\qnum{4}+1}}$&
$\frac{q^{-\frac{1}{4}}\sqrt{\qnum{3}}}{2\sqrt{\qnum{4}+1}}$&
&
$-\frac{\sqrt{\qnum{2}+1}}{\sqrt{2(\qnum{4}+1)}}$&
\\
$(1,-2)$ &
$\frac{q^{-\frac{1}{2}}\sqrt{\qnum{2}+1}}{\sqrt{2(\qnum{4}+1)}}$&
$X$&
&
$\frac{-q^{-\frac{1}{4}}-q^{\frac{1}{4}}+q^{-\frac{3}{4}}}{2\sqrt{\qnum{4}+1}}$&
$-\frac{q^{-\frac{1}{4}}\sqrt{\qnum{3}}}{2\sqrt{\qnum{4}+1}}$&
$-\frac{\sqrt{\qnum{2}+1}}{\sqrt{2(\qnum{4}+1)}}$&
&
\\
$(-1,-1)$ &
$X$&
$\frac{q^{-\frac{1}{2}}\sqrt{\qnum{2}+1}}{\sqrt{2(\qnum{4}+1)}}$&
$\frac{q^{-\frac{1}{2}}\sqrt{\qnum{2}+1}}{\sqrt{2(\qnum{4}+1)}}$&
$\frac{q^{-\frac{3}{4}}}{\sqrt{\qnum{4}+1}}$&
$0$&
&
&
\\
\multicolumn{9}{c}{\phantom{$5$}}
\\
\multicolumn{9}{c}{The $q$-CG coefficients for the ${\bf 8}$ (when
$m \neq (0,0)_\pm$) in ${\bf 8}\times{\bf 8}$.}
\end{tabular}
}
\hspace*{5mm}
\rotatebox{90}{%
\begin{tabular}{c|cccccccc}
& $(1,1)$ & $(-1,2)$ & $(2,-1)$ & $(0,0)_+$ & $(0,0)_-$ & $(-2,1)$ & $(1,-2)$ & $(-1,-1)$ \\
\hline
$(1,1)$ &
&
&
&
&
&
&
&
$\frac{q^{\frac{1}{4}}}{\sqrt{\qnum{4}+1}}$\\
$(-1,2)$ &
&
&
&
&
&
&
$\frac{-q^{-\frac{1}{4}}+q^{\frac{1}{4}}+q^{\frac{3}{4}}}{2\sqrt{\qnum{4}+1}}$&
\\
$(2,-1)$ &
&
&
&
&
&
$\frac{-q^{-\frac{1}{4}}+q^{\frac{1}{4}}+q^{\frac{3}{4}}}{2\sqrt{\qnum{4}+1}}$&
&
\\
$(0,0)_+$ &
&
&
&
$\frac{-2(q^{-\frac{1}{4}}+q^{\frac{1}{4}})+q^{\frac{3}{4}}+q^{-\frac{3}{4}}}{2\sqrt{\qnum{4}+1}}$&
$0$&
&
&
\\
$(0,0)_-$ &
&
&
&
$0$&
$\frac{q^{\frac{3}{4}}+q^{-\frac{3}{4}}}{2\sqrt{\qnum{4}+1}}$&
&
&
\\
$(-2,1)$ &
&
&
$\frac{q^{-\frac{1}{4}}-q^{\frac{1}{4}}+q^{-\frac{3}{4}}}{2\sqrt{\qnum{4}+1}}$&
&
&
&
&
\\
$(1,-2)$ &
&
$\frac{q^{-\frac{1}{4}}-q^{\frac{1}{4}}+q^{-\frac{3}{4}}}{2\sqrt{\qnum{4}+1}}$&
&
&
&
&
&
\\
$(-1,-1)$ &
$\frac{q^{-\frac{1}{4}}}{\sqrt{\qnum{4}+1}}$ &
&
&
&
&
&
&
\\
\multicolumn{9}{c}{The $q$-CG coefficients for the ${\bf 8}$ (for $m = (0,0)_+$)
in ${\bf 8}\times{\bf 8}$.}
\end{tabular}
}

\begin{tabular}{c|cccccccc}
& $(1,1)$ & $(-1,2)$ & $(2,-1)$ & $(0,0)_+$ & $(0,0)_-$ & $(-2,1)$ & $(1,-2)$ & $(-1,-1)$ \\
\hline
$(1,1)$ &
&
&
&
&
&
&
&
$0$\\
$(-1,2)$ &
&
&
&
&
&
&
$\frac{q^{\frac{1}{4}}\sqrt{\qnum{3}}}{2\sqrt{\qnum{4}+1}}$&
\\
$(2,-1)$ &
&
&
&
&
&
$-\frac{q^{\frac{1}{4}}\sqrt{\qnum{3}}}{2\sqrt{\qnum{4}+1}}$&
&
\\
$(0,0)_+$ &
&
&
&
$0$&
$\frac{q^{\frac{3}{4}}+q^{-\frac{3}{4}}}{2\sqrt{\qnum{4}+1}}$&
&
&
\\
$(0,0)_-$ &
&
&
&
$\frac{q^{\frac{3}{4}}+q^{-\frac{3}{4}}}{2\sqrt{\qnum{4}+1}}$&
$0$&
&
&
\\
$(-2,1)$ &
&
&
$-\frac{q^{-\frac{1}{4}}\sqrt{\qnum{3}}}{2\sqrt{\qnum{4}+1}}$&
&
&
&
&
\\
$(1,-2)$ &
&
$\frac{q^{-\frac{1}{4}}\sqrt{\qnum{3}}}{2\sqrt{\qnum{4}+1}}$&
&
&
&
&
&
\\
$(-1,-1)$ &
$0$&
&
&
&
&
&
&
\\
\multicolumn{9}{c}{The $q$-CG coefficients for the ${\bf 8}$ (for $m = (0,0)_-$) in
${\bf 8}\times{\bf 8}$.}
\end{tabular}\\~~\\

\begin{tabular}{c|cccccccc}
& $(1,1)$ & $(-1,2)$ & $(2,-1)$ & $(0,0)_+$ & $(0,0)_-$ & $(-2,1)$ & $(1,-2)$ & $(-1,-1)$ \\
\hline
$(1,1)$ &
&
&
&
&
&
&
&
$\frac{q}{\sqrt{\qnum{2}\qnum{4}}}$\\
$(-1,2)$ &
&
&
&
&
&
&
$-\frac{q^{\frac{1}{2}}}{\sqrt{\qnum{2}\qnum{4}}}$&
\\
$(2,-1)$ &
&
&
&
&
&
$-\frac{q^{\frac{1}{2}}}{\sqrt{\qnum{2}\qnum{4}}}$&
&
\\
$(0,0)_+$ &
&
&
&
$\frac{1}{\sqrt{\qnum{2}\qnum{4}}}$&
$0$&
&
&
\\
$(0,0)_-$ &
&
&
&
$0$&
$\frac{1}{\sqrt{\qnum{2}\qnum{4}}}$&
&
&
\\
$(-2,1)$ &
&
&
$-\frac{q^{-\frac{1}{2}}}{\sqrt{\qnum{2}\qnum{4}}}$&
&
&
&
&
\\
$(1,-2)$ &
&
$-\frac{q^{-\frac{1}{2}}}{\sqrt{\qnum{2}\qnum{4}}}$&
&
&
&
&
&
\\
$(-1,-1)$ &
$\frac{q^{-1}}{\sqrt{\qnum{2}\qnum{4}}}$&
&
&
&
&
&
&
\\
\multicolumn{9}{c}{The $q$-CG coefficients for the $\id$ in ${\bf 8}\times{\bf 8}$.}
\end{tabular}

\subsubsection{The relation between the $q$-CG coefficients}
\label{su32abc}

In table \ref{su32sym}, we specify the coefficients $s_1$, $s_2$ and $s_3$,
which appear in the symmetry relations
between the various $q$-CG coefficients as explained in section \ref{symandconv}.

\begin{table}[h]
\begin{center}
\begin{tabular}{|c|c|c|r|r|r||c|c|c|r|r|r|}
\hline
$j_1$ & $j_2$ & $j$ & $s_1$ & $s_2$ & $s_3$ & $j_1$ & $j_2$ & $j$ & $s_1$ & $s_2$ & $s_3$\\
\hline
\hline
${\bf 3}$ & ${\bf 3}$ & $\bar{\bf 3}$ & $1$ & $1$ & $0$ &
${\bf 6}$ & ${\bf 6}$ & $\bar{\bf 6}$ & $0$ & $0$ & $0$ \\
${\bf 3}$ & ${\bf 3}$ & ${\bf 6}$ & $0$ & $0$ & $0$ &
${\bf 6}$ & $\bar{\bf 6}$ & $\id$ & $0$ & $0$ & $0$ \\
${\bf 3}$ & $\bar{\bf 3}$ & $\id$ & $0$ & $0$ & $0$ &
${\bf 6}$ & ${\bf 8}$ & $\bar{\bf 3}$ & $1$ & $1$ & $0$ \\
${\bf 3}$ & $\bar{\bf 3}$ & ${\bf 8}$ & $0$ & $0$ & $0$ &
${\bf 8}$ & ${\bf 8}$ & $\id$ & $0$ & $0$ & $0$ \\
${\bf 3}$ & ${\bf 6}$ & ${\bf 8}$ & $1$ & $1$ & $0$ &
${\bf 8}$ & ${\bf 8}$ & ${\bf 8}$ & $0$ & $0$ & $0$ \\
${\bf 3}$ & $\bar{\bf 6}$ & $\bar{\bf 3}$ & $0$ & $0$ & $0$ &&&&&&\\
${\bf 3}$ & ${\bf 8}$ & ${\bf 3}$ & $0$ & $0$ & $0$ &&&&&&\\
${\bf 3}$ & ${\bf 8}$ & $\bar{\bf 6}$ & $1$ & $1$ & $0$ &&&&&&\\
\hline
\end{tabular}
\end{center}
\caption{The parameters $s_1$, $s_2$ and $s_3$ in the symmetry relations between the $q$-CG
coefficients for $su(3)_2$, as explained in section \ref{symandconv}.}
\label{su32sym}
\end{table}

\subsection{The $F$-symbols for $su(3)_2$}

In total, there are 405 F-symbols.
Out of these, there are 147 symbols which have a $\id$ on at least one of the outer lines. In the
gauge we chose, all these symbols are equal to one. Moreover, because the F-symbols are
invariant under the operation of taking the conjugate representation of all 6 indices, namely
$\bigl(\fs{a}{b}{c}{d}\bigr)_{e,f} = \big(\fs{\bar a}{\bar b}{\bar c}{\bar d}\bigr)_{{\bar e},{\bar f}}$,
we will only list those symbols which either have $a={\bf 3}$, $a={\bf 6}$ or $a={\bf 8}$.

All those symbols which have at least one ${\bf 6}$ or $\overline{\bf 6}$ on the
outer lines will correspond to a one dimensional transformation, and reduce to $\pm 1$
for $q=e^{2 \pi i/5}$, but we will list these  symbols for general $q$. It turns out that there are
only 10 independent values.
\begin{align}
&\fs{\bf 3}{\bf 3}{\bar{\bf 3}}{\bar{\bf 6}} =
\fs{\bf 3}{\bf 3}{{\bf 6}}{{\bf 3}} =
\fs{{\bf 3}}{{\bf 3}}{{\bf 6}}{\bar{\bf 6}} = 
\fs{{\bf 3}}{\bar{\bf 3}}{{\bf 3}}{\bar{\bf 6}} = 
\fs{{\bf 3}}{\bar{\bf 3}}{\bar{\bf 3}}{{\bf 6}} = 
\fs{{\bf 3}}{\bar{\bf 3}}{{\bf 6}}{\bar{\bf 3}} = 
\fs{{\bf 3}}{\bar{\bf 3}}{\bar{\bf 6}}{{\bf 3}} = 
\fs{{\bf 3}}{{\bf 6}}{{\bf 3}}{{\bf 3}} =
\fs{{\bf 3}}{{\bf 6}}{{\bf 3}}{\bar{\bf 6}} =  \\ \nonumber
&\fs{{\bf 3}}{{\bf 6}}{\bar{\bf 3}}{\bar{\bf 3}} = 
\fs{{\bf 3}}{{\bf 6}}{{\bf 6}}{\bar{\bf 3}} = 
\fs{{\bf 3}}{\bar{\bf 6}}{\bar{\bf 3}}{{\bf 3}} = 
\fs{{\bf 6}}{{\bf 3}}{{\bf 3}}{{\bf 3}} = 
\fs{{\bf 6}}{{\bf 3}}{{\bf 3}}{\bar{\bf 6}} =
\fs{{\bf 6}}{{\bf 3}}{\bar{\bf 3}}{\bar{\bf 3}} = 
\fs{{\bf 6}}{{\bf 3}}{{\bf 6}}{\bar{\bf 3}} = 
\fs{{\bf 6}}{\bar{\bf 3}}{{\bf 3}}{\bar{\bf 3}} =
\fs{{\bf 6}}{{\bf 6}}{{\bf 3}}{\bar{\bf 3}} = 1 \\
%\end{align}
%\begin{align}
&\fs{{\bf 3}}{{\bf 3}}{\bar{\bf 6}}{{\bf 8}} =
\fs{{\bf 3}}{{\bf 3}}{{\bf 8}}{{\bf 6}} =
\fs{{\bf 3}}{\bar{\bf 6}}{\bar{\bf 6}}{{\bf 8}} =
\fs{{\bf 3}}{\bar{\bf 6}}{{\bf 8}}{\bar{\bf 3}} =
\fs{{\bf 3}}{{\bf 8}}{\bar{\bf 6}}{\bar{\bf 3}} =
\fs{{\bf 3}}{{\bf 8}}{\bar{\bf 6}}{{\bf 6}} =
\fs{{\bf 6}}{\bar{\bf 3}}{\bar{\bf 3}}{{\bf 8}} = 
\fs{{\bf 6}}{\bar{\bf 3}}{{\bf 8}}{\bar{\bf 6}} = 
\fs{{\bf 6}}{{\bf 6}}{\bar{\bf 3}}{{\bf 8}} = 
\\ \nonumber
&\fs{{\bf 6}}{{\bf 6}}{{\bf 8}}{{\bf 3}} =  
\fs{{\bf 6}}{{\bf 8}}{\bar{\bf 3}}{{\bf 3}} = 
\fs{{\bf 6}}{{\bf 8}}{\bar{\bf 3}}{\bar{\bf 6}} = 
\fs{{\bf 8}}{{\bf 3}}{{\bf 3}}{{\bf 6}} = 
\fs{{\bf 8}}{{\bf 3}}{\bar{\bf 6}}{{\bf 6}} =
\fs{{\bf 8}}{\bar{\bf 3}}{\bar{\bf 3}}{\bar{\bf 6}} = 
\fs{{\bf 8}}{\bar{\bf 3}}{{\bf 6}}{\bar{\bf 6}} = 
\fs{{\bf 8}}{{\bf 6}}{\bar{\bf 3}}{{\bf 3}} = 
\fs{{\bf 8}}{{\bf 6}}{{\bf 6}}{{\bf 3}} =
\\ \nonumber 
&\fs{{\bf 8}}{\bar{\bf 6}}{{\bf 3}}{\bar{\bf 3}} = 
\fs{{\bf 8}}{\bar{\bf 6}}{\bar{\bf 6}}{\bar{\bf 3}}
= -\frac{\qnum{3}}{\qnum{2}\sqrt{\qnum{5}+1}} \rightarrow -1 \\
%\end{align}
%\begin{align}
&\fs{{\bf 3}}{\bar{\bf 3}}{{\bf 6}}{{\bf 6}} =
\fs{{\bf 3}}{\bar{\bf 6}}{{\bf 6}}{{\bf 3}} =
\fs{{\bf 6}}{\bar{\bf 3}}{{\bf 3}}{{\bf 6}} =
\fs{{\bf 6}}{{\bf 6}}{\bar{\bf 6}}{{\bf 6}} =
\fs{{\bf 6}}{\bar{\bf 6}}{{\bf 3}}{{\bf 3}} =
\fs{{\bf 6}}{\bar{\bf 6}}{\bar{\bf 6}}{\bar{\bf 6}}
= \frac{1}{\sqrt{\qnum{5}+1}} \rightarrow 1 \\
%\end{align}
%\begin{align}
&\fs{{\bf 3}}{{\bf 6}}{\bar{\bf 3}}{{\bf 6}} =
\fs{{\bf 3}}{\bar{\bf 6}}{\bar{\bf 3}}{\bar{\bf 6}} =
\fs{{\bf 6}}{{\bf 3}}{\bar{\bf 6}}{{\bf 3}} =
\fs{{\bf 6}}{\bar{\bf 3}}{\bar{\bf 6}}{\bar{\bf 3}}
= -\frac{1}{\sqrt{\qnum{5}+1}} \rightarrow -1 \\
%\end{align}
%\begin{align}
&\fs{{\bf 3}}{\bar{\bf 3}}{\bar{\bf 6}}{\bar{\bf 6}} =
\fs{{\bf 3}}{{\bf 6}}{\bar{\bf 6}}{{\bf 3}} =
\fs{{\bf 6}}{{\bf 3}}{\bar{\bf 3}}{{\bf 6}} =
\fs{{\bf 6}}{\bar{\bf 6}}{\bar{\bf 3}}{\bar{\bf 3}} 
= -\frac{\qnum{2}}{\qnum{3}} \rightarrow -1 \\
%\end{align}
%\begin{align}
&\fs{{\bf 3}}{{\bf 6}}{{\bf 8}}{{\bf 8}} =
\fs{{\bf 3}}{{\bf 8}}{{\bf 8}}{\bar{\bf 6}} =
\fs{{\bf 6}}{{\bf 3}}{{\bf 8}}{{\bf 8}} =
\fs{{\bf 6}}{{\bf 8}}{{\bf 8}}{\bar{\bf 3}} =
\fs{{\bf 8}}{{\bf 3}}{{\bf 6}}{{\bf 8}} =
\fs{{\bf 8}}{\bar{\bf 3}}{\bar{\bf 6}}{{\bf 8}} = \\ \nonumber 
&\fs{{\bf 8}}{{\bf 6}}{{\bf 3}}{{\bf 8}} =
\fs{{\bf 8}}{\bar{\bf 6}}{\bar{\bf 3}}{{\bf 8}} =
\fs{{\bf 8}}{{\bf 8}}{{\bf 3}}{\bar{\bf 6}} =
\fs{{\bf 8}}{{\bf 8}}{\bar{\bf 3}}{{\bf 6}} =
\fs{{\bf 8}}{{\bf 8}}{{\bf 6}}{\bar{\bf 3}} =
\fs{{\bf 8}}{{\bf 8}}{\bar{\bf 6}}{{\bf 3}} 
= -\frac{1}{\sqrt{2}\qnum{2}}\sqrt{\frac{\qnum{2}+\qnum{3}}{\qnum{3}-1}} \rightarrow -1 \\
%\end{align}
%\begin{align}
&\fs{{\bf 3}}{\bar{\bf 6}}{{\bf 3}}{{\bf 8}} =
\fs{{\bf 3}}{\bar{\bf 6}}{{\bf 8}}{{\bf 6}} =
\fs{{\bf 3}}{{\bf 8}}{{\bf 3}}{{\bf 6}} =
\fs{{\bf 6}}{\bar{\bf 3}}{{\bf 6}}{{\bf 8}} =
\fs{{\bf 6}}{\bar{\bf 3}}{{\bf 8}}{{\bf 3}} =
\fs{{\bf 6}}{{\bf 8}}{{\bf 6}}{{\bf 3}} =
\fs{{\bf 8}}{{\bf 3}}{\bar{\bf 6}}{\bar{\bf 3}} =
\fs{{\bf 8}}{\bar{\bf 3}}{{\bf 6}}{{\bf 3}} =
\fs{{\bf 8}}{{\bf 6}}{\bar{\bf 3}}{\bar{\bf 6}} =
\\ \nonumber 
&\fs{{\bf 8}}{\bar{\bf 6}}{{\bf 3}}{{\bf 6}}
=\frac{1}{\qnum{4}} \rightarrow 1 \\
%\end{align}
%\begin{align}
&\fs{{\bf 6}}{\bar{\bf 6}}{{\bf 8}}{{\bf 8}} =
\fs{{\bf 6}}{{\bf 8}}{{\bf 8}}{{\bf 6}} =
\fs{{\bf 8}}{{\bf 6}}{\bar{\bf 6}}{{\bf 8}} =
\fs{{\bf 8}}{\bar{\bf 6}}{{\bf 6}}{{\bf 8}} =
\fs{{\bf 8}}{{\bf 8}}{{\bf 6}}{{\bf 6}} =
\fs{{\bf 8}}{{\bf 8}}{\bar{\bf 6}}{\bar{\bf 6}}
=-\frac{1}{\qnum{4}} \rightarrow -1 \\
%\end{align}
%\begin{align}
&\fs{{\bf 3}}{{\bf 8}}{{\bf 6}}{{\bf 8}} =
\fs{{\bf 6}}{{\bf 8}}{{\bf 3}}{{\bf 8}} =
\fs{{\bf 6}}{{\bf 8}}{\bar{\bf 6}}{{\bf 8}} =
\fs{{\bf 8}}{{\bf 3}}{{\bf 8}}{\bar{\bf 6}} =
\fs{{\bf 8}}{\bar{\bf 3}}{{\bf 8}}{{\bf 6}} =
\\ \nonumber 
&\fs{{\bf 8}}{{\bf 6}}{{\bf 8}}{\bar{\bf 3}} =
\fs{{\bf 8}}{{\bf 6}}{{\bf 8}}{{\bf 6}} =
\fs{{\bf 8}}{\bar{\bf 6}}{{\bf 8}}{{\bf 3}} =
\fs{{\bf 8}}{\bar{\bf 6}}{{\bf 8}}{\bar{\bf 6}}
=\frac{\qnum{3}}{\qnum{3}+\qnum{5}} \rightarrow 1 \\
%\end{align}
%\begin{align}
&\fs{{\bf 6}}{\bar{\bf 6}}{{\bf 6}}{{\bf 6}} 
=\frac{1}{\qnum{5}+1} \rightarrow 1
\end{align}

The remaining 72 symbols correspond to 18 different labelings of the external lines. We give
these 18 matrices explicitly below. Again, the internal labels $e$ and $f$ have to be inferred from the
others, and the ordering in the matrices is always according to
$\id,{\bf 3}, \bar{\bf 3}, {\bf 6},\bar{\bf 6},{\bf 8}$.

\begin{align}
&\fs{\bf 3}{\bar{\bf 3}}{\bar{\bf 3}}{\bar{\bf 3}} =
\fs{\bf 3}{\bf 3}{\bar{\bf 3}}{\bf 3} =
\begin{pmatrix}
\frac{-1}{\sqrt{\qnum{3}}} & \frac{\qnum{3}-1}{\sqrt{\qnum{5}+1}} \\
\frac{\qnum{3}-1}{\sqrt{\qnum{5}+1}} & \frac{1}{\sqrt{\qnum{3}}} 
\end{pmatrix}
\rightarrow
\begin{pmatrix}
-1/\sqrt{\phi} & 1/\phi \\
1/\phi & 1/\sqrt{\phi} \\
\end{pmatrix}
\\
&\fs{\bf 3}{\bar{\bf 3}}{\bf 3}{\bf 3} =
\begin{pmatrix}
\frac{1}{\qnum{3}} & \frac{\qnum{2}\sqrt{\qnum{3}-1}}{\qnum{3}} \\
\frac{\qnum{2}\sqrt{\qnum{3}-1}}{\qnum{3}}  & \frac{-1}{\qnum{3}} 
\end{pmatrix}
\rightarrow
\begin{pmatrix}
1/\phi & 1/\sqrt{\phi} \\
1/\sqrt{\phi} & -1/\phi \\
\end{pmatrix}
\\
&\fs{\bf 3}{\bar{\bf 3}}{\bf 8}{\bf 8} = 
\fs{\bf 8}{\bf 8}{\bar{\bf 3}}{\bar{\bf 3}} =
\fs{\bf 8}{\bf 8}{\bf 3}{\bf 3} =
\begin{pmatrix}
\frac{1}{\sqrt{\qnum{3}+\qnum{5}}} & \frac{-1}{\qnum{2}} \\
\frac{\qnum{4}+1}{\sqrt{2}\qnum{2}\sqrt{\qnum{5}+\qnum{4}+1}}
&
\frac{\sqrt{\qnum{2}+\qnum{3}}}{\sqrt{2}\qnum{2}}
\end{pmatrix}
\rightarrow
\begin{pmatrix}
1/\sqrt{\phi} & -1/\phi \\
1/\phi & 1/\sqrt{\phi} \\
\end{pmatrix}
\\
&\fs{\bf 3}{\bf 8}{\bf 8}{\bf 3} =
\fs{\bf 8}{\bar{\bf 3}}{\bf 3}{\bf 8} =
\fs{\bf 8}{\bf 3}{\bar{\bf 3}}{\bf 8} =
\begin{pmatrix}
\frac{1}{\sqrt{\qnum{3}+\qnum{5}}} & 
\frac{\qnum{4}+1}{\sqrt{2}\qnum{2}\sqrt{\qnum{5}+\qnum{4}+1}}
\\
\frac{-1}{\qnum{2}} &
\frac{\sqrt{\qnum{2}+\qnum{3}}}{\sqrt{2}\qnum{2}}
\end{pmatrix} 
\rightarrow
\begin{pmatrix}
1/\sqrt{\phi} & 1/\phi \\
-1/\phi & 1/\sqrt{\phi} \\
\end{pmatrix}
\\
&\fs{\bf 3}{\bf 3}{\bf 3}{\bf 8} =
\fs{\bf 3}{\bf 3}{\bf 8}{\bar{\bf 3}} =
\fs{\bf 3}{\bf 8}{\bf 3}{\bar{\bf 3}} =
\fs{\bf 8}{\bar{\bf 3}}{\bar{\bf 3}}{\bf 3} =
\fs{\bf 8}{\bf 3}{\bf 3}{\bar{\bf 3}} =
\begin{pmatrix}
\frac{-1}{\qnum{2}} & \frac{\sqrt{\qnum{3}}}{\qnum{2}} \\
\frac{\sqrt{\qnum{3}}}{\qnum{2}} & \frac{1}{\qnum{2}}
\end{pmatrix}
\rightarrow
\begin{pmatrix}
-1/\phi & 1/\sqrt{\phi} \\
1/\sqrt{\phi} & 1/\phi \\
\end{pmatrix}
\\
&\fs{\bf 3}{\bf 8}{\bar{\bf 3}}{\bf 8} =
\fs{\bf 8}{\bar{\bf 3}}{\bf 8}{\bar{\bf 3}} =
\fs{\bf 8}{\bf 3}{\bf 8}{\bf 3} =
\begin{pmatrix}
\frac{-1}{\qnum{3}+\qnum{5}} &
-\frac{\qnum{3}}{\qnum{2}^2} \sqrt{\frac{\qnum{3}}{\qnum{5}+1}} \\
-\frac{\qnum{3}}{\qnum{2}^2} \sqrt{\frac{\qnum{3}}{\qnum{5}+1}} &
\frac{\qnum{3}}{\qnum{2}^2} 
\end{pmatrix}
\rightarrow
\begin{pmatrix}
-1/\phi & -1/\sqrt{\phi} \\
-1/\sqrt{\phi} & 1/\phi \\
\end{pmatrix}
\\
&\fs{\bf 8}{\bf 8}{\bf 8}{\bf 8} =
\begin{pmatrix}
\frac{1}{\qnum{3}+\qnum{5}} & \frac{1}{\sqrt{\qnum{3}+\qnum{5}}} \\
\frac{1}{\sqrt{\qnum{3}+\qnum{5}}} & \frac{\qnum{3}-\qnum{2}-4}{2(\qnum{2}+\qnum{3})}
\end{pmatrix}
\rightarrow
\begin{pmatrix}
1/\phi & 1/\sqrt{\phi} \\
1/\sqrt{\phi} & -1/\phi \\
\end{pmatrix}
\end{align}

\subsection{The $R$-symbols for $su(3)_2$}

We give the $R$-symbols below. We use the following
symmetries to shorten the list
\begin{align}
R^{a,b}_c &= R^{b,a}_c & R^{a,b}_c &= R^{\bar{a},\bar{b}}_{\bar{c}} \ .
\end{align}
With our conventions, we have
\begin{align}
R^{\id,a}_{a} &= 1 \\
R^{{\bf 3},\bar{\bf 3}}_{\id} &= q^{-\frac{4}{3}} &
R^{{\bf 6},\bar{\bf 6}}_{\id} &= q^{-\frac{10}{3}} &
R^{{\bf 8},{\bf 8}}_{\id} &= q^{-3} \\
R^{\bar{\bf 3},\bar{\bf 3}}_{\bf 3} &= -q^{-\frac{2}{3}} &
R^{\bar{\bf 3},{\bf 6}}_{\bf 3} &= q^{-\frac{5}{3}} &
R^{{\bf 3},{\bf 8}}_{\bf 3} &= q^{-\frac{3}{2}} &
R^{\bar{\bf 6},{\bf 8}}_{\bf 3} &= -q^{-\frac{5}{2}} \\
R^{{\bf 3},{\bf 3}}_{\bf 6} &= q^{\frac{1}{3}} &
R^{\bar{\bf 3},{\bf 8}}_{\bf 6} &= -q^{-\frac{1}{2}} &
R^{\bar{\bf 6},\bar{\bf 6}}_{\bf 6} &= q^{-\frac{5}{3}} \\
R^{{\bf 3},\bar{\bf 3}}_{\bf 8} &= q^{\frac{1}{6}} &
R^{{\bf 3},{\bf 6}}_{\bf 8} &= -q^{\frac{5}{6}} &
R^{{\bf 8},{\bf 8}}_{\bf 8} &= q^{\frac{3}{2}} \ .
\end{align}

\section{The case $su(3)_3/Z_3$}
\label{appsu33z3}

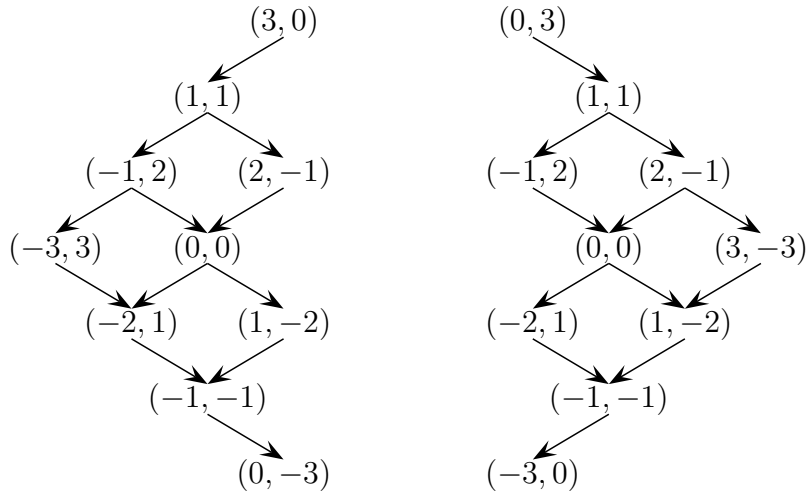
\begin{figure}[h]
\begin{center}
\psset{unit=1mm,linewidth=.2mm,dimen=middle,arrowsize=4pt 4}
\begin{pspicture}(-5,-2)(35,62)
\rput(30,60){$(3,0)$}
\rput(20,50){$(1,1)$}
\rput(10,40){$(-1,2)$}
\rput(30,40){$(2,-1)$}
\rput(0,30){$(-3,3)$}
\rput(20,30){$(0,0)$}
\rput(10,20){$(-2,1)$}
\rput(30,20){$(1,-2)$}
\rput(20,10){$(-1,-1)$}
\rput(30,0){$(0,-3)$}
\psline{->}(30,58)(20,52)
\psline{->}(20,48)(30,42)
\psline{->}(20,48)(10,42)
\psline{->}(10,38)(0,32)
\psline{->}(10,38)(20,32)
\psline{->}(30,38)(20,32)
\psline{->}(0,28)(10,22)
\psline{->}(20,28)(10,22)
\psline{->}(20,28)(30,22)
\psline{->}(10,18)(20,12)
\psline{->}(30,18)(20,12)
\psline{->}(20,8)(30,2)
\end{pspicture}
\hspace{2 cm}
\begin{pspicture}(-5,-2)(35,62)
\rput(0,60){$(0,3)$}
\rput(10,50){$(1,1)$}
\rput(0,40){$(-1,2)$}
\rput(20,40){$(2,-1)$}
\rput(30,30){$(3,-3)$}
\rput(10,30){$(0,0)$}
\rput(0,20){$(-2,1)$}
\rput(20,20){$(1,-2)$}
\rput(10,10){$(-1,-1)$}
\rput(0,0){$(-3,0)$}
\psline{->}(0,58)(10,52)
\psline{->}(10,48)(0,42)
\psline{->}(10,48)(20,42)
\psline{->}(20,38)(30,32)
\psline{->}(20,38)(10,32)
\psline{->}(0,38)(10,32)
\psline{->}(30,28)(20,22)
\psline{->}(10,28)(20,22)
\psline{->}(10,28)(0,22)
\psline{->}(20,18)(10,12)
\psline{->}(0,18)(10,12)
\psline{->}(10,8)(0,2)
\end{pspicture}
\end{center}
\caption{The weights of the $su(3)$ representations ${\bf 10}$ and ${\bf \overline{10}}$.}
\label{fig:1010b}
\end{figure}

In this appendix, we give the topological data for the theory $su(3)_3/Z_3$, namely the
$q$-CG coefficients, as well as the $F$ and $R$ symbols. This theory contains four particles,
which we will denote by $\id$, ${\bf 8}$, ${\bf 10}$ and ${\bf \overline{10}}$, where the last
two correspond to the $su(3)$ representations $(3,0)$ and $(0,3)$ respectively;
the weights of these representations are given in figure \ref{fig:1010b}.
We note that
this theory is a 'sub-theory' of the full, modular $su(3)_3$ theory. However, the theory under
consideration is not modular.

The fusion rules are given in table \ref{su33fr}.
\begin{table}[h]
\begin{center}
\begin{tabular}{|c|c|c|c|}
\hline
$\times$ & ${\bf 8}$ & ${\bf 10}$ & ${\bf \overline{10}}$ \\
\hline
${\bf 8}$ & $\id + {\bf 8} + {\bf 8} + {\bf 10} + {\bf \overline{10}}$ &&\\
${\bf 10}$ & ${\bf 8}$ & ${\bf \overline{10}}$ & \\
${\bf \overline{10}}$ & ${\bf 8}$ & $\id$ & ${\bf 10}$\\
\hline
\end{tabular}
\end{center}
\caption{The fusion rules of $su(3)_3/Z_3$.}
\label{su33fr}
\end{table}

The quantum dimensions of the particles are given by
$d_{\id} = 1$,
$d_{\bf 8} = \qnum{3}+\qnum{5} \rightarrow 3$,
$d_{\bf 10} = d_{\bf \overline{10}} = (\qnum{3}-1) \qnum{5} \rightarrow 1$,
where the numerical value is for $q=e^{2 \pi i/6}$. 
The twist factors are given by $\theta_\id=1$, $\theta_{\bf 8} = q^3$
and $\theta_{\bf 10} = \theta_{\bf \overline{10}} = q^6$.
The Frobenius-Schur indicator of the self-dual particles $\id$ and ${\bf 8}$ is the same as for
$su(3)_2$, namely ${\rm fb}_\id = {\rm fb}_{\bf 8} = 1$, while also
${\rm fb}_{\bf 10} = {\rm fb}_{\bf \overline{10}} = 1$.

%\vfill

\subsection{$q$-CG coefficients for $su(3)_3/Z_3$}

In this section, we will give the $q$-CG coefficients in the case of $su(3)_3/Z_3$.
The Clebsch-Gordan coefficients for the
$\id$ and ${\bf 8}$ in the ${\bf 8}\times{\bf 8}$ are given in the section containing
the Clebsch-Gordan coefficients for $su(3)_2$.

\subsubsection{${\bf 8}\times {\bf 8} = \id + {\bf 8}+ {\bf 8'} + {\bf 10} + \overline{\bf 10}$}

\rotatebox{90}{%
\begin{tabular}{c|cccccccc}
& $(1,1)$ & $(-1,2)$ & $(2,-1)$ & $(0,0)_+$ & $(0,0)_-$ & $(-2,1)$ & $(1,-2)$ & $(-1,-1)$ \\
\hline
$(1,1)$ &
&
&
&
$0$&
$\frac{q^{\frac{3}{4}}}{\sqrt{\qnum{4}-1}}$&
$\frac{q^{\frac{1}{2}}\sqrt{\qnum{2}-1}}{\sqrt{2(\qnum{4}-1)}}$&
$-\frac{q^{\frac{1}{2}}\sqrt{\qnum{2}-1}}{\sqrt{2(\qnum{4}-1)}}$&
$X$\\
$(-1,2)$ &
&
&
$-\frac{\sqrt{\qnum{2}-1}}{\sqrt{2(\qnum{4}-1)}}$&
$-\frac{q^{\frac{1}{4}}\sqrt{\qnum{3}}}{2\sqrt{\qnum{4}-1}}$&
$\frac{-q^{-\frac{1}{4}}+q^{\frac{1}{4}}+q^{\frac{3}{4}}}{2\sqrt{\qnum{4}-1}}$&
&
$X$&
$-\frac{q^{\frac{1}{2}}\sqrt{\qnum{2}-1}}{\sqrt{2(\qnum{4}-1)}}$\\
$(2,-1)$ &
&
$\frac{\sqrt{\qnum{2}-1}}{\sqrt{2(\qnum{4}-1)}}$&
&
$\frac{q^{\frac{1}{4}}\sqrt{\qnum{3}}}{2\sqrt{\qnum{4}-1}}$&
$\frac{-q^{-\frac{1}{4}}+q^{\frac{1}{4}}+q^{\frac{3}{4}}}{2\sqrt{\qnum{4}-1}}$&
$X$&
&
$\frac{q^{\frac{1}{2}}\sqrt{\qnum{2}-1}}{\sqrt{2(\qnum{4}-1)}}$\\
$(0,0)_+$ &
$0$&
$\frac{q^{-\frac{1}{4}}\sqrt{\qnum{3}}}{2\sqrt{\qnum{4}-1}}$&
$-\frac{q^{-\frac{1}{4}}\sqrt{\qnum{3}}}{2\sqrt{\qnum{4}-1}}$&
$X$&
$X$&
$\frac{q^{\frac{1}{4}}\sqrt{\qnum{3}}}{2\sqrt{\qnum{4}-1}}$&
$-\frac{q^{\frac{1}{4}}\sqrt{\qnum{3}}}{2\sqrt{\qnum{4}-1}}$&
$0$\\
$(0,0)_-$ &
$-\frac{q^{-\frac{3}{4}}}{\sqrt{\qnum{4}-1}}$&
$\frac{-q^{-\frac{3}{4}}-q^{-\frac{1}{4}}+q^{\frac{1}{4}}}{2\sqrt{\qnum{4}-1}}$&
$\frac{-q^{-\frac{3}{4}}-q^{-\frac{1}{4}}+q^{\frac{1}{4}}}{2\sqrt{\qnum{4}-1}}$&
$X$&
$X$&
$\frac{-q^{-\frac{1}{4}}+q^{\frac{1}{4}}+q^{\frac{3}{4}}}{2\sqrt{\qnum{4}-1}}$&
$\frac{-q^{-\frac{1}{4}}+q^{\frac{1}{4}}+q^{\frac{3}{4}}}{2\sqrt{\qnum{4}-1}}$&
$\frac{q^{\frac{3}{4}}}{\sqrt{\qnum{4}-1}}$\\
$(-2,1)$ &
$-\frac{q^{-\frac{1}{2}}\sqrt{\qnum{2}-1}}{\sqrt{2(\qnum{4}-1)}}$&
&
$X$&
$-\frac{q^{-\frac{1}{4}}\sqrt{\qnum{3}}}{2\sqrt{\qnum{4}-1}}$&
$\frac{-q^{-\frac{3}{4}}-q^{-\frac{1}{4}}+q^{\frac{1}{4}}}{2\sqrt{\qnum{4}-1}}$&
&
$-\frac{\sqrt{\qnum{2}-1}}{\sqrt{2(\qnum{4}-1)}}$&
\\
$(1,-2)$ &
$\frac{q^{-\frac{1}{2}}\sqrt{\qnum{2}-1}}{\sqrt{2(\qnum{4}-1)}}$&
$X$&
&
$\frac{q^{-\frac{1}{4}}\sqrt{\qnum{3}}}{2\sqrt{\qnum{4}-1}}$&
$\frac{-q^{-\frac{3}{4}}-q^{-\frac{1}{4}}+q^{\frac{1}{4}}}{2\sqrt{\qnum{4}-1}}$&
$\frac{\sqrt{\qnum{2}-1}}{\sqrt{2(\qnum{4}-1)}}$&
&
\\
$(-1,-1)$ &
$X$&
$\frac{q^{-\frac{1}{2}}\sqrt{\qnum{2}-1}}{\sqrt{2(\qnum{4}-1)}}$&
$-\frac{q^{-\frac{1}{2}}\sqrt{\qnum{2}-1}}{\sqrt{2(\qnum{4}-1)}}$&
$0$&
$-\frac{q^{-\frac{3}{4}}}{\sqrt{\qnum{4}-1}}$&
&
&
\\
\multicolumn{9}{c}{\phantom{$\sqrt{5}$}}
\\
\multicolumn{9}{c}{The $q$-CG coefficients for the ${\bf 8'}$ (when
$m \neq (0,0)_\pm$) in ${\bf 8}\times{\bf 8}$.}
\end{tabular}
}

\hspace*{-10mm}\rotatebox{90}{%
\begin{tabular}{c|cccccccc}
& $(1,1)$ & $(-1,2)$ & $(2,-1)$ & $(0,0)_+$ & $(0,0)_-$ & $(-2,1)$ & $(1,-2)$ & $(-1,-1)$ \\
\hline
$(1,1)$ &
&
&
&
&
&
&
&
$-\frac{q^{\frac{1}{4}}}{\sqrt{\qnum{4}-1}}$\\
$(-1,2)$ &
&
&
&
&
&
&
$\frac{q^{-\frac{1}{4}}+q^{\frac{1}{4}}-q^{\frac{3}{4}}}{2\sqrt{\qnum{4}-1}}$&
\\
$(2,-1)$ &
&
&
&
&
&
$\frac{q^{-\frac{1}{4}}+q^{\frac{1}{4}}-q^{\frac{3}{4}}}{2\sqrt{\qnum{4}-1}}$&
&
\\
$(0,0)_+$ &
&
&
&
$\frac{-q^{-\frac{3}{4}}+q^{\frac{3}{4}}}{2\sqrt{\qnum{4}-1}}$&
$0$&
&
&
\\
$(0,0)_-$ &
&
&
&
$0$&
$\frac{2(-q^{-\frac{1}{4}}+q^{\frac{1}{4}})-q^{-\frac{3}{4}}+q^{\frac{3}{4}}}{2\sqrt{\qnum{4}-1}}$&
&
&
\\
$(-2,1)$ &
&
&
$\frac{-q^{-\frac{1}{4}}-q^{\frac{1}{4}}+q^{-\frac{3}{4}}}{2\sqrt{\qnum{4}-1}}$&
&
&
&
&
\\
$(1,-2)$ &
&
$\frac{-q^{-\frac{1}{4}}-q^{\frac{1}{4}}+q^{-\frac{3}{4}}}{2\sqrt{\qnum{4}-1}}$&
&
&
&
&
&
\\
$(-1,-1)$ &
$\frac{q^{-\frac{1}{4}}}{\sqrt{\qnum{4}-1}}$ &
&
&
&
&
&
&
\\
\multicolumn{9}{c}{The $q$-CG coefficients for the ${\bf 8'}$ (for $m = (0,0)_{-}$)
in ${\bf 8}\times{\bf 8}$.}
\end{tabular}
}
\hspace*{5mm}
\rotatebox{90}{%
\begin{tabular}{c|cccccccc}
& $(1,1)$ & $(-1,2)$ & $(2,-1)$ & $(0,0)_+$ & $(0,0)_-$ & $(-2,1)$ & $(1,-2)$ & $(-1,-1)$ \\
\hline
$(1,1)$ &
&
&
&
&
&
&
&
$0$\\
$(-1,2)$ &
&
&
&
&
&
&
$-\frac{q^{\frac{1}{4}}\sqrt{\qnum{3}}}{2\sqrt{\qnum{4}-1}}$&
\\
$(2,-1)$ &
&
&
&
&
&
$\frac{q^{\frac{1}{4}}\sqrt{\qnum{3}}}{2\sqrt{\qnum{4}-1}}$&
&
\\
$(0,0)_+$ &
&
&
&
$0$&
$\frac{-q^{-\frac{3}{4}}+q^{\frac{3}{4}}}{2\sqrt{\qnum{4}-1}}$&
&
&
\\
$(0,0)_-$ &
&
&
&
$\frac{-q^{-\frac{3}{4}}+q^{\frac{3}{4}}}{2\sqrt{\qnum{4}-1}}$&
$0$&
&
&
\\
$(-2,1)$ &
&
&
$-\frac{q^{-\frac{1}{4}}\sqrt{\qnum{3}}}{2\sqrt{\qnum{4}-1}}$&
&
&
&
&
\\
$(1,-2)$ &
&
$\frac{q^{-\frac{1}{4}}\sqrt{\qnum{3}}}{2\sqrt{\qnum{4}-1}}$&
&
&
&
&
&
\\
$(-1,-1)$ &
$0$&
&
&
&
&
&
&
\\
\multicolumn{9}{c}{The $q$-CG coefficients for the ${\bf 8'}$ (for $m = (0,0)_+$) in
${\bf 8}\times{\bf 8}$.}
\end{tabular}
}

\rotatebox{90}{%
\begin{tabular}{c|cccccccc}
& $(1,1)$ & $(-1,2)$ & $(2,-1)$ & $(0,0)_+$ & $(0,0)_-$ & $(-2,1)$ & $(1,-2)$ & $(-1,-1)$ \\
\hline
$(1,1)$ &
&
&
$\frac{q^{\frac{1}{4}}}{\sqrt{\qnum{2}}}$&
$\frac{\sqrt{\qnum{2}+1}}{\sqrt{2\qnum{2}\qnum{3}}}$&
$\frac{\sqrt{\qnum{2}-1}}{\sqrt{2\qnum{2}\qnum{3}}}$&
$\frac{q^{-\frac{1}{4}}}{\sqrt{\qnum{2}\qnum{3}}}$&
$\frac{q^{-\frac{1}{4}}}{\sqrt{\qnum{2}\qnum{3}}}$&
$\frac{q^{-\frac{1}{2}}}{\qnum{2}\sqrt{\qnum{3}}}$\\
$(-1,2)$ &
&
&
$\frac{q^{\frac{3}{4}}}{\sqrt{\qnum{2}\qnum{3}}}$&
$\frac{q^{\frac{1}{2}}\sqrt{\qnum{2}+1}}{\sqrt{2\qnum{2}\qnum{3}}}$&
$\frac{q^{\frac{1}{2}}\sqrt{\qnum{2}-1}}{\sqrt{2\qnum{2}\qnum{3}}}$&
$\frac{q^{\frac{1}{4}}}{\sqrt{\qnum{2}}}$&
$\frac{1}{\qnum{2}\sqrt{\qnum{3}}}$&
$\frac{q^{-\frac{1}{4}}}{\sqrt{\qnum{2}\qnum{3}}}$\\
$(2,-1)$ &
$-\frac{q^{-\frac{1}{4}}}{\sqrt{\qnum{2}}}$&
$-\frac{q^{-\frac{3}{4}}}{\sqrt{\qnum{2}\qnum{3}}}$&
&
$\frac{(1-q^{-\frac{1}{2}})\sqrt{\qnum{2}+1}}{\sqrt{2\qnum{2}\qnum{3}}}$&
$\frac{(1+q^{-\frac{1}{2}})\sqrt{\qnum{2}-1}}{\sqrt{2\qnum{2}\qnum{3}}}$&
$\frac{1}{\qnum{2}\sqrt{\qnum{3}}}$&
&
$\frac{q^{-\frac{1}{4}}}{\sqrt{\qnum{2}\qnum{3}}}$\\
$(0,0)_+$ &
$-\frac{\sqrt{\qnum{2}+1}}{\sqrt{2\qnum{2}\qnum{3}}}$&
$-\frac{q^{-\frac{1}{2}}\sqrt{\qnum{2}+1}}{\sqrt{2\qnum{2}\qnum{3}}}$&
$\frac{(-1+q^{\frac{1}{2}})\sqrt{\qnum{2}+1}}{\sqrt{2\qnum{2}\qnum{3}}}$&
$\frac{(q^{\frac{1}{2}}-q^{-\frac{1}{2}})(\qnum{2}+1)}{2\qnum{2}\sqrt{\qnum{3}}}$&
$\frac{1}{2}$&
$\frac{q^{\frac{1}{2}}\sqrt{\qnum{2}+1}}{\sqrt{2\qnum{2}\qnum{3}}}$&
$\frac{(1-q^{-\frac{1}{2}})\sqrt{\qnum{2}+1}}{\sqrt{2\qnum{2}\qnum{3}}}$&
$\frac{\sqrt{\qnum{2}+1}}{\sqrt{2\qnum{2}\qnum{3}}}$\\
$(0,0)_-$ &
$-\frac{\sqrt{\qnum{2}-1}}{\sqrt{2\qnum{2}\qnum{3}}}$&
$-\frac{q^{-\frac{1}{2}}\sqrt{\qnum{2}-1}}{\sqrt{2\qnum{2}\qnum{3}}}$&
$\frac{(-1-q^{\frac{1}{2}})\sqrt{\qnum{2}-1}}{\sqrt{2\qnum{2}\qnum{3}}}$&
$-\frac{1}{2}$&
$\frac{(q^{-\frac{1}{2}}-q^{\frac{1}{2}})(\qnum{2}-1)}{2\qnum{2}\sqrt{\qnum{3}}}$&
$-\frac{q^{\frac{1}{2}}\sqrt{\qnum{2}-1}}{\sqrt{2\qnum{2}\qnum{3}}}$&
$\frac{(-1-q^{-\frac{1}{2}})\sqrt{\qnum{2}-1}}{\sqrt{2\qnum{2}\qnum{3}}}$&
$-\frac{\sqrt{\qnum{2}-1}}{\sqrt{2\qnum{2}\qnum{3}}}$\\
$(-2,1)$ &
$-\frac{q^{\frac{1}{4}}}{\sqrt{\qnum{2}\qnum{3}}}$&
$-\frac{q^{-\frac{1}{4}}}{\sqrt{\qnum{2}}}$&
$-\frac{1}{\qnum{2}\sqrt{\qnum{3}}}$&
$-\frac{q^{-\frac{1}{2}}\sqrt{\qnum{2}+1}}{\sqrt{2\qnum{2}\qnum{3}}}$&
$\frac{q^{-\frac{1}{2}}\sqrt{\qnum{2}-1}}{\sqrt{2\qnum{2}\qnum{3}}}$&
&
$-\frac{q^{-\frac{3}{4}}}{\sqrt{\qnum{2}\qnum{3}}}$&
\\
$(1,-2)$ &
$-\frac{q^{\frac{1}{4}}}{\sqrt{\qnum{2}\qnum{3}}}$&
$-\frac{1}{\qnum{2}\sqrt{\qnum{3}}}$&
&
$\frac{(-1+q^{\frac{1}{2}})\sqrt{\qnum{2}+1}}{\sqrt{2\qnum{2}\qnum{3}}}$&
$\frac{(1+q^{\frac{1}{2}})\sqrt{\qnum{2}-1}}{\sqrt{2\qnum{2}\qnum{3}}}$&
$\frac{q^{\frac{3}{4}}}{\sqrt{\qnum{2}\qnum{3}}}$&
&
$\frac{q^{\frac{1}{4}}}{\sqrt{\qnum{2}}}$\\
$(-1,-1)$ &
$-\frac{q^{\frac{1}{2}}}{\qnum{2}\sqrt{\qnum{3}}}$&
$-\frac{q^{\frac{1}{4}}}{\sqrt{\qnum{2}\qnum{3}}}$&
$-\frac{q^{\frac{1}{4}}}{\sqrt{\qnum{2}\qnum{3}}}$&
$-\frac{\sqrt{\qnum{2}+1}}{\sqrt{2\qnum{2}\qnum{3}}}$&
$\frac{\sqrt{\qnum{2}-1}}{\sqrt{2\qnum{2}\qnum{3}}}$&
&
$-\frac{q^{-\frac{1}{4}}}{\sqrt{\qnum{2}}}$&
\\
\multicolumn{9}{c}{\phantom{$\sqrt{5}$}}
\\
\multicolumn{9}{c}{The $q$-CG coefficients for the ${\bf 10}$ in
${\bf 8}\times{\bf 8}$.}
\end{tabular}
}

\subsubsection{${\bf 10}\times {\bf 8} = {\bf 8}$}
%For layout reasons, we give the Clebsch-Gordan coefficients for the ${\bf 8}$ in
%${\bf 10}\times {\bf 8}$, instead of ${\bf 8}\times {\bf 10}$. The coefficients in the
%later case are obtained from the former via the usual symmetry relations. 

\rotatebox{90}{%
\begin{tabular}{c|cccccccc}
& $(1,1)$ & $(-1,2)$ & $(2,-1)$ & $(0,0)_+$ & $(0,0)_-$ & $(-2,1)$ & $(1,-2)$ & $(-1,-1)$ \\
\hline
$(3,0)$ &
&
&
&
&
&
$\frac{q^{\frac{3}{4}}\sqrt{\qnum{2}}}{\sqrt{\qnum{5}}}$&
&
$\frac{q^{\frac{3}{4}}\sqrt{\qnum{2}}}{\sqrt{\qnum{5}}}$
\\
$(1,1)$ &
&
&
&
$-\frac{\sqrt{\qnum{2}(\qnum{2}+1)}}{\sqrt{2\qnum{3}\qnum{5}}}$&
$-\frac{\sqrt{\qnum{2}(\qnum{2}-1)}}{\sqrt{2\qnum{3}\qnum{5}}}$&
$\frac{q^{\frac{5}{4}}\sqrt{\qnum{2}}}{\sqrt{\qnum{3}\qnum{5}}}$&
$-\frac{q^{-\frac{1}{4}}\sqrt{\qnum{2}}}{\sqrt{\qnum{3}\qnum{5}}}$&
$X$
\\
$(-1,2)$ &
&
&
$\frac{q^{-\frac{3}{4}}\sqrt{\qnum{2}}}{\sqrt{\qnum{3}\qnum{5}}}$&
$-\frac{q^{\frac{1}{2}}\sqrt{\qnum{2}(\qnum{2}+1)}}{\sqrt{2\qnum{3}\qnum{5}}}$&
$-\frac{q^{\frac{1}{2}}\sqrt{\qnum{2}(\qnum{2}-1)}}{\sqrt{2\qnum{3}\qnum{5}}}$&
&
$X$&
$\frac{q^{\frac{5}{4}}\sqrt{\qnum{2}}}{\sqrt{\qnum{3}\qnum{5}}}$
\\
$(2,-1)$ &
&
$\frac{q^{-\frac{3}{4}}\sqrt{\qnum{2}}}{\sqrt{\qnum{3}\qnum{5}}}$&
&
$\frac{(-1+q^{-\frac{1}{2}})\sqrt{\qnum{2}(\qnum{2}+1)}}{\sqrt{2\qnum{3}\qnum{5}}}$&
$\frac{(-1-q^{-\frac{1}{2}})\sqrt{\qnum{2}(\qnum{2}-1)}}{\sqrt{2\qnum{3}\qnum{5}}}$&
$X$&
&
$\frac{q^{\frac{5}{4}}\sqrt{\qnum{2}}}{\sqrt{\qnum{3}\qnum{5}}}$
\\
$(-3,3)$ &
&
&
$\frac{q^{-\frac{1}{4}}\sqrt{\qnum{2}}}{\sqrt{\qnum{5}}}$&
&
&
&
-$\frac{q^{\frac{1}{4}}\sqrt{\qnum{2}}}{\sqrt{\qnum{5}}}$&
\\
$(0,0)$ &
$-\frac{q^{-\frac{3}{2}}}{\sqrt{\qnum{3}\qnum{5}}}$&
$\frac{q^{-\frac{1}{2}}}{\sqrt{\qnum{3}\qnum{5}}}$&
$\frac{q^{-\frac{1}{2}}}{\sqrt{\qnum{3}\qnum{5}}}$&
$X$&
$X$&
$-\frac{q^{\frac{1}{2}}}{\sqrt{\qnum{3}\qnum{5}}}$&
$-\frac{q^{\frac{1}{2}}}{\sqrt{\qnum{3}\qnum{5}}}$&
$\frac{q^{\frac{3}{2}}}{\sqrt{\qnum{3}\qnum{5}}}$
\\
$(-2,1)$ &
$-\frac{q^{-\frac{5}{4}}\sqrt{\qnum{2}}}{\sqrt{\qnum{3}\qnum{5}}}$&
&
$X$&
$\frac{q^{-\frac{1}{2}}\sqrt{\qnum{2}(\qnum{2}+1)}}{\sqrt{2\qnum{3}\qnum{5}}}$&
$-\frac{q^{-\frac{1}{2}}\sqrt{\qnum{2}(\qnum{2}-1)}}{\sqrt{2\qnum{3}\qnum{5}}}$&
&
$-\frac{q^{\frac{3}{4}}\sqrt{\qnum{2}}}{\sqrt{\qnum{3}\qnum{5}}}$&
\\
$(1,-2)$ &
$-\frac{q^{-\frac{5}{4}}\sqrt{\qnum{2}}}{\sqrt{\qnum{3}\qnum{5}}}$&
$X$&
&
$\frac{(1-q^{\frac{1}{2}})\sqrt{\qnum{2}(\qnum{2}+1)}}{\sqrt{2\qnum{3}\qnum{5}}}$&
$\frac{(-1-q^{\frac{1}{2}})\sqrt{\qnum{2}(\qnum{2}-1)}}{\sqrt{2\qnum{3}\qnum{5}}}$&
$-\frac{q^{\frac{3}{4}}\sqrt{\qnum{2}}}{\sqrt{\qnum{3}\qnum{5}}}$&
&
\\
$(-1,-1)$ &
$X$&
$-\frac{q^{-\frac{5}{4}}\sqrt{\qnum{2}}}{\sqrt{\qnum{3}\qnum{5}}}$&
$\frac{q^{\frac{1}{4}}\sqrt{\qnum{2}}}{\sqrt{\qnum{3}\qnum{5}}}$&
$\frac{\sqrt{\qnum{2}(\qnum{2}+1)}}{\sqrt{2\qnum{3}\qnum{5}}}$&
$-\frac{\sqrt{\qnum{2}(\qnum{2}-1)}}{\sqrt{2\qnum{3}\qnum{5}}}$&
&
&
\\
$(-3,0)$ &
$-\frac{q^{-\frac{3}{4}}\sqrt{\qnum{2}}}{\sqrt{\qnum{5}}}$&
$-\frac{q^{-\frac{3}{4}}\sqrt{\qnum{2}}}{\sqrt{\qnum{5}}}$&
&
&
&
&
&
\\
\multicolumn{9}{c}{\phantom{$\sqrt{5}$}}
\\
\multicolumn{9}{c}{The $q$-CG coefficients for the ${\bf 8}$ (for $m \neq (0,0)_{\pm}$)
in ${\bf 10}\times{\bf 8}$.}
\end{tabular}
}

\hspace{-2cm}
\rotatebox{90}{%
\begin{tabular}{c|cccccccc}
& $(1,1)$ & $(-1,2)$ & $(2,-1)$ & $(0,0)_+$ & $(0,0)_-$ & $(-2,1)$ & $(1,-2)$ & $(-1,-1)$ \\
\hline
$(3,0)$ &
&
&
&
&
&
&
&
\\
$(1,1)$ &
&
&
&
&
&
&
&
$\frac{q\sqrt{\qnum{2}(\qnum{2}+1)}}{\sqrt{2\qnum{3}\qnum{5}}}$
\\
$(-1,2)$ &
&
&
&
&
&
&
$-\frac{\sqrt{\qnum{2}(\qnum{2}+1)}}{\sqrt{2\qnum{3}\qnum{5}}}$&
\\
$(2,-1)$ &
&
&
&
&
&
$\frac{(-1+q^{\frac{3}{2}})\sqrt{\qnum{2}}}{\sqrt{2\qnum{3}\qnum{5}(\qnum{2}+1)}}$&
&
\\
$(-3,3)$ &
&
&
&
&
&
&
&
\\
$(0,0)$ &
&
&
&
$\frac{q^{-1}+q^{-\frac{1}{2}}-q^{\frac{1}{2}}-q}{2\sqrt{\qnum{3}\qnum{5}}}$&
$-\frac{\qnum{2}}{2\sqrt{\qnum{5}}}$&
&
&
\\
$(-2,1)$ &
&
&
$\frac{\sqrt{\qnum{2}(\qnum{2}+1)}}{\sqrt{2\qnum{3}\qnum{5}}}$&
&
&
&
&
\\
$(1,-2)$ &
&
$\frac{(1-q^{-\frac{3}{2}})\sqrt{\qnum{2}}}{\sqrt{2\qnum{3}\qnum{5}(\qnum{2}+1)}}$&
&
&
&
&
&
\\
$(-1,-1)$ &
$-\frac{q^{-1}\sqrt{\qnum{2}(\qnum{2}+1)}}{\sqrt{2\qnum{3}\qnum{5}}}$&
&
&
&
&
&
&
\\
$(-3,0)$ &
&
&
&
&
&
&
&
\\
\multicolumn{9}{c}{The $q$-CG coefficients for the ${\bf 8}$ (for $m = (0,0)_{+}$)
in ${\bf 10}\times{\bf 8}$.}
\end{tabular}
}
\hspace{.5cm}
\rotatebox{90}{%
\begin{tabular}{c|cccccccc}
& $(1,1)$ & $(-1,2)$ & $(2,-1)$ & $(0,0)_+$ & $(0,0)_-$ & $(-2,1)$ & $(1,-2)$ & $(-1,-1)$ \\
\hline
$(3,0)$ &
&
&
&
&
&
&
&
\\
$(1,1)$ &
&
&
&
&
&
&
&
$\frac{q\sqrt{\qnum{2}(\qnum{2}-1)}}{\sqrt{2\qnum{3}\qnum{5}}}$
\\
$(-1,2)$ &
&
&
&
&
&
&
$-\frac{\sqrt{\qnum{2}(\qnum{2}-1)}}{\sqrt{2\qnum{3}\qnum{5}}}$&
\\
$(2,-1)$ &
&
&
&
&
&
$\frac{(-1-q^{\frac{3}{2}})\sqrt{\qnum{2}}}{\sqrt{2\qnum{3}\qnum{5}(\qnum{2}-1)}}$&
&
\\
$(-3,3)$ &
&
&
&
&
&
&
&
\\
$(0,0)$ &
&
&
&
$\frac{\qnum{2}}{2\sqrt{\qnum{5}}}$&
$\frac{-q^{-1}+q^{-\frac{1}{2}}-q^{\frac{1}{2}}+q}{2\sqrt{\qnum{3}\qnum{5}}}$&
&
&
\\
$(-2,1)$ &
&
&
$-\frac{\sqrt{\qnum{2}(\qnum{2}-1)}}{\sqrt{2\qnum{3}\qnum{5}}}$&
&
&
&
&
\\
$(1,-2)$ &
&
$\frac{(-1-q^{-\frac{3}{2}})\sqrt{\qnum{2}}}{\sqrt{2\qnum{3}\qnum{5}(\qnum{2}-1)}}$&
&
&
&
&
&
\\
$(-1,-1)$ &
$\frac{q^{-1}\sqrt{\qnum{2}(\qnum{2}-1)}}{\sqrt{2\qnum{3}\qnum{5}}}$&
&
&
&
&
&
&
\\
$(-3,0)$ &
&
&
&
&
&
&
&
\\
\multicolumn{9}{c}{The $q$-CG coefficients for the ${\bf 8}$ (for $m = (0,0)_{-}$)
in ${\bf 10}\times{\bf 8}$.}
\end{tabular}
}

\subsubsection{${\bf 10}\times{\bf 10} = \overline{\bf 10}$}

\rotatebox{90}{%
\begin{tabular}{c|cccccccccc}
& $(3,0)$ & $(1,1)$ & $(-1,2)$ & $(2,-1)$ & $(-3,3)$ & $(0,0)$ & $(-2,1)$ & $(1,-2)$ & $(-1,-1)$ & $(-3,0)$\\
\hline
$(3,0)$ &
&
&
&
&
$\frac{q^{\frac{3}{4}}}{\sqrt{\qnum{4}}}$&
&
$\frac{q^{\frac{3}{4}}}{\sqrt{\qnum{4}}}$&
&
$\frac{q^{\frac{3}{4}}}{\sqrt{\qnum{4}}}$&
$\frac{q^{\frac{3}{4}}}{\sqrt{\qnum{4}}}$\\
$(1,1)$ &
&
&
$\frac{-q^{\frac{1}{4}}}{\sqrt{\qnum{4}}}$&
&
&
$\frac{-\sqrt{\qnum{2}}}{\sqrt{\qnum{3}\qnum{4}}}$&
$\frac{q^{\frac{5}{4}}}{\sqrt{\qnum{3}\qnum{4}}}$&
$\frac{-q^{-1}}{\sqrt{\qnum{3}\qnum{4}}}$&
$\frac{q\sqrt{\qnum{2}}}{\sqrt{\qnum{3}\qnum{4}}}$&
$\frac{q^{\frac{3}{4}}}{\sqrt{\qnum{4}}}$\\
$(-1,2)$ &
&
$\frac{q^{-\frac{1}{4}}}{\sqrt{\qnum{4}}}$&
&
$\frac{q^{-\frac{3}{4}}}{\sqrt{\qnum{3}\qnum{4}}}$&
&
$\frac{-q^{\frac{1}{2}}\sqrt{\qnum{2}}}{\sqrt{\qnum{3}\qnum{4}}}$&
&
$\frac{-\sqrt{\qnum{2}}}{\sqrt{\qnum{3}\qnum{4}}}$&
$\frac{q^{\frac{5}{4}}}{\sqrt{\qnum{3}\qnum{4}}}$&
$\frac{q^{\frac{3}{4}}}{\sqrt{\qnum{4}}}$\\
$(2,-1)$ &
&
&
$\frac{-q^{\frac{3}{4}}}{\sqrt{\qnum{3}\qnum{4}}}$&
&
$\frac{-q^{\frac{1}{4}}}{\sqrt{\qnum{4}}}$&
$\frac{-q^{\frac{1}{2}}\sqrt{\qnum{2}}}{\sqrt{\qnum{3}\qnum{4}}}$&
$\frac{-\sqrt{\qnum{2}}}{\sqrt{\qnum{3}\qnum{4}}}$&
$\frac{-q^{\frac{1}{4}}}{\sqrt{\qnum{4}}}$&
$\frac{-q^{-\frac{1}{4}}}{\sqrt{\qnum{3}\qnum{4}}}$&
\\
$(-3,3)$ &
$\frac{-q^{-\frac{3}{4}}}{\sqrt{\qnum{4}}}$&
&
&
$\frac{q^{-\frac{1}{4}}}{\sqrt{\qnum{4}}}$&
&
&
&
$\frac{-q^{\frac{1}{4}}}{\sqrt{\qnum{4}}}$&
&
$\frac{q^{\frac{3}{4}}}{\sqrt{\qnum{4}}}$\\
$(0,0)$ &
&
$\frac{\sqrt{\qnum{2}}}{\sqrt{\qnum{3}\qnum{4}}}$&
$\frac{q^{-\frac{1}{2}}\sqrt{\qnum{2}}}{\sqrt{\qnum{3}\qnum{4}}}$&
$\frac{q^{-\frac{1}{2}}\sqrt{\qnum{2}}}{\sqrt{\qnum{3}\qnum{4}}}$&
&
$\frac{(q^{-\frac{1}{2}}-q^{\frac{1}{2}})\sqrt{\qnum{2}}}{\sqrt{\qnum{3}\qnum{4}}}$&
$\frac{-q^{\frac{1}{2}}\sqrt{\qnum{2}}}{\sqrt{\qnum{3}\qnum{4}}}$&
$\frac{-q^{\frac{1}{2}}\sqrt{\qnum{2}}}{\sqrt{\qnum{3}\qnum{4}}}$&
$\frac{-\sqrt{\qnum{2}}}{\sqrt{\qnum{3}\qnum{4}}}$&
\\
$(-2,1)$ &
$\frac{-q^{-\frac{3}{4}}}{\sqrt{\qnum{4}}}$&
$\frac{-q^{-\frac{5}{4}}}{\sqrt{\qnum{3}\qnum{4}}}$&
&
$\frac{\sqrt{\qnum{2}}}{\sqrt{\qnum{3}\qnum{4}}}$&
&
$\frac{q^{-\frac{1}{2}}\sqrt{\qnum{2}}}{\sqrt{\qnum{3}\qnum{4}}}$&
&
$\frac{-q^{\frac{3}{4}}}{\sqrt{\qnum{3}\qnum{4}}}$&
$\frac{-q^{\frac{1}{4}}}{\sqrt{\qnum{4}}}$&
\\
$(1,-2)$ &
&
$\frac{q^{\frac{1}{4}}}{\sqrt{\qnum{3}\qnum{4}}}$&
$\frac{\sqrt{\qnum{2}}}{\sqrt{\qnum{3}\qnum{4}}}$&
$\frac{q^{-\frac{1}{4}}}{\sqrt{\qnum{4}}}$&
$\frac{q^{-\frac{1}{4}}}{\sqrt{\qnum{4}}}$&
$\frac{q^{-\frac{1}{2}}\sqrt{\qnum{2}}}{\sqrt{\qnum{3}\qnum{4}}}$&
$\frac{q^{-\frac{3}{4}}}{\sqrt{\qnum{3}\qnum{4}}}$&
&
&
\\
$(-1,-1)$ &
$\frac{-q^{-\frac{3}{4}}}{\sqrt{\qnum{4}}}$&
$\frac{-q^{-1}\sqrt{\qnum{2}}}{\sqrt{\qnum{3}\qnum{4}}}$&
$\frac{-q^{-\frac{5}{4}}}{\sqrt{\qnum{3}\qnum{4}}}$&
$\frac{q^{\frac{1}{4}}}{\sqrt{\qnum{3}\qnum{4}}}$&
&
$\frac{\sqrt{\qnum{2}}}{\sqrt{\qnum{3}\qnum{4}}}$&
$\frac{q^{-\frac{1}{4}}}{\sqrt{\qnum{4}}}$&
&
&
\\
$(-3,0)$ &
$\frac{-q^{-\frac{3}{4}}}{\sqrt{\qnum{4}}}$&
$\frac{-q^{-\frac{3}{4}}}{\sqrt{\qnum{4}}}$&
$\frac{-q^{-\frac{3}{4}}}{\sqrt{\qnum{4}}}$&
&
$\frac{-q^{-\frac{3}{4}}}{\sqrt{\qnum{4}}}$&
&
&
&
&
\\
\multicolumn{11}{c}{\phantom{$\sqrt{5}$}}
\\
\multicolumn{11}{c}{The $q$-CG coefficients for the $\overline{\bf 10}$ in ${\bf 10}\times{\bf 10}$.}
\end{tabular}
}

\subsubsection{${\bf 10}\times\overline{\bf 10} = \id$}
%\newpage
%\rotatebox{90}{%
\hspace{-1cm}
\begin{tabular}{c|cccccccccc}
& $(0,3)$ & $(1,1)$ & $(-1,2)$ & $(2,-1)$ & $(0,0)$ & $(3,-3)$ & $(-2,1)$ & $(1,-2)$ & $(-1,-1)$ & $(0,-3)$\\
\hline
$(3,0)$ &
&
&
&
&
&
&
&
&
&
$q^{\frac{3}{2}}$
\\
$(1,1)$ &
&
&
&
&
&
&
&
&
$-q$&
\\
$(-1,2)$ &
&
&
&
&
&
&
&
$q^{\frac{1}{2}}$&
&
\\
$(2,-1)$ &
&
&
&
&
&
&
$q^{\frac{1}{2}}$&
&
&
\\
$(-3,3)$ &
&
&
&
&
&
$-1$&
&
&
&
\\
$(0,0)$ &
&
&
&
&
$-1$&
&
&
&
&
\\
$(-2,1)$ &
&
&
&
$q^{-\frac{1}{2}}$&
&
&
&
&
&
\\
$(1,-2)$ &
&
&
$q^{-\frac{1}{2}}$&
&
&
&
&
&
&
\\
$(-1,-1)$ &
&
$-q^{-1}$&
&
&
&
&
&
&
&
\\
$(-3,0)$ &
$q^{-\frac{3}{2}}$&
&
&
&
&
&
&
&
&
\\
\multicolumn{9}{c}{\phantom{5}}
\\
\multicolumn{11}{c}{The $q$-CG coefficients for the $\id$ in ${\bf 10}\times\overline{\bf 10}$. Each symbol
has an additional factor $\frac{1}{\sqrt{\qnum{5}(\qnum{3}-1)}}$}
\end{tabular}
%}

\subsection{The relation between the $q$-CG coefficients}
\label{su33abc}

In table \ref{su33sym}, we specify the coefficients $s_1$, $s_2$ and $s_3$, which appear in
the symmetry relations
between the various $q$-CG coefficients as explained in section \ref{symandconv}.

\begin{table}[ht]
\begin{center}
\begin{tabular}{|c|c|c|r|r|r|}
\hline
$j_1$ & $j_2$ & $j$ & $s_1$ & $s_2$ & $s_3$ \\
\hline
\hline
${\bf 8}$ & ${\bf 8}$ & $\id$ & $0$ & $0$ & $0$ \\
${\bf 8}$ & ${\bf 8}$ & ${\bf 8}$ & $0$ & $0$ & $0$ \\
${\bf 8}$ & ${\bf 8}$ & ${\bf 8'}$ & $1$ & $0$ & $1$ \\
${\bf 8}$ & ${\bf 8}$ & ${\bf 10}$ & $1$ & $1$ & $0$ \\
${\bf 8}$ & ${\bf 10}$ & ${\bf 8}$ & $1$ & $1$ & $0$ \\
${\bf 10}$ & ${\bf 10}$ & ${\bf \overline{10}}$ & $1$ & $1$ & $0$ \\
${\bf 10}$ & $\overline{{\bf 10}}$ & $\id$ & $0$ & $0$ & $0$ \\
\hline
\end{tabular}
\end{center}
\caption{The parameters $s_1$, $s_2$ and $s_3$ in the symmetry relations between the $q$-CG
coefficients for $su(3)_3/Z_3$.}
\label{su33sym}
\end{table}

\subsection{The $F$ and $R$-symbols for $su(3)_3/Z_3$}

Because the expressions are rather involved, we will not give the $F$-symbols for
$su(3)_3/Z_3$ as a function of $q$, but only for the specialization $q=e^{2 \pi i /6}$. In principle,
it is straightforward to obtain them as a function of $q$ from the $q$-Clebsch-Gordan coefficients.

All the one-dimensional $F$-symbols with a $\id$ on an outer line are one. However, in case of the presence of a vertex with three ${\bf 8}$'s, one has two-dimensional objects: 
\begin{equation}
\fs{\id}{\bf 8}{\bf 8}{\bf 8}=
\fs{\bf 8}{\id}{\bf 8}{\bf 8}=
\fs{\bf 8}{\bf 8}{\id}{\bf 8}=
\fs{\bf 8}{\bf 8}{\bf 8}{\id}=
\begin{pmatrix} 1 & 0 \\ 0 & 1 \end{pmatrix}
\end{equation}
The following non-tivial symbols are $1$:
\begin{align}
&\fs{\bf 8}{\bf 10}{\bf 8}{\bf 10}=
\fs{\bf 8}{\bf 10}{\bf 8}{\bf \overline{10}}=
\fs{\bf 8}{\bf \overline{10}}{\bf 8}{\bf 10}=
\fs{\bf 8}{\bf \overline{10}}{\bf 8}{\bf \overline{10}}=
\fs{\bf 10}{\bf 8}{\bf 10}{\bf 8}=
\fs{\bf 10}{\bf 8}{\bf \overline{10}}{\bf 8}=
\fs{\bf \overline{10}}{\bf 8}{\bf 10}{\bf 8}=
\fs{\bf \overline{10}}{\bf 8}{\bf \overline{10}}{\bf 8}=\\ \nonumber
&\fs{\bf 10}{\bf \overline{10}}{\bf 10}{\bf 10}=
\fs{\bf \overline{10}}{\bf 10}{\bf \overline{10}}{\bf \overline{10}}= 1
\end{align}
The following symbols are $-1$.
\begin{align}
&\fs{\bf 8}{\bf 8}{\bf 10}{\bf 10}=
\fs{\bf 8}{\bf 8}{\bf 10}{\bf \overline{10}}=
\fs{\bf 8}{\bf 8}{\bf \overline{10}}{\bf 10}=
\fs{\bf 8}{\bf 8}{\bf \overline{10}}{\bf \overline{10}}=
\fs{\bf 8}{\bf 10}{\bf 10}{\bf 8}=
\fs{\bf 8}{\bf 10}{\bf \overline{10}}{\bf 8}=
\fs{\bf 8}{\bf \overline{10}}{\bf 10}{\bf 8}=
\fs{\bf 8}{\bf \overline{10}}{\bf \overline{10}}{\bf 8}=\\ \nonumber
&\fs{\bf 10}{\bf 10}{\bf 8}{\bf 8}=
\fs{\bf 10}{\bf \overline{10}}{\bf 8}{\bf 8}=
\fs{\bf \overline{10}}{\bf 10}{\bf 8}{\bf 8}=
\fs{\bf \overline{10}}{\bf \overline{10}}{\bf 8}{\bf 8}=
\fs{\bf 10}{\bf 8}{\bf 8}{\bf 10}=
\fs{\bf 10}{\bf 8}{\bf 8}{\bf \overline{10}}=
\fs{\bf \overline{10}}{\bf 8}{\bf 8}{\bf 10}=
\fs{\bf \overline{10}}{\bf 8}{\bf 8}{\bf \overline{10}}=\\ \nonumber
&\fs{\bf 10}{\bf 10}{\bf \overline{10}}{\bf 10}=
\fs{\bf 10}{\bf \overline{10}}{\bf \overline{10}}{\bf \overline{10}}=
\fs{\bf \overline{10}}{\bf 10}{\bf 10}{\bf 10}=
\fs{\bf \overline{10}}{\bf \overline{10}}{\bf 10}{\bf \overline{10}}= -1
\end{align}
In addition, we have
\begin{equation}
\fs{\bf 8}{\bf 8}{\bf 8}{\bf 10}=
\fs{\bf 8}{\bf 8}{\bf 10}{\bf 8}=
\fs{\bf 8}{\bf \overline{10}}{\bf 8}{\bf 8}=
\fs{\bf 10}{\bf 8}{\bf 8}{\bf 8}=
\begin{pmatrix} -\frac{1}{2} & -\frac{\sqrt{3}}{2} \\ \frac{\sqrt{3}}{2} & -\frac{1}{2} \end{pmatrix} 
\end{equation}

\begin{equation}
\fs{\bf 8}{\bf 8}{\bf 8}{\bf \overline{10}}=
\fs{\bf 8}{\bf 8}{\bf \overline{10}}{\bf 8}=
\fs{\bf 8}{\bf 10}{\bf 8}{\bf 8}=
\fs{\bf \overline{10}}{\bf 8}{\bf 8}{\bf 8}=
\begin{pmatrix} -\frac{1}{2} & \frac{\sqrt{3}}{2} \\ -\frac{\sqrt{3}}{2} & -\frac{1}{2} \end{pmatrix} 
\end{equation}
Finally, the most interesting $F$-symbol reads
\begin{equation}
\fs{\bf 8}{\bf 8}{\bf 8}{\bf 8} =
\begin{pmatrix}
\frac{1}{3} & \frac{1}{\sqrt{3}} & 0 & 0 & \frac{1}{\sqrt{3}} & -\frac{1}{3} & -\frac{1}{3} \\
\frac{1}{\sqrt{3}} & -\frac{1}{2} & 0 & 0 & \frac{1}{2} & \frac{1}{\sqrt{12}} & \frac{1}{\sqrt{12}} \\
0 & 0 & \frac{1}{2} & \frac{1}{2} & 0 & \frac{1}{2} & -\frac{1}{2} \\
0 & 0 & \frac{1}{2} & \frac{1}{2} & 0 & -\frac{1}{2} & \frac{1}{2} \\
\frac{1}{\sqrt{3}} & \frac{1}{2} & 0 & 0 & -\frac{1}{2} & \frac{1}{\sqrt{12}} & \frac{1}{\sqrt{12}} \\
-\frac{1}{3} & \frac{1}{\sqrt{12}} & -\frac{1}{2} & \frac{1}{2} & \frac{1}{\sqrt{12}} & \frac{1}{3} & \frac{1}{3} \\
-\frac{1}{3} & \frac{1}{\sqrt{12}} & \frac{1}{2} & -\frac{1}{2} & \frac{1}{\sqrt{12}} & \frac{1}{3} & \frac{1}{3} \\
\end{pmatrix}
\end{equation}
Here, we used the following basis. The first row corresponds to the case for which $j_{12}= \id$.
The next four rows correspond to $j_{12}= {\bf 8}$. In this case, there are two vertices with three external
${\bf 8}$ lines, which each are two-dimensional. The second and fifth row correspond to the cases in
which we took the vertices to be `the same', while the third and fourth row correspond to the `off-diagonal'
cases. Finally, row six and seven correspond to the ${\bf 10}$ and ${\bf \overline{10}}$ respectively.
Note that this matrix is not symmetric, but
${\fs{\bf 8}{\bf 8}{\bf 8}{\bf 8} \cdot \fs{\bf 8}{\bf 8}{\bf 8}{\bf 8}}^T = \id$, as it should.
The $R$ symbols read
\begin{align}
R^{\id,x}_x &= 1 \\
R^{{\bf 8},{\bf 8}}_{\id} &= q^{-3} \rightarrow -1 &
R^{{\bf 8},{\bf 8}}_{{\bf 8},1} &= q^{-\frac{3}{2}} \rightarrow -i & 
R^{{\bf 8},{\bf 8}}_{{\bf 8},2} &= -q^{-\frac{3}{2}} \rightarrow i &
R^{{\bf 8},{\bf 8}}_{\bf 10} &= -1 \\
R^{{\bf 8},{\bf 10}}_{{\bf 8}} &= -q^{-3} \rightarrow 1 \\
R^{{\bf 10},{\bf 10}}_{{\bf \overline{10}}} &= -q^{-3} \rightarrow 1 \\
R^{{\bf 10},{\bf \overline{10}}}_{\id} &= q^{-6} \rightarrow 1 
\end{align}
We verified that the symbols above satisfy the pentagon and both hexagon equations.

Because of the relevance for the quantum dimension and Frobenius-Schur indicator, we
give the following $F$-symbols as a function of $q$:
\begin{align}
\bigl(F^{{\bf 8},{\bf 8}}_{{\bf 8},{\bf 8}}\bigr)_{\id,\id} &= \frac{1}{\qnum{3}+\qnum{5}} &
\bigl(F^{{\bf 10},{\overline{\bf 10}}}_{{\bf 10},{\bf 10}}\bigr)_{\id,\id} &= \frac{1}{(\qnum{3}-1)\qnum{5}} \ .
\end{align}

\section{The case $so(5)_1$}
\label{appso51}

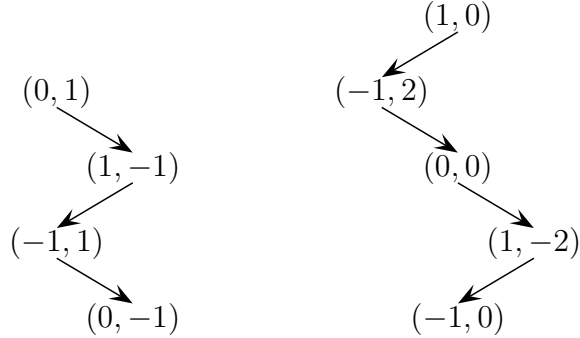
\begin{figure}[ht]
\begin{center}
\psset{unit=1mm,linewidth=.2mm,dimen=middle,arrowsize=4pt 4}
\begin{pspicture}(-5,-2)(15,32)
\rput(0,30){$(0,1)$}
\rput(10,20){$(1,-1)$}
\rput(0,10){$(-1,1)$}
\rput(10,0){$(0,-1)$}
\psline{->}(0,28)(10,22)
\psline{->}(10,18)(0,12)
\psline{->}(0,8)(10,2)
\end{pspicture}
\hspace{2 cm}
\begin{pspicture}(-5,-2)(25,42)
\rput(10,40){$(1,0)$}
\rput(0,30){$(-1,2)$}
\rput(10,20){$(0,0)$}
\rput(20,10){$(1,-2)$}
\rput(10,0){$(-1,0)$}
\psline{->}(10,38)(0,32)
\psline{->}(0,28)(10,22)
\psline{->}(10,18)(20,12)
\psline{->}(20,8)(10,2)
\end{pspicture}
\end{center}
\caption{The weights of the $so(5)_1$ representations ${\bf  4}$ and ${\bf 5}$.}
\label{so51ws}
\end{figure}

The Cartan matrix, its inverse and the quadratic form matrix for $so(5)$ read
\begin{align}
A &= \begin{pmatrix} 2 & -2 \\ -1 & 2 \end{pmatrix} & 
A^{-1} &= \frac{1}{2} \begin{pmatrix} 2 & 2 \\ 1 & 2 \end{pmatrix} & 
F_{\rm qf} &= \frac{1}{2} \begin{pmatrix} 2 & 1 \\ 1 & 1 \end{pmatrix} \ .
\end{align}
The second root is short, and we have $t_1=1$ and $t_2 = 2$, where
$t_i = \frac{2}{(\alpha_i,\alpha_i)}$. For $so(5)$, one has that
$\overline{\Lambda} = \Lambda$. The weight spaces of the representations ${\bf 4} = (0,1)$ and
${\bf 5} = (1,0)$, relevant for $so(5)_1$, are given in figure \ref{so51ws}.

The non-trivial fusion rules of $so(5)_1$ are
\begin{align}
{\bf 4} \times {\bf 4} &= \id + {\bf 5} &
{\bf 4} \times {\bf 5} &= {\bf 4} &
{\bf 5} \times {\bf 5} &= \id \ ,
\end{align}
namely the fusion rules of $su(2)_2$ or the Ising conformal field theory.

Continuing with the topological data, the quantum dimensions are given by
$d_{\id}=1$, $d_{\bf 4} = \qnum{5}_2-1\rightarrow \sqrt{2}$ and
$d_{\bf 5} = \qnum{4}+1\rightarrow 1$, where the numerical values are obtained by
setting $q=e^{2 \pi i /4}$.
The twist factors are given by
$\theta_{\bf 4} = q^{\frac{5}{4}}$ and $\theta_{\bf 5} = q^{2}$ and the
Frobenius-Schur indicators by 
${\rm fb}_\id = {\rm fb}_{\bf 5} = 1$ and ${\rm fb}_{\bf 4} = -1$.
Finally, the central charge is $\frac{5}{2}$.

\subsection{The $q$-CG coefficients for $so(5)_1$}.
In this section, we give the $q$-CG coefficients for $so(5)_1$.

\subsubsection{$\bf{4}\times{4}=\id + \bf{5}$}

\noindent 
\begin{tabular}{c|cccc}
& $(0,1)$ & $(1,-1)$ & $(-1,1)$ & $(0,-1)$ \\
\hline
$(0,1)$ & &
$\frac{q^{\frac{1}{8}}}{\sqrt{\qnum{2}_3}}$ &
$\frac{q^{\frac{1}{8}}}{\sqrt{\qnum{2}_3}}$ &
$\frac{1}{\qnum{2}_3}$ \\
$(1,-1)$ & $\frac{-q^{-\frac{1}{8}}}{\sqrt{\qnum{2}_3}}$ & &
$\frac{q^{\frac{1}{4}}}{\qnum{2}_3}$ &
$\frac{q^{\frac{1}{8}}}{\sqrt{\qnum{2}_3}}$\\
$(-1,1)$ & $\frac{-q^{-\frac{1}{8}}}{\sqrt{\qnum{2}_3}}$ &
$\frac{-q^{-\frac{1}{4}}}{\qnum{2}_3}$ & &
$\frac{q^{\frac{1}{8}}}{\sqrt{\qnum{2}_3}}$\\
$(0,-1)$ & $\frac{-1}{\qnum{2}_3}$ & 
$\frac{-q^{-\frac{1}{8}}}{\sqrt{\qnum{2}_3}}$ &
$\frac{-q^{-\frac{1}{8}}}{\sqrt{\qnum{2}_3}}$ & \\
\multicolumn{5}{c}{\phantom{$\sqrt{5}$}}
\\
\multicolumn{5}{c}{The $q$-CG coefficients for the $\bf{5}$ in ${\bf 4}\times{\bf 4}$.}
\end{tabular}
%
%\noindent
\begin{tabular}{c|cccc}
& $(0,1)$ & $(1,-1)$ & $(-1,1)$ & $(0,-1)$ \\
\hline
$(0,1)$ & & & & $\frac{q^{\frac{1}{2}}}{\sqrt{\qnum{5}_2-1}}$ \\
$(1,-1)$ & & & $\frac{-q^{\frac{1}{4}}}{\sqrt{\qnum{5}_2-1}}$ & \\
$(-1,1)$ & & $\frac{q^{-\frac{1}{4}}}{\sqrt{\qnum{5}_2-1}}$ & & \\
$(0,-1)$ & $\frac{-q^{-\frac{1}{2}}}{\sqrt{\qnum{5}_2-1}}$ & & & \\
\multicolumn{5}{c}{\phantom{$\sqrt{5}$}}
\\
\multicolumn{5}{c}{The $q$-CG coefficients for the $\id$ in ${\bf 4}\times{\bf 4}$.}
\end{tabular}

\subsubsection{$\bf{4}\times{5}=\bf{4}$}

\noindent
\begin{tabular}{c|ccccc}
& $(1,0)$ & $(-1,2)$ & $(0,0)$ & $(1,2)$ & $(-1,0)$ \\
\hline
$(0,1)$ & & & $\frac{q^{\frac{1}{2}}}{\sqrt{\qnum{5}_2}}$ &
$\frac{q^{\frac{3}{8}}\sqrt{\qnum{2}_2}}{\sqrt{\qnum{5}_2}}$ &
$\frac{q^{\frac{3}{8}}\sqrt{\qnum{2}_2}}{\sqrt{\qnum{5}_2}}$  \\
$(1,-1)$ & & $\frac{-q^{\frac{1}{8}}\sqrt{\qnum{2}_2}}{\sqrt{\qnum{5}_2}}$ &
$\frac{-1}{\sqrt{\qnum{5}_2}}$ & &
$\frac{q^{\frac{3}{8}}\sqrt{\qnum{2}_2}}{\sqrt{\qnum{5}_2}}$\\
$(-1,1)$ & $\frac{q^{\frac{-3}{8}}\sqrt{\qnum{2}_2}}{\sqrt{\qnum{5}_2}}$ & &
$\frac{-1}{\sqrt{\qnum{5}_2}}$ &
$\frac{-q^{\frac{-1}{8}}\sqrt{\qnum{2}_2}}{\sqrt{\qnum{5}_2}}$ & \\
$(0,-1)$ & $\frac{q^{\frac{-3}{8}}\sqrt{\qnum{2}_2}}{\sqrt{\qnum{5}_2}}$ &
$\frac{q^{\frac{-3}{8}}\sqrt{\qnum{2}_2}}{\sqrt{\qnum{5}_2}}$ &
$\frac{q^{\frac{-1}{2}}}{\sqrt{\qnum{5}_2}}$ & & \\
\multicolumn{5}{c}{\phantom{$5$}}
\\
\multicolumn{5}{c}{The $q$-CG coefficients for the $\bf{4}$ in ${\bf 4}\times{\bf 5}$.}
\end{tabular}\\

\noindent We note that the coefficients for $\bf{5}\times\bf{4}=\bf{4}$ are obtained by using the relation
\[\cg{\bf 4}{m_1}{\bf 5}{m_2}{\bf 4}{m_1+m_2}{q} = 
\cg{\bf 5}{m_2}{\bf 4}{m_1}{\bf 4}{m_1+m_2}{\frac{1}{q}}
\]
from section \ref{symandconv}.

\subsubsection{$\bf{5}\times{5}=\id$}

\noindent
\begin{tabular}{c|ccccc}
& $(1,0)$ & $(-1,2)$ & $(0,0)$ & $(1,-2)$ & $(-1,0)$ \\
\hline
$(1,0)$ & & & & & $\frac{q^{\frac{3}{4}}}{\sqrt{\qnum{4}+1}}$ \\
$(-1,2)$ & & & & $\frac{-q^{\frac{1}{4}}}{\sqrt{\qnum{4}+1}}$ & \\
$(0,0)$ & & & $\frac{1}{\sqrt{\qnum{4}+1}}$ & & \\
$(1,-2)$ & & $\frac{-q^{-\frac{1}{4}}}{\sqrt{\qnum{4}+1}}$ & & & \\
$(-1,0)$ & $\frac{q^{-\frac{3}{4}}}{\sqrt{\qnum{4}+1}}$ & & & & \\
\multicolumn{5}{c}{\phantom{$5$}}
\\
\multicolumn{5}{c}{The $q$-CG coefficients for the $\id$ in ${\bf 5}\times{\bf 5}$.}
\end{tabular}

\subsection{The $F$ and $R$-symbols for $so(5)_1$.}

The F-symbols which are not unity are give by
\begin{align}
\fs{\bf{4}}{\bf{4}}{\bf{5}}{\bf{5}} &=
\fs{\bf{4}}{\bf{5}}{\bf{5}}{\bf{4}} =
\fs{\bf{5}}{\bf{4}}{\bf{4}}{\bf{5}} =
\fs{\bf{5}}{\bf{5}}{\bf{4}}{\bf{4}} =
\frac{1}{\sqrt{\qnum{4}+1}} \rightarrow1\\
\fs{\bf{4}}{\bf{5}}{\bf{4}}{\bf{5}} &=
\fs{\bf{5}}{\bf{4}}{\bf{5}}{\bf{4}} =
-\frac{\qnum{3}_2}{\qnum{5}_2} \rightarrow-1\\
\fs{\bf{5}}{\bf{5}}{\bf{5}}{\bf{5}} &=
\frac{1}{\qnum{4}+1} \rightarrow 1 
\end{align}
In addition, we have
\begin{equation}
\fs{\bf{4}}{\bf{4}}{\bf{4}}{\bf{4}} = \begin{pmatrix}
\frac{-1}{\qnum{5}_2-1} & -\frac{\sqrt{\qnum{5}_2}}{\qnum{2}_2^2\sqrt{\qnum{4}+1}} \\
-\frac{\sqrt{\qnum{5}_2}}{\qnum{2}_2^2\sqrt{\qnum{4}+1}} & \frac{\qnum{3}_2}{(\qnum{2}_2)^{2}}
\end{pmatrix} \rightarrow
\begin{pmatrix} -1/\sqrt{2} & -1/\sqrt{2} \\ -1/\sqrt{2} & 1/\sqrt{2} \end{pmatrix}
\end{equation}

The values of the $R$-matrix are as follows
\begin{align}
R^{\id,a}_{a} &= R^{a,\id}_a = 1 \\
R^{{\bf 4},{\bf 4}}_{\bf 5} &= -q^{-\frac{1}{4}} & R^{{\bf 4},{\bf 4}}_{\id} &= -q^{-\frac{5}{4}} \\
R^{{\bf 4},{\bf 5}}_{\bf 4} &=  R^{{\bf 5},{\bf 4}}_{\bf 4} = q^{-1} \\
R^{{\bf 5},{\bf 5}}_{\id} &= q^{-2}
\end{align}

\section{The case $G_2$, $k=1$}
\label{appg21}

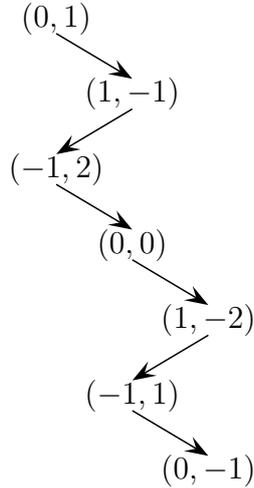
\begin{figure}[ht]
\begin{center}
\psset{unit=1mm,linewidth=.2mm,dimen=middle,arrowsize=4pt 4}
\begin{pspicture}(-5,-2)(25,62)
\rput(0,60){$(0,1)$}
\rput(10,50){$(1,-1)$}
\rput(0,40){$(-1,2)$}
\rput(10,30){$(0,0)$}
\rput(20,20){$(1,-2)$}
\rput(10,10){$(-1,1)$}
\rput(20,0){$(0,-1)$}
\psline{->}(0,58)(10,52)
\psline{->}(10,48)(0,42)
\psline{->}(0,38)(10,32)
\psline{->}(10,28)(20,22)
\psline{->}(20,18)(10,12)
\psline{->}(10,8)(20,2)
\end{pspicture}
\end{center}
\caption{The weights of the $G_2$ representation ${\bf 7}$.}
\label{g21ws}
\end{figure}

The Cartan matrix, its inverse and the quadratic form matrix for $G_2$ read
\begin{align}
A &= \begin{pmatrix} 2 & -3 \\ -1 & 2 \end{pmatrix} & 
A^{-1} &= \begin{pmatrix} 2 & 3 \\ 1 & 2 \end{pmatrix} & 
F_{\rm qf} &= \frac{1}{3} \begin{pmatrix} 6 & 3 \\ 3 & 2 \end{pmatrix} \ ,
\end{align}
thus, the second root is short, and we have $t_1=1$ and $t_2 = 3$, where
$t_i = \frac{2}{(\alpha_i,\alpha_i)}$. For $G_2$, one has that
$\overline{\Lambda} = \Lambda$. The weight spaces of the representation ${\bf 7} = (0,1)$
is given in figure \ref{g21ws}.

The only non-trivial fusion rule is ${\bf 7} \times {\bf 7} = \id \times {\bf 7}$, i.e. the fibonacci
fusion rule.  Furthermore,
$d_{\bf 7} = \qnum{11}_3 - \qnum{7}_3 + \qnum{3}_3 = \qnum{7}_3(\qnum{5}_3-\qnum{3}_3-1)
\rightarrow \phi$, where we used $q=e^{2 \pi i/5}$. The twist factor is $\theta_{\bf 7} = q^{2}$.
In addition, we find ${\rm fb}_\id = {\rm fb}_{\bf 7} = 1$, while the central charge is $\frac{14}{5}$.

\subsection{The $q$-CG coefficients for $G_2$, $k=1$}

\noindent
\hspace{-1cm}
\begin{tabular}{c|ccccccc}
& $(0,1)$ & $(1,-1)$ & $(-1,2)$ & $(0,0)$ & $(1,-2)$ &$(-1,1)$ & $(0,-1)$ \\
\hline
$(0,1)$ & & & & $q^{\frac{1}{2}}$ & $q^{\frac{5}{12}}\sqrt{\qnum{2}_3}$ & $q^{\frac{5}{12}}\sqrt{\qnum{2}_3}$ & $q^{\frac{1}{3}}$ \\ 
$(1,-1)$ & & & $-q^{\frac{1}{4}}\sqrt{\qnum{2}_3}$ & $-q^{\frac{1}{6}}$ & & $q^{\frac{1}{2}}$ & $q^{\frac{5}{12}}\sqrt{\qnum{2}_3}$ \\ 
$(-1,2)$ & & $q^{-\frac{1}{4}}\sqrt{\qnum{2}_3}$ & & $-q^{\frac{1}{6}}$ & $-1$ & & $q^{\frac{5}{12}}\sqrt{\qnum{2}_3}$ \\ 
$(0,0)$ & $-q^{-\frac{1}{2}}$ & $q^{-\frac{1}{6}}$ & $q^{-\frac{1}{6}}$ & $q^{-\frac{1}{6}}-q^{\frac{1}{6}}$ & $-q^{\frac{1}{6}}$ & $-q^{\frac{1}{6}}$ & $q^{\frac{1}{2}}$ \\ 
$(1,-2)$ & $-q^{-\frac{5}{12}}\sqrt{\qnum{2}_3}$ & & $1$ & $q^{-\frac{1}{6}}$ & & $-q^{\frac{1}{4}}\sqrt{\qnum{2}_3}$ & \\ 
$(-1,1)$ & $-q^{-\frac{5}{12}}\sqrt{\qnum{2}_3}$ & $-q^{-\frac{1}{2}}$ & & $q^{-\frac{1}{6}}$ & $q^{-\frac{1}{4}}\sqrt{\qnum{2}_3}$ & & \\ 
$(0,-1)$ & $-q^{-\frac{1}{3}}$ & $-q^{-\frac{5}{12}}\sqrt{\qnum{2}_3}$ & $-q^{-\frac{5}{12}}\sqrt{\qnum{2}_3}$ & $-q^{-\frac{1}{2}}$ & & & \\ 
\multicolumn{8}{c}{\phantom{$5$}}
\\
\multicolumn{8}{c}{The $q$-CG coefficients for the $\bf{7}$ in $\bf{7}\times\bf{7}$. Each
coefficient is to multiplied by $\frac{1}{\sqrt{\qnum{7}_3-1}}$.}
\end{tabular}\\ \\

\noindent 
\hspace{-1cm}
\begin{tabular}{c|ccccccc}
& $(0,1)$ & $(1,-1)$ & $(-1,2)$ & $(0,0)$ & $(1,-2)$ &$(-1,1)$ & $(0,-1)$ \\
\hline
$(0,1)$ & & & & & & & $q^{\frac{5}{6}}$ \\ 
$(1,-1)$ & & & & & & $-q^{\frac{2}{3}}$ & \\ 
$(-1,2)$ & & & & & $q^{\frac{1}{6}}$ & & \\ 
$(0,0)$ & & & & $-1$ & & & \\ 
$(1,-2)$ & & & $q^{-\frac{1}{6}}$ & & & & \\ 
$(-1,1)$ & & $-q^{-\frac{2}{3}}$ & & & & & \\ 
$(0,-1)$ & $q^{-\frac{5}{6}}$ & & & & & & \\ 
\multicolumn{8}{c}{\phantom{$5$}}
\\
\multicolumn{8}{c}{The $q$-CG coefficients for the $\id$ in $\bf{7}\times\bf{7}$.
%Each coefficient is to multiplied by 
Factor:
$\frac{1}{\sqrt{\qnum{11}_3-\qnum{7}_3+\qnum{3}_3}}$.}
\end{tabular}\\

\subsection{The $F$ and $R$-symbols}

The only F-symbols which are not equal to one are the following
(the numerical values are for $q=e^{2 \pi i/5}$)
\begin{equation}
\fs{\bf 7}{\bf 7}{\bf 7}{\bf 7} = \begin{pmatrix}
\frac{1}{\qnum{11}_3 - \qnum{7}_3 + \qnum{3}_3} &
-\frac{1}{\sqrt{\qnum{11}_3 - \qnum{7}_3 + \qnum{3}_3}} \\
-\frac{1}{\sqrt{\qnum{11}_3 - \qnum{7}_3 + \qnum{3}_3}} &
-\frac{\qnum{3}_3-2}{\qnum{5}_3-\qnum{3}_3}
\end{pmatrix}
\rightarrow
\begin{pmatrix}
1/\phi & -1/\sqrt{\phi} \\ -1/\sqrt{\phi} & -1/\phi
\end{pmatrix}
\end{equation}

In addition, we obtain the following $R$-symbols.
\begin{align}
R^{\id,a}_{a} &= R^{a,\id}_{a} = 1 &
R^{{\bf 7},{\bf 7}}_{\id} &= q^{-2} & 
R^{{\bf 7},{\bf 7}}_{\bf 7} &= -q^{-1}
\end{align}

\section{The case $su(2)_k$}
\label{appsu2k}

For completeness, we will give an explicit expression for the $q$-Clebsch-Gordan coefficients, as
well as the $F$ and $R$-symbols in the case of $su(2)_k$ (see, for instance, \cite{kr88}), using the same
basis conventions as we used throughout this paper. In particular, this formula yields $1$ for any
$F$-symbol with an identity on any of the outer lines. We will use Dynkin notation to denote the
particles, thus, the labels of the particles take the values $a=0,1,2,\ldots k$, in the case of $su(2)_k$.

The fusion rules can be written as
\begin{equation}
a \times b = \sideset{}{'}\sum_{c=|a-b|}^{\min(a+b,2k-a-b)} c \ ,
\end{equation}
where $c$ increases in steps of two. The quantum dimensions simply read
$d_a = \qnum{a+1}$; the twist factors are $\theta_a = q^{a(a+2)/4}$, 
while the Frobenius-Schur indicators are ${\rm fb}_a = (-1)^a$. 

To write the $q$-Clebsch-Gordan coefficients,
we define $\qnum{n}! = \qnum{n} \qnum{n-1} \cdots \qnum{1}$, and $\qnum{0}!=1$.
Furthermore, for $a\leq b+c$, $b\leq a+c$, $c\leq a+b$ and $a+b+c = 0 \bmod 2$, we define
\begin{equation}
\Delta (a,b,c) = \sqrt{\frac{\qnum{(a+b-c)/2}!\qnum{(a-b+c)/2}!\qnum{(-a+b+c)/2}!}{\qnum{(a+b+c+2)/2}!}}
\end{equation}
With  these conventions, the $q$-Clebsch-Gordan coefficients can, for instance, be written in the
following way  (see \cite{gkk90} for this particular form, and various others, as well as
\cite{stk91})
\begin{equation}
\begin{split}
&\qcg{a}{k}{b}{l}{c}{m} = 
q^{((a+b-c)(a+b+c+2)+2(al-bk))/16}\Delta(a,b,c)\times\\
&\sqrt{\qnum{(a-k)/2}!\qnum{(a+k)/2}!\qnum{(b-l)/2}!\qnum{(b+l)/2}!\qnum{(c-m)/2}!\qnum{(c+m)/2}!
\qnum{c+1}}\\
&\sideset{}{'}\sum_{n} \frac{(-1)^{n/2} q^{-n(a+b+c+2)/8}}{\qnum{n/2}!\qnum{(a-k-n)/2}!\qnum{(b+l-n)/2}!}\times\\
&\frac{1}{\qnum{(a+b-c-n)/2}!\qnum{(c-b+k+n)/2}!\qnum{(c-a-l+n)/2}!} \ ,
\end{split}
\end{equation}
where the sum over $n$ is over (non-negative) even integers, such that
$\max(0,-(c-b+k),-(c-a-l))\leq n \leq \min(a+b-c,a-k,b+l)$, in order that the arguments of the
$q$-factorials are non-negative integers.

One way of writing the $F$-symbols is as follows \cite{kr88}
\begin{equation}
\begin{split}
&\bigl(\fs{a}{b}{c}{d}\bigr)_{e,f} = 
(-1)^{(a+b+c+d)/2} \Delta(a,b,e)\Delta(c,d,e)\Delta(b,c,f)\Delta(a,d,f)
\sqrt{\qnum{e+1}}\sqrt{\qnum{f+1}}\\
&\sideset{}{'}\sum_{n}
\frac{(-1)^{n/2}\qnum{(n+2)/2}!}
{\qnum{(a+b+c+d-n)/2}!\qnum{(a+c+e+f-n)/2}!\qnum{(b+d+e+f-n}!}\\
&\times \frac{1}
{\qnum{(n-a-b-e)/2}!\qnum{(n-c-d-e)/2}!\qnum{(n-b-c-f)/2}!\qnum{(n-a-d-f)/2}!} \ ,
\end{split}
\end{equation}
where the sum over $n$ is over (non-negative) even integers, such that
$\max(a+b+e,c+d+e,b+c+f,a+d+f)\leq n \leq \min(a+b+c+d,a+c+e+f,b+d+e+f)$.

In addition, the $R$-symbols read as follow
\begin{equation}
R^{a,b}_{c} = (-1)^{(a+b-c)/2} q^{\frac{1}{8}(c(c+2)-a(a+2)-b(b+2))} \ .
\end{equation}

\section{The pentagon and hexagon equations}
\label{apppenthex}

We note that in the presence of fusion multiplicities, the vertices carry an additional
label. For completeness, we give the appropriate form of the pentagon and hexagon
equations here. The greek index $\alpha$, associated to a vertex with labels $(a,b,c)$,
runs over the values $1,2,\ldots,N_{a,b}^{c}$, where $n_{a,b}^{c}$ is the number of times
$c$ appears in the fusion product of $a\times b$. With this notation, we have the following
condition, for all possible $j_1,j_2,j_3,j_4,j,j_{12},j_{123},j_{34},j_{234}$ and
$\alpha,\beta,\gamma,\eta,\iota,\kappa$
\begin{equation}
\begin{split}
\sum_{j_{23},\delta,\epsilon,\zeta}
\bigl( F^{j_{1},j_{2},j_{3}}_{j_{123}}\bigr)_{(j_{12},\alpha,\beta;j_{23},\delta,\epsilon)}
\bigl( F^{j_{1},j_{23},j_{4}}_{j}\bigr)_{(j_{123},\epsilon,\gamma;j_{234},\zeta,\eta)}
\bigl( F^{j_{2},j_{3},j_{4}}_{j_{234}}\bigr)_{(j_{23},\delta,\zeta;j_{34},\iota,\kappa)} = \\
\sum_{\lambda}
\bigl( F^{j_{12},j_{3},j_{4}}_{j}\bigr)_{(j_{123},\beta,\gamma;j_{34},\iota,\lambda)}
\bigl( F^{j_{1},j_{2},j_{34}}_{j}\bigr)_{(j_{12},\alpha,\lambda;j_{234},\kappa,\eta)}
\end{split}
\end{equation}
To write down the hexagon equations, we will assume that twisting a vertex does not
change the internal label, or in other words, the fusion channel. With that restriction, the
hexagon equations read, for all $j_1,j_2,j_3,j,j_{12},j_{13}$ and
$\alpha,\beta,\gamma,\delta$
\begin{equation}
R^{j_{1},j_{2}}_{j_{12}}
\bigl( F^{j_{2},j_{1},j_{3}}_{j} \bigr)_{(j_{12},\alpha,\beta;j_{13},\gamma,\delta)}
R^{j_1,j_3}_{j_{13}} =
\sum_{j_{23},\epsilon,\zeta}
\bigl( F^{j_{2},j_{1},j_{3}}_{j} \bigr)_{(j_{12},\alpha,\beta;j_{23},\epsilon,\zeta)}
R^{j_{1},j_{23}}_{j}
\bigl( F^{j_{2},j_{3},j_{1}}_{j} \bigr)_{(j_{23},\epsilon,\zeta;j_{13},\gamma,\delta)} \ ,
\end{equation}
and a similar equation for $R^{-1}$.

%\bibliographystyle{unsrt}
%\bibliography{sixj}

\begin{thebibliography}{99}
\addtolength{\itemsep}{-3mm}

\bibitem{Kitaev03}
A.Yu.~Kitaev,
%{\it Fault-tolerant quantum computation by anyons},
Ann. Phys. {\bf 303}, 2 (2003).

\bibitem{Nayak08} C.~Nayak,  S.~H.~Simon, 
A.~Stern, M.~Freedman, S.~Das Sarma, Rev. Mod. Phys. \textbf{80}, 1083 (2008). 

\bibitem{hh07}
T.J.~Hagge, S.-M~Hong,
%{\it Some non-braided fusion categories of rank 3},
%arXiv:0704.0208[math.GT].
Comm. Cont. Math. {\bf 11}, 615 (2009).

\bibitem{Turaev}
V.G.~Turaev,
{\it Quantum Invariants of Knots and 3-Manifolds},
Watler de Gruyter, Berlin, New York (1994).

\bibitem{Kitaev06} 
A.~Kitaev,  Annals of 
Physics, \textbf{321}, 2 (2006).

\bibitem{Preskill-lectures}
J.~Preskill,
{\it Topological Quantum Computation}, Chapter 9 (2004),\\
\href{http://www.theory.caltech.edu/~preskill/ph219/topological.pdf}
{http://www.theory.caltech.edu/\raisebox{-2mm}{\~{ }}preskill/ph219/topological.pdf}

\bibitem{mslectures}
G.~Moore, N.~Seiberg,
{\it Lectures on RCFT},
in H.C.~Lee, ed.,
{\it Physics, Geometry and Topology},
p  263, Plenum Press, New York,  (1990).

\bibitem{Kassel95}
C.~Kassel,
{\it Quantum Groups},
Springer-Verlag, New-York, Berlin, Heidelberg (1995).

\bibitem{kac}
The program Kac is written by A.N.~Schellekens, see\\
\href{http://web.me.com/aschellekens/Site/Kac.html}
{http://web.me.com/aschellekens/Site/Kac.html}

\bibitem{bl81}
L.C.~Biedenharn, J.D.~Louck,
{\it The Racah-Wigner algebra in quantum theory},
Addison-Wesley, Advanced Book Program, (1981).


\bibitem{fms}
P.~Di~Francesco, P.~Mathieu, D. S\'en\'echal,
{\it Conformal field theory},
Springer, New York (1999).

\bibitem{fs}
J.~Fuchs, C.~Schweigert,
{\it Symmetries, lie algebras and representations: a graduate course for
physicists},
Cambridge monographs on mathematical physics,
Cambridge Univ. Press, (1997).

\bibitem{cftqgconnection}
J.~Fuchs,
{\it Affine Lie algebras and quantum groups},
Cambridge monographs on mathematical physics,
Cambridge University Press (1992);\\
L.~Alvarez-Gaum\'{e}, M.~Ruiz-Altaba, G.~Sierra,
{\it Quantum groups in tow-dimensional physics},
Cambridge monographs on mathematical physics,
Cambridge University Press (1996).

\bibitem{bonderson-thesis}
P.~Bonderson, Ph.D thesis, Caltech, (2007).


\bibitem{topintmr}
E.~Fradkin, C.~Nayak, A.~Tsvelik, F.~Wilczek,
%{\it A Chern-Simons effective field theory for the Pfaffian quantum Hall state},
Nucl. Phys. B {\bf 516}, 704 (1998);\\
A.~Stern, B.I.~Halperin,
%{\it Proposed Experiments to Probe the Non-Abelian $?=5/2$ Quantum Hall State},
Phys. Rev. Lett. {\bf 96}, 016802 (2006);\\
P.~Bonderson, A.~Kitaev, K.~Shtengel,
%{\it Detecting Non-Abelian Statistics in the $\nu=5/2$ Fractional Quantum Hall State},
Phys. Rev. Lett. {\bf 96}, 016803 (2006).

\bibitem{mr91}
G.~Moore, N.~Read,
%{\it Nonabelions in the fractional quantum Hall effect},
Nucl. Phys. B {\bf360}, 362 (1991).

\bibitem{topintrr}
P.~Bonderson, K.~Shtengel, J.K.~Slingerland,
%{\it Probing Non-Abelian Statistics with Quasiparticle Interferometry},
Phys. Rev. Lett. {\bf 97}, 016401 (2006);\\
S.B.~Chung, M.~Stone,
%{\it Proposal for reading out anyon qubits in non-abelian $\nu = 12/5$ quantum Hall state},
Phys. Rev. B {\bf 73}, 245311 (2006).

\bibitem{rr99}
N.~Read, E.~Rezayi,
%{\it Beyond paired quantum Hall states: Parafermions and incompressible states in the first excited Landau level},
Phys. Rev. B {\bf 59}, 8084 (1999).

\bibitem{geninterferometry}
B.J.~Overbosch, F.A.~Bais,
%{\it Inequivalent classes of interference experiments with non-abelian anyons},
Phys. Rev. A {\bf 64}, 062107 (2001);\\
P.~Bonderson, K.~Shtengel, J.K.~Slingerland,
%{\it Interferometry of non-Abelian Anyons},
Ann. Phys. {\bf 323}, 2709 (2008).

\bibitem{Willett09}  R.~L.~Willett,
 L.~N.~Pfeiffer,  K.~W.~West, Proceedings of the National Academy of Science, 106, 8853 (2009).

\bibitem{Dolev08} M.~Dolev, M.~Heiblum, 
V.~Umansky, A.~Stern, D.~Mahalu, Nature 452, 829 (2008).

\bibitem{Radu08} I.~P.~Radu, J.~B.~Miller,  C.~M.~Marcus, M.~A.~Kastner,   L.~N.~Pfeiffer,  K.~W.~West, Science, 320, 899 (2008).

\bibitem{sb01}
J.K.~Slingerland, F.A.~Bais,
%{\it Quantum groups and nonabelian braiding in quantum Hall systems},
Nucl. Phys. B {\bf 612}, 229 (2001).

\bibitem{nw96}
C.~Nayak, F.~Wilczek,
%{\it $2n$-quasihole states realize $2^{n-1}$-dimensional spinor braiding statistics in paired quantum Hall states},
Nucl. Phys. B {\bf 479}, 529 (1996).

\bibitem{as07}
E.~Ardonne, K.~Schoutens,
%{\it Wave functions for quantum registers},
Ann. Phys. {\bf 322}, 201 (2007).

\bibitem{as99}
E.~Ardonne, K.~Schoutens,
%{\it New class of non-abelian spin-singlet quantum Hall states},
Phys. Rev. Lett. {\bf 82}, 5096 (1999).

\bibitem{anyonicchains}
A.~Feiguin, S.~Trebst, A.W.W.~Ludwig, M.~Troyer, A.~Kitaev, Z.~Wang, M.H.~Freedman,
%{\it Interacting anyons in topological quantum liquids: The golden chain},
Phys. Rev. Lett. {\bf 98}, 160409 (2007);\\
S.~Trebst, E.~Ardonne, A.~Feiguin, D.A.~Huse, A.W.W.~Ludwig, M.~Troyer,
%{\it Collective states of interacting Fibonacci anyons},
Phys. Rev. Lett. {\bf 101}, 050401 (2008);\\
C.~Gils, E.~Ardonne, S.~Trebst, A.W.W.~Ludwig, M.~Troyer, Z.~Wang,
%{\it Topological stability of anyonic quantum spin chains and formation of new topological liquids},
%arXiv:0810.2277.
Phys. Rev. Lett. {\bf 103}, 070401 (2009).

\bibitem{disorderedanyons}
N.E.~Bonesteel, K.~Yang,
%{\it Infinite-Randomness Fixed Points for Chains of Non-Abelian Quasiparticles},
Phys. Rev. Lett. {\bf 99} 140405 (2007);\\
L.~Fidkowski, G.~Refael, N.E.~Bonesteel, J.E.~Moore, 
%{\it c-theorem violation for effective central charge of infinite-randomness fixed points},
Phys. Rev. B {\bf 78}, 224204 (2008);\\
L.~Fidkowski, G.~Refael, H.H.~Lin, P.~Titum,
%{\it Permutation Symmetric Critical Phases in Disordered Non-Abelian Anyonic Chains},
%arXiv:0812.3158.
Phys. Rev. B {\bf 79}, 155120 (2009).


\bibitem{ttw08}
S.~Trebst, M.~Troyer, Z.~Wang, A.W.W.~Ludwig,
%{\it A short introduction to Fibonacci anyon models},
Prog. Theor. Phys. Supp. {\bf 176}, 384 (2008).

\bibitem{bs-pentagon}
P. Bonderson, J.K.~Slingerland, in preparation.

\bibitem{kr88}
A.N.~Kirillov, N.Y.~Reshetikhin,
{\it Representations of the algebra $U_q(sl(2))$, $q$-orthogonal polynomials and invariants of links},
in V.G.~Kac, ed., {\it Infinite dimensional Lie algebras and groups, Proceedings of the conference held at CIRM, Luminy, Marseille}, p. 285, World Scientific, Singapore (1988).  

\bibitem{gkk90}
V.A.~Groza, I.I.~Kachurik, A.U.~Klimyk,
%{\it On Clebsch-Gordan coefficients and matrix elements of representations of the quantum algebra $U_q(su_2)$}
J. Math. Phys. {\bf 31}, 2769 (1990).


\bibitem{stk91}
Y.F.~Smirnov, V.N.~Tolstoi, Y.I.~Kharitonov,
%{METHOD OF PROJECTION OPERATORS AND THE Q ANALOG OF THE QUANTUM-THEORY OF ANGULAR-MOMENTUM - CLEBSCH-GORDAN-COEFFICIENTS AND IRREDUCIBLE TENSOR-OPERATORS},
Sov. J. Nucl. Phys. {\bf 53}, 593 (1991).



\end{thebibliography}

\end{document}